\documentclass[a4paper,10pt,oneside]{article}

\usepackage{t1enc}
\usepackage[latin1]{inputenc}
\usepackage{graphics}
\usepackage{fix-cm}

\usepackage{graphicx}
\usepackage{amsfonts}
\usepackage{eufrak}
\usepackage[all]{xypic}
\usepackage{amssymb}
\usepackage{amsmath}
\usepackage{amsthm}
\usepackage{mathrsfs}
\usepackage{makeidx}

\theoremstyle{plain}
\newtheorem{thm}{Theorem}[section]
\newtheorem{lemma}[thm]{Lemma}
\newtheorem{prop}[thm]{Proposition}
\newtheorem{cor}[thm]{Corollary}

\theoremstyle{definition}
\newtheorem{df}[thm]{Definition}
\newtheorem{rem}[thm]{Remark}
\newtheorem{ex}[thm]{Example}

\newtheorem{open}[thm]{Open Question}

\renewcommand{\a}{\mathcal{A}}
\renewcommand{\b}{\mathcal{B}}

\renewcommand{\d}{\mathcal{D}}

\providecommand{\gm}{\Gamma}
\providecommand{\h}{\mathcal{H}}
\providecommand{\cl}{\mathcal{C}}

\title{\bf Crossed products by Hecke pairs II: $C^*$-completions}
\author{Rui Palma}
\date{}

\makeindex

\begin{document}
\maketitle

\begin{abstract}
 In this second article on crossed products by ``actions'' of Hecke pairs we study their different $C^*$-completions, namely we show how reduced and full $C^*$-crossed products can be defined. We also establish that our construction coincides with that of Laca, Larsen and Neshveyev \cite{phase} whenever they are both definable. As an application of our theory, we prove a Stone-von Neumann theorem for Hecke pairs which encompasses the work of an Huef, Kaliszewski and Raeburn \cite{cov} and we lay down the foundations for obtaining a form of Katayama duality with respect to the Echterhoff-Quigg crossed product \cite{full}, \cite{cov}. \\
\end{abstract}

\tableofcontents

{\renewcommand{\thefootnote}{}
\footnotetext{\emph{Date:} \today}}

\newpage

\section*{Introduction}
\addcontentsline{toc}{section}{Introduction}

The present work is a continuation of \cite{palma crossed 11}, being the last of two articles whose goal is the development of a theory of crossed products by Hecke pairs with a view towards applications in non-abelian $C^*$-duality.

A \emph{Hecke pair} $(G, \gm)$ consists of a group $G$ and a subgroup $\gm \subseteq G$ for which every double coset $\gm g \gm$ is the union of finitely many left cosets. In this case $\gm$ is also said to be a \emph{Hecke subgroup} of $G$. Examples of Hecke subgroups include finite subgroups, finite-index subgroups and normal subgroups. It is in fact many times insightful to think of this definition as a generalization of the notion of normality of a subgroup.

Given a Hecke pair $(G,\gm)$ the \emph{Hecke algebra} $\h(G, \gm)$ is a $^*$-algebra of functions over the set of double cosets $\gm \backslash G / \gm$, with a suitable convolution product and involution. It generalizes the definition of the group algebra $\mathbb{C}(G / \gm)$ of the quotient group when $\gm$ is a normal subgroup.

A crossed product of an algebra $A$ by a Hecke pair $(G,\gm)$ should be thought of as a crossed product (in the usual sense) of $A$ by an ``action'' of $G / \gm$, even though $G / \gm$ is not necessarily a group. The term ``crossed product by a Hecke pair'' was first used by Tzanev \cite{tzanev talk} in order to give another perspective on the work of Connes and Marcolli \cite{connes marcolli}. This point of view was later formalized by Laca, Larsen and Neshveyev in \cite{phase}, where they defined a $C^*$-algebra which can be interpreted as a reduced $C^*$-crossed product of a commutative $C^*$-algebra by a Hecke pair.

In a different approach, we developed in \cite{palma crossed 11} a $^*$-algebraic theory of crossed products by Hecke pairs, which makes sense when the algebra $A$ is certain algebra of sections of a Fell bundle over a discrete groupoid. The main ingredients behind our construction were the following: we started with a Hecke pair $(G, \gm)$, a Fell bundle $\a$ over a discrete groupoid $X$ and an action $\alpha$ of $G$ on $\a$  satisfying some ``nice'' properties. From this we could naturally give the space $\a / \gm$ of $\gm$-orbits of $\a$ a Fell bundle structure over the orbit space $X / \gm$, which under our assumptions on the action $\alpha$ is in fact a groupoid. We then defined a $^*$-algebra $C_c(\a / \gm) \times^{alg} G / \gm$ which we justifiably called the \emph{$^*$-algebraic crossed product} of $C_c(\a / \gm)$ by the Hecke pair $(G, \gm)$.

In this second article on the subject of crossed products by Hecke pairs we will study the different $C^*$-completions of these $^*$-algebraic crossed products, with special emphasis on the reduced case which is technically more challenging to define, and explore some connections with non-abelian $C^*$-duality.

Reduced $C^*$-crossed products by groups are defined via the so-called regular representations. We will introduce a notion of a regular representation in the Hecke pair case by using the regular representation of the Hecke algebra on $\ell^2(G / \gm)$. The main novelty here is that we will have to start with a representation of a certain direct limit of algebras of the form $C_c(\a / H)$, where $H$ is a finite intersection of conjugates of the Hecke subgroup $\gm$. In case $\gm$ is normal, this direct limit is simply $C_c(\a / \gm)$ itself and we recover the usual notion of a covariant representation. This crucial direct limit is defined in Section \ref{direc limit section} and regular representations are studied in Section \ref{regular representations chapter}.

 From regular representations it is then possible to define \emph{reduced} $C^*$-crossed products. Since the algebra $C_c(\a / \gm)$ admits several $C^*$-completions there are several reduced $C^*$-crossed products that one can form, as for example $C^*_r(\a / \gm) \times_{\alpha, r} G / \gm$ and $C^*(\a) \times_{\alpha, r} G / \gm$, each of these thought of as the reduced $C^*$-crossed product of $C^*_r(\a / \gm)$, respectively $C^*(\a/\gm)$, by the Hecke pair $(G, \gm)$. These reduced $C^*$-crossed products have always a faithful conditional expectation onto $C^*_r(\a / \gm)$ (respectively, $C^*(\a / \gm)$), a property that determines the reduced crossed product uniquely, just like in the case of groups. This is achieved in Section \ref{reduced crossed products chapter}.

Our construction of reduced $C^*$-crossed products by Hecke pairs is different from that of Laca, Larsen and Neshveyev in \cite{phase}, being more particular in some sense (since we treat only discrete spectrum), but also more general (since it makes sense for certain non-commutative $C^*$-algebras). What we are going to show is that both constructions agree whenever they are both definable. This will be done in Section \ref{comparison with larsen laca neshveyev section}.

Complementing the reduced setting, one would like to form different \emph{full} $C^*$-crossed products, as for example $C^*_r(\a / \gm) \times_{\alpha} G / \gm$ and $C^*(\a / \gm) \times_{\alpha} G / \gm$, but in general their existence is not assured. They will always exist, however, if the Hecke algebra is a $BG^*$-algebra, which is a property that is satisfied by several classes of Hecke pairs, including most of those studied in the literature for which a full Hecke $C^*$-algebra is known to exist (see \cite{palma}). Full $C^*$-crossed products are studied in detail in Section \ref{other completions chapter}.

This theory of crossed products by Hecke pairs is intended for applications in non-abelian duality theory. We develop completely a Stone-von Neumann type theorem for Hecke pairs which encompasses the work of an Huef, Kaliszewski and Raeburn \cite{cov}, and we envisage for future work a form of Katayama duality with respect to Echterhoff-Quigg's   ``crossed product'' \cite{full}.

The Stone-von Neumann theorem, in the language of crossed products by groups, states that for the action of translation of $G$ on $C_0(G)$ we have
\begin{align*}
 C_0(G) \times G \;\cong\; C_0(G) \times_r G \;\cong \; \mathcal{K}(\ell^2(G))\,. 
\end{align*}
In \cite{cov} an Huef, Kaliszweski and Raeburn introduced the notion of \emph{covariant pairs} of representations of $C_0(G / \gm)$ and $\h(G, \gm)$, for a Hecke pair $(G, \gm)$, and proved that all covariant pairs are amplifications of a certain ``regular'' covariant pair. Their result was proven without using or defining crossed products, and can also be thought of as a Stone-von Neumann theorem for Hecke pairs. Using our construction we express their result in the language of crossed products (Section \ref{Stone von Neumann chapter}). We will show that the full crossed product $C_0(G / \gm) \times G /\gm$ always exists and one has
\begin{align*}
 C_0(G / \gm) \times G/ \gm \;\cong \; C_0(G / \gm) \times_r G/ \gm \;\cong \; \mathcal{K}(\ell^2(G / \gm))\,.
\end{align*}
Moreover, our notion of a covariant representation coincides with the notion of a covariant pair of \cite{cov}, and an Huef, Kaliszewski and Raeburn's result follows as a direct corollary of the above isomorphisms.

Our construction was very much influenced and developed with the wish of obtaining a form of Katayama duality for homogeneous spaces (those arising from Hecke pairs). Even though this has been left for future work, we shall nevertheless explain in Section \ref{towards katayama duality chapter} what we have in mind and how our set up is suitable for tackling this problem.

Katayama's duality theorem \cite{kat} is an analogue for coactions of the duality theorem of Imai and Takai. One version of it states the following: given a coaction $\delta$ of a group $G$ on a $C^*$-algebra $A$ and denoting by $A \times_{\delta} G$ the corresponding crossed product, we have a canonical isomorphism $A \times_{\delta} G \times_{\widehat{\delta}, \omega} G \cong A \otimes \mathcal{K}(\ell^2(G))$, for some crossed product by the dual action $\widehat{\delta}$ of $G$. We would like to extend this result to homogeneous spaces coming from Hecke pairs. In spirit we hope to obtain an isomorphism of the type:
\begin{align*}
A \times_{\delta} G / \gm \times_{\widehat{\delta}, \omega} G / \gm \cong A \otimes \mathcal{K}(\ell^2(G / \gm))\,.
\end{align*}
The $C^*$-algebra $A \times_{\delta} G / \gm$ should be a crossed product by a coaction of the homogeneous space $G / \gm$, while the second crossed product should be by the ``dual action'' of the Hecke pair $(G, \gm)$ in our sense. It does not make sense in general for a homogeneous space to coact on a $C^*$-algebra, but it is many times possible to define $C^*$-algebras which can be thought of as crossed products by coactions of homogeneous spaces (\cite{echt kal rae}, \cite{full}).

It is our point of view that $A \times_{\delta} G / \gm$ should be a certain $C^*$-completion of the $^*$-algebra $C_c(\a \times G / \gm)$ defined by Echterhoff and Quigg \cite{full}, which we dub the Echterhoff and Quigg's crossed product (a terminology used in \cite{cov} for $C^*(\a \times G / \gm)$ in case of a maximal coaction). We explain in Section \ref{towards katayama duality chapter} how our set up for defining crossed products by Hecke pairs is suitable for achieving such a Katayama duality result for Echterhoff and Quigg's crossed product, and can therefore bring insight into the emerging theory of crossed products by coactions of homogeneous spaces.

The present work is based on the author's Ph.D. thesis \cite{palmathesis} written at the University of Oslo. The author would like to thank his advisor Nadia Larsen for the very helpful discussions, suggestions and comments during the elaboration of this work. A word of appreciation goes also to John Quigg, Dana Williams and Erik B\'edos for some very helpful comments.\\

\section{Preliminaries}

In this section we set up the conventions, notation, and background results which will be used throughout this work. References where the reader can find more details will be indicated.

Since this article is a continuation of \cite{palma crossed 11} we will follow all of its conventions and notation, most of the times without reference. \\

\subsection{Hecke algebras}

\begin{df}
 Let $G$ be a group and $\gm$ a subgroup. The pair $(G , \gm)$\index{G2  gm@$(G, \gm)$} is called a \emph{Hecke pair} if every double coset $\gm g\gm$ is the union of finitely many right (and left) cosets. In this case, $\gm$ is also called a \emph{Hecke subgroup} of $G$.\\
\end{df}

Given a Hecke pair $(G, \gm)$ we will denote by $L$\index{L(g)@$L(g)$} and $R$\index{R(g)@$R(g)$}, respectively, the left and right coset counting functions, i.e.
\begin{align}
 L(g):= |\gm g \gm / \gm| = [\gm : \gm^g] < \infty\\
 R(g) :=|\gm \backslash \gm g \gm| = [\gm : \gm^{g^{-1}}] < \infty\,.
\end{align}
We recall that $L$ and $R$ are $\gm$-biinvariant functions which satisfy $L(g) = R(g^{-1})$ for all $g \in G$. Moreover, the function $\Delta: G \to \mathbb{Q^+}$\index{D1elta@$\Delta(g)$} given by
\begin{align}
\label{def modular function Hecke algebra}
 \Delta(g) := \frac{L(g)}{R(g)}\,,
\end{align}
is a group homomorphism, usually called the \emph{modular function} of $(G, \gm)$.\\

\begin{df}
 Given a Hecke pair $(G, \gm)$, the \emph{Hecke algebra} $\h(G, \gm)$\index{H(G,gm)@$\h(G, \gm)$} is the $^*$-algebra of finitely supported $\mathbb{C}$-valued functions on the double coset space $\gm \backslash G / \gm$ with the product and involution defined by
\begin{align}
 (f_1*f_2)(\gm g \gm) & := \sum_{h\gm \in G / \gm} f_1(\gm h \gm)f_2(\gm h^{-1}g\gm)\,,\\
f^*(\gm g\gm) & := \Delta(g^{-1}) \overline{f(\gm g^{-1} \gm)}\,.\\ \notag
\end{align}

\end{df}

The Hecke algebra has a natural basis, as a vector space, given by the characteristic functions of double cosets. As in \cite{palma crossed 11} we identify a characteristic function of a double coset $1_{\gm g \gm}$ with the double coset $\gm g \gm$ itself.

As it is known, group algebras have two canonical $C^*$-completions, the reduced group $C^*$-algebra $C^*_r(G)$ and the full group $C^*$-algebra $C^*(G)$. For Hecke algebras the situation becomes more complicated, there being essentially four canonical $C^*$-completions. We will briefly review these completions in this subsection, but first we need to recall the definitions and basic facts about regular representations of Hecke algebras and $L^1$-norms.\\

\begin{df}
 Let $(G, \gm)$ be a Hecke pair. The mapping $\rho:\h(G, \gm) \to B(\ell^2(G / \gm))$\index{rho@$\rho$} defined, for $f \in \h(G, \gm)$, $\xi \in \ell^2(G / \gm)$ and $g \gm \in G / \gm$, by
\begin{align}
\label{right regular rep of Hecke algebra main expression}
 (\rho(f)\xi)\,(g \gm):= \sum_{[h] \in G / \gm} \Delta(h)^{\frac{1}{2}} f(\gm h \gm) \xi(gh\gm)\,,
\end{align}
is called the \emph{right regular representation} of $\h(G, \gm)$.\\
\end{df}

It can be checked that $\rho$ does define a $^*$-representation of $\h(G, \gm)$. For the canonical vectors $\delta_{r \gm} \in \ell^2(G / \gm)$, expression (\ref{right regular rep of Hecke algebra main expression}) becomes:
\begin{align}
\label{right regular rep of Hecke algebra second expression}
 \rho(f)\delta_{r \gm} = \sum_{[g] \in G / \gm} \Delta(g^{-1}r)^{\frac{1}{2}} f(\gm g^{-1}r\gm) \delta_{g \gm}\,,
\end{align}
and furthermore for $f$ of the form $f:=\gm d \gm$ we obtain:
\begin{align}
 \rho(\gm d \gm)\delta_{r \gm} = \sum_{t\gm \subseteq \gm d^{-1} \gm} \Delta(d)^{\frac{1}{2}}  \delta_{rt \gm} \;=\; \Delta(d)^{\frac{1}{2}} \delta_{r\gm d^{-1} \gm}\,.
\end{align}
 
It can be easily checked, applying (\ref{right regular rep of Hecke algebra second expression}) to the vector $\delta_{\gm}$ for example, that $\rho$ always defines a faithful $^*$-representation.

One could in a similar fashion define a left regular representation of $\h(G, \gm)$, but in this work, however, it is the right regular representation the one that will play a central role.

We now recall the definition of the $L^1$-norm in a Hecke algebra (from \cite{schl}):\\

\begin{df}
 The \emph{$L^1$-norm} on $\h(G, \gm)$, denoted $\| \cdot \|_{L^1}$, is given by
\begin{align}
 \| f \|_{L^1} := \sum_{\gm g \gm \in \gm \backslash G / \gm} |f(\gm g \gm)|\, L(g) = \sum_{g\gm \in G / \gm} |f(\gm g \gm)|\,.
\end{align}
We will denote by $L^1(G, \gm)$\index{L1 G gm@$L^1(G, \gm)$} the completion of $\h(G, \gm)$ under this norm, which is a Banach $^*$-algebra.\\
\end{df}

There are several canonical $C^*$-completions of $\h(G, \gm)$ (\cite{schl}, \cite{tzanev}) These are:

\begin{itemize}
 \item $C^*_r(G, \gm)$\index{C(G gm)r@$C^*_r(G, \gm)$} - Called the \emph{reduced Hecke $C^*$-algebra}, it is the completion of $\h(G, \gm)$ under the $C^*$-norm arising from the right regular representation.
 \item $pC^*(\overline{G})p$ - The corner of the full group $C^*$-algebra $C^*(\overline{G})$ of the Schlichting completion $(\overline{G}, \overline{\gm})$ of the pair $(G, \gm)$, by the projection $p := 1_{\overline{\gm}}$. We will not describe this construction here since it is well documented in the literature (see \cite{schl}, for example) and because we will not make use of this $C^*$-completion in this work.
\item $C^*(L^1(G, \gm))$ - The enveloping $C^*$-algebra of $L^1(G, \gm)$.
\item $C^*(G, \gm)$\index{C(G gm)@$C^*(G, \gm)$} - The enveloping $C^*$-algebra (if it exists!) of $\h(G, \gm)$. When it exists, it is usually called the \emph{full Hecke $C^*$-algebra}.
\end{itemize}

The various $C^*$-completions of $\h(G, \gm)$ are related in the following way, through canonical surjective maps:
\begin{equation*}
 C^*(G,\gm)  \dashrightarrow C^*(L^1(G, \gm)) \longrightarrow pC^*(\overline{G})p \longrightarrow C^*_r(G, \gm)\,.
\end{equation*}

As was pointed out by Hall in \cite[Proposition 2.21]{hall}, the full Hecke $C^*$-algebra $C^*(G, \gm)$ does not have to exist in general, with the Hecke algebra of the pair $(SL_2(\mathbb{Q}_p), SL_2(\mathbb{Z}_p))$ being one such example, where $p$ is a prime number and $\mathbb{Q}_p$, $\mathbb{Z}_p$ denote respectively the field of $p$-adic numbers and the ring of $p$-adic integers.\\

\subsection{Cross sectional algebras of Fell bundles}

Let $X$ be a discrete groupoid. As in \cite{palma crossed 11} we denote by $X^0$ the unit space of $X$ and by $\mathbf{s}$\index{source@$\mathbf{s}(x)$} and $\mathbf{r}$\index{range@$\mathbf{r}(x)$} the source and range functions, respectively.

Given a Fell bundle $\a$ over a groupoid $X$ we will denote by $\a^0$\index{A0@$\a^0$} the restricted bundle $\a|_{X^0}$ over the unit space $X^0$. Since the fibers of $\a$ over $X^0$ are $C^*$-algebras, $\a^0$ is a $C^*$-bundle over $X^0$.

In this subsection we will briefly recall how the full and the reduced cross sectional algebras of a Fell bundle $\a$ over a groupoid $X$ are defined.

The \emph{full cross sectional algebra} of $\a$, denoted $C^*(\a)$, is defined as the enveloping $C^*$-algebra of $C_c(\a)$, which is known to exist (see for example \cite[Proposition 2.1]{full}).  If the groupoid $X$ is just a set, in which case $\a$ is a $C^*$-bundle, we will use the notation $C_0(\a)$ instead of $C^*(\a)$.

We now recall, from \cite{kum}, how the \emph{reduced cross sectional algebra} $C_r^*(\a)$ is defined. We see $C_c(\a)$ as a pre-Hilbert $C_0(\a^0)$-module, where the inner product is defined by
\begin{align*}
 \langle f_1, f_2 \rangle_{C_c(\a^0)} := (f_1^*\cdot f_2)|_{X^0}\,, \qquad\qquad f_1, f_2 \in C_c(\a)\,.
\end{align*}
Its completion is a full Hilbert $C_0(\a^0)$-module, which we denote by $L^2(\a)$. Now, the algebra $C_c(\a)$ acts on itself by left multiplication, and moreover this action is continuous with respect to the norm induced by the inner product above, hence we get an injective $^*$-homomorphism
\begin{align}
\label{homomorphism of CcB into LL2B}
 C_c(\a) \to \mathcal{L}(L^2(\a))\,.
\end{align}
The reduced cross sectional algebra $C^*_r(\a)$ is defined as the completion of $C_c(\a)$ with respect to the operator norm in $\mathcal{L}(L^2(\a))$, and in this way we get a right-Hilbert bimodule $_{C_r^*(\a)} L^2(\a) _{C_0(\a^0)}$.

Since $C^*_r(\a)$ is a completion of $C_c(\a)$ we immediately get a canonical map $\Lambda:C^*(\a) \to C^*_r(\a)$. Also, the $^*$-homomorphism above in (\ref{homomorphism of CcB into LL2B}) always completes to a $^*$-homomorphism $C^*(\a) \to \mathcal{L}(L^2(\a))$, and therefore gives rise to a right-Hilbert bimodule $_{C^*(\a)} L^2(\a) _{C_0(\a^0)}$. The image of $C^*(\a)$ on $\mathcal{L}(L^2(\a))$ is then isomorphic to $C_r^*(\b)$, or in other words, the kernel of the map $C^*(\a) \to \mathcal{L}(L^2(\a))$ is the same as the kernel of the canonical map $\Lambda : C^*(\a) \to C_r^*(\a)$.\\

\subsection{$^*$-Algebraic crossed products by Hecke pairs}

In \cite{palma crossed 11} we layed down the foundations for defining certain $^*$-algebraic crossed products by Hecke pairs. We will now just briefly recall the main assumptions behind our construction and the main concepts involved. The reader is referred to \cite{palma crossed 11} for the details. We will also provide a result concerning unitary equivalence which is going to be used later in this article.

The main ingredients in \cite{palma crossed 11} are the following: a Hecke pair $(G, \gm)$; a Fell bundle $\a$ over a discrete groupoid $X$; an action $\alpha$ of $G$ on $\a$ which is assumed to be $\gm$-good and to satisfy the $\gm$-intersection property (see \cite[Standing Assumnption 3.1]{palma crossed 11}).

From this assumptions we can give the orbit space $X / \gm$ a groupoid structure, and the orbit space $\a / \gm$ a Fell bundle structure over $X / \gm$ (see \cite[Section 2.1]{palma crossed 11}). Moreover, the action $\alpha$ of $G$ on $\a$ gives naturally rise to an action $\overline{\alpha}$ of $G$ on the algebra of sections $C_c(\a)$, which extends to the multiplier algebra $M(C_c(\a))$. We then define the \emph{$^*$-algebraic crossed product} of $C_c(\a / \gm)$ by the Hecke pair $(G, \gm)$ as the $^*$-algebra of functions $f:G / \gm \to M(C_c(\a))$ such that
\begin{itemize}
 \item[i)] $f(g \gm) \in C_c(\a / \gm^g)$, and
 \item[ii)] $f(\gamma g \gm) = \overline{\alpha}_{\gamma}(f(g\gm))$, for any $\gamma \in \gm$,
\end{itemize}
with the product and involution operations defined by
\begin{align*}
  (f_1 * f_2)(g \gm) & := \sum_{[h] \in G/ \gm} f_1(h \gm)\, \overline{\alpha}_{h}(f_2(h^{-1} g\gm))\,,\\
 (f^*)\,(g \gm) & := \Delta(g^{-1})\overline{\alpha}_{g}(f(g^{-1} \gm))^*\,.
\end{align*}
We denote this $^*$-algebraic crossed product by $C_c(\a / \gm) \times_{\alpha}^{alg} G / \gm$, and the reader is referred to \cite[Section 3.1]{palma crossed 11} for the details.

The algebra $C_c(\a /\gm)$ embeds canonically in $C_c(\a / \gm) \times_{\alpha}^{alg} G / \gm$. Moreover, the Hecke algebra $\h(G, \gm)$ and the algebra $C_c(X^0/ \gm)$ both embed canonically in the multiplier algebra $M(C_c(\a / \gm) \times_{\alpha}^{alg} G / \gm)$. As we showed in \cite[Theorem 3.14]{palma crossed 11}, the crossed product $C_c(\a / \gm) \times_{\alpha}^{alg} G / \gm$ is spanned by elements of the form
\begin{align}
\label{canonical spanning set of elements}
 [a]_{x\gm} * \gm g \gm * 1_{\mathbf{s}(x)g\gm}\,,
\end{align}
where $[a]_{x\gm} \in C_c(\a / \gm)$, $\gm g \gm \in \h(G, \gm)$ and $1_{\mathbf{s}(x)g\gm} \in C_c(X^0 / \gm)$.

Looking at how products of the form $\gm g \gm * [a]_{x\gm} * \gm s \gm$ decompose as a sum of the canonical elements (\ref{canonical spanning set of elements}), we can then define \emph{covariant pre-$^*$-representations}, which are pairs $(\pi, \mu)$ consisting of a nondegenerate $^*$-representation $\pi$ of $C_c(\a / \gm)$ on a Hilbert space $\mathscr{H}$ and a unital pre-$^*$-representation $\mu$ of $\mathcal{H}(G, \gm)$ on the dense inner product space $\pi(C_c(\a / \gm)) \mathscr{H}$, that preserve the mentioned product decomposition (see \cite[Definition 3.24]{palma crossed 11}).

As we showed in \cite[Theorem 3.40]{palma crossed 11} there is a bijective correspondence between covariant pre-$^*$-representations and nondegenerate $^*$-representations of $C_c(\a / \gm) \times_{\alpha}^{alg} G / \gm$, given by $(\pi, \mu) \mapsto \pi \times \mu$, where $\pi \times \mu$ is the \emph{integrated form} of $(\pi, \mu)$ as defined in \cite[Definition 3.33]{palma crossed 11}.

We will now show that the bijective correspondence between covariant pre-$^*$-representations and nondegenerate $^*$-representations of the crossed product behaves as expected regarding unitary equivalence. First however we make the following remark/definition:\\

\begin{rem}
Let $(\pi, \mu)$ be a covariant pre-$^*$-representation on a Hilbert space $\mathscr{H}$. If $\mathscr{H}_0$ is another Hilbert space and $U:\mathscr{H} \to \mathscr{H}_0$ is a unitary, then it is easily seen that $(U \pi U^*, U \mu U^*)$ is also a covariant pre-$^*$-representation. We will henceforward say that two covariant pre-$^*$-representations $(\pi_1, \mu_1)$ and $(\pi_2, \mu_2)$ are \emph{unitarily equivalent} if there exists a unitary $U$ between the underlying Hilbert spaces such that $(\pi_1, \mu_1)=(U\pi_2 U^*, U \mu_2U^*)$.\\
\end{rem}

\begin{prop}
\label{unitary equivalence and bijection}
Suppose that $(\pi_1, \mu_1)$ and $(\pi_2, \mu_2)$ are two covariant pre-$^*$-representations. Then $(\pi_1, \mu_1)$ is unitarily equivalent to $(\pi_2, \mu_2)$ if and only if $\pi_1 \times \mu_1$ is unitarily equivalent to $\pi_2 \times \mu_2$.\\
\end{prop}

{\bf \emph{Proof:}} $(\Longrightarrow)$ This direction is straightforward from the definition of the integrated form and from the following computation, where $U$ is a unitary which establishes an equivalence between $(\pi_1, \mu_1)$ and $(\pi_2, \mu_2)$ :
\begin{eqnarray*}
& & [U(\pi_1 \times \mu_1)U^*]([a]_{x\gm} * \gm g \gm * 1_{\mathbf{s}(x)g\gm}) \;\; = \;\;\\
 & = & U\pi_1([a]_{x\gm}) \mu_1(\gm g \gm )\widetilde{\pi_1}( 1_{\mathbf{s}(x)g\gm})U^*\\
& = & U\pi_1([a]_{x\gm})U^* U \mu_1(\gm g \gm )U^* U\widetilde{\pi_1}( 1_{\mathbf{s}(x)g\gm})U^*\\
& = & [U\pi_1 U^* \times U \mu_1 U^*]([a]_{x\gm} * \gm g \gm * 1_{\mathbf{s}(x)g\gm})\\
& = & [\pi_2 \times  \mu_2 ]([a]_{x\gm} * \gm g \gm * 1_{\mathbf{s}(x)g\gm})\,.
\end{eqnarray*}
$(\Longleftarrow)$ Suppose that $\pi_1 \times \mu_1$ and $\pi_2 \times \mu_2$ are unitarily equivalent and let $U$ be a unitary which establishes this equivalence. Then, since $\pi_1$ and $\pi_2$ are just the restrictions of, respectively, $\pi_1 \times \mu_1$ and $\pi_2 \times \mu_2$ we automatically have that $U\pi_1 U^* = \pi_2$. To see that $U\mu_1 U^* = \mu_2$ we just note that $U$ canonically establishes a unitary equivalence between the associated pre-$^*$-representations $\widetilde{\pi_1 \times \mu_1}$ and $\widetilde{\pi_2 \times \mu_2}$ of the multiplier algebra $M(C_c(\a / \gm) \times_{\alpha}^{alg} G / \gm)$. \qed\\

\section{Direct limits of sectional algebras}
\label{direc limit section}

In this section we will see how, for a finite index subgroup inclusion $K \subseteq H$, the algebra $C_c(\a / H)$ embedds canonically inside $C_c(\a / K)$. All these inlcusions are compatible with each other, so that we are able to form a certain direct limit $\mathcal{D}(\a)$ which will play an essential role for defining $C^*$-crossed products by Hecke pairs.

Throughout this section $\a$ denotes a Fell bundle over a discrete groupoid $X$, endowed with an action $\alpha$ of a group $G$. Will always assume in every statement of this section that the subgroup denoted by $H \subseteq G$ is such that the action $\alpha$ is $H$-good. We recall once again that we are using all the conventions and notation of \cite{palma crossed 11}.\\

\begin{prop}
\label{algebraic embeddings of CcA H into Cc A K}
 Suppose $K \subseteq H \subseteq G$ are subgroups such that $[H:K] < \infty$. Then, there is an embedding of $C_c(\a / H)$ into $C_c(\a / K)$ determined by
\begin{align*}
 [a]_{xH} \longmapsto \sum_{[h] \in \mathcal{S}_x \backslash H / K} [\alpha_{h^{-1}}(a)]_{x h K}\,.
\end{align*}

\end{prop}

\begin{rem}
 We have shown in \cite[Proposition 2.22]{palma crossed 11} that inside the multiplier algebra $M(C_c(\a))$ the element $[a]_{xH}$ decomposes as a sum of elements of $C_c(\a / K)$ as above. The point of Proposition \ref{algebraic embeddings of CcA H into Cc A K} is that this decomposition really defines an embedding of $C_c(\a / H)$ into $C_c(\a / K)$. Moreover, here we are not working inside $M(C_c(\a))$ anymore. Nevertheless this embedding of $C_c(\a / H)$ into $C_c(\a / K)$ is compatible with the embeddings of these algebras into  $M(C_c(\a))$ as we will see at the end of this section.\\
\end{rem}

{\bf \emph{Proof of Propostion \ref{algebraic embeddings of CcA H into Cc A K}:}} It is clear that the expression above is well-defined, since $[H:K] < \infty$, and it determines a linear map $\Phi :C_c(\a / H) \to C_c(\a / K)$. Moreover, it follows directly from \cite[Proposition 1.23]{palma crossed 11} that this map is injective. The fact that $\Phi$ preserves the involution follows from the following computation
\begin{eqnarray*}
 \Phi(([a]_{xH})^*) & = & \Phi([a^*]_{x^{-1}H})\;\; = \;\; \sum_{[h] \in \mathcal{S}_{x^{-1}} \backslash H / K} [\alpha_{h^{-1}}(a^*)]_{x^{-1} h K}\\
& = &   \sum_{[h] \in \mathcal{S}_x  \backslash H / K} [\alpha_{h^{-1}}(a^*)]_{x^{-1} h K} \;\; = \;\; \big( \sum_{[h] \in \mathcal{S}_x  \backslash H / K} [\alpha_{h^{-1}}(a)]_{x h K}\; \big)^* \\
& = & \Phi([a]_{xH})^*\,.
\end{eqnarray*}
Let us now check that $\Phi$ preserves products. If the pair $(xH, yH)$ is not composable, then no pair of the form $(xuK, ytK)$, with $u, t \in H$, is composable. Hence, in this case we have
\begin{align*}
 \Phi([a]_{xH}[b]_{yH}) \;=\; 0\; = \;\Phi([a]_{xH})\Phi([b]_{yH})\,.
\end{align*}
Suppose now the pair $(xH, yH)$ is composable, and let $\widetilde{h} \in H_{x,y}$. We have
\begin{eqnarray*}
 \Phi([a]_{xH}) \Phi([b]_{yH}) & = & \Big(\sum_{[u] \in \mathcal{S}_x \backslash H / K} [\alpha_{u^{-1}}(a)]_{x u K} \Big) \Big(\sum_{[t] \in \mathcal{S}_y \backslash H / K} [\alpha_{t^{-1}}(b)]_{y t K} \Big)\\
& = & \sum_{[t] \in \mathcal{S}_y \backslash H / K} \sum_{[u] \in \mathcal{S}_x \backslash H / K} [\alpha_{u^{-1}}(a)]_{x u K}[\alpha_{t^{-1}}(b)]_{y t K}\,.
\end{eqnarray*}
We now claim that
\begin{align*}
 \sum_{[u] \in \mathcal{S}_x \backslash H / K} [\alpha_{u^{-1}}(a)]_{x u K}[\alpha_{t^{-1}}(b)]_{y t K} = [\alpha_{t^{-1}\widetilde{h}^{-1}}(a)\alpha_{t^{-1}}(b)]_{(x\widetilde{h}t)(yt)K}\,.
\end{align*}
 To see this we notice that for $u = \widetilde{h}t$ we do have that the pair $(x\widetilde{h}tK, ytK)$ is composable and $[\alpha_{t^{-1}\widetilde{h}^{-1}}(a)]_{x \widetilde{h} t K}[\alpha_{t^{-1}}(b)]_{ytK} = [\alpha_{t^{-1}\widetilde{h}^{-1}}(a)\alpha_{t^{-1}}(b)]_{(x\widetilde{h}t)(yt)K}$. Now if $[u] \in  \mathcal{S}_x \backslash H / K$ is such that the pair $(xuK, ytK)$ is composable, then $\mathbf{s}(x)uK = \mathbf{r}(y) t K$. Since the pair $(x\widetilde{h}tK, ytK)$ is composable we also have $\mathbf{s}(x)\widetilde{h}tK = \mathbf{r}(y) t K$. Thus, $\mathbf{s}(x)uK = \mathbf{s}(x)\widetilde{h}tK$, i.e. $[u] = [\widetilde{h}t]$ by \cite[Proposition 1.23]{palma crossed 11}. This proves our claim and therefore we get
\begin{eqnarray*}
 \Phi([a]_{xH}) \Phi([b]_{yH}) & = & \sum_{[t] \in \mathcal{S}_y \backslash H / K} \sum_{[u] \in \mathcal{S}_x \backslash H / K} [\alpha_{u^{-1}}(a)]_{x u K}[\alpha_{t^{-1}}(b)]_{y t K}\\
 & = & \sum_{[t] \in \mathcal{S}_y \backslash H / K} [\alpha_{t^{-1}\widetilde{h}^{-1}}(a)\alpha_{t^{-1}}(b)]_{(x\widetilde{h}t)(yt)K}\\
 & = & \sum_{[t] \in \mathcal{S}_y \backslash H / K} [\alpha_{t^{-1}}(\alpha_{\widetilde{h}^{-1}}(a)b)]_{(x\widetilde{h}y)tK}\,.
\end{eqnarray*}
Recall that since the $G$-action on $\a$ is $H$-good we have
\begin{align*}
 \mathcal{S}_{y} \cap H =\mathcal{S}_{\mathbf{s}(y)} \cap H = \mathcal{S}_{x\widetilde{h} y} \cap H\,.
\end{align*}
Hence, using \cite[Proposition 1.22]{palma crossed 11}, we have bijections
\begin{align*}
 \mathcal{S}_y \backslash H / K \; \cong \; (\mathcal{S}_y \cap H) \backslash H / K \; \cong  \; (\mathcal{S}_{x\widetilde{h} y} \cap H) \backslash H / K \; \cong \; \mathcal{S}_{x\widetilde{h} y}  \backslash H / K\,,
\end{align*}
determined by the maps $[t] \to [t]$. Therefore we get
\begin{eqnarray*}
\Phi([a]_{xH}) \Phi([b]_{yH}) & = & \sum_{[t] \in \mathcal{S}_{x\widetilde{h} y} \backslash H / K} [\alpha_{t^{-1}}(\alpha_{\widetilde{h}^{-1}}(a)b)]_{(x\widetilde{h}y)tK}\\
 & = & \Phi([\alpha_{\widetilde{h}^{-1}}(a)b]_{x\widetilde{h}yH})\\
 & = & \Phi([a]_{xH} [b]_{yH})\,.
\end{eqnarray*}
Hence, $\Phi$ is an embedding of $C_c(\a / H)$ into $C_c(\a / K)$. \qed\\

The canonical embeddings described in Proposition \ref{algebraic embeddings of CcA H into Cc A K} are all compatible, as the following result shows:\\

\begin{prop}
\label{compatibility of embeddings CcAH CcAK}
 Suppose that $L \subseteq K \subseteq H$ are subgroups of $G$ such that $[H:L] < \infty$. The canonical embedding of $C_c(\a / H)$ into $C_c(\a / L)$ factors through the canonical embeddings of $C_c(\a / H)$ into $C_c(\a / K)$, and $C_c(\a / K)$ into $C_c(\a / L)$. In other words, the following diagram of canonical embeddings commutes:
\begin{displaymath}
\xymatrix{ C_c(\a / H) \ar[r] \ar@/_1.5pc/[rr] &  C_c(\a / K) \ar[r] & C_c(\a / L)}\,. 
\end{displaymath}\\
\end{prop}

{\bf \emph{Proof:}} Let us denote by $\Phi_1 : C_c(\a / H) \to C_c(\a / K)$, $\Phi_2 : C_c(\a / K) \to C_c(\a / L)$ and $\Phi_3 : C_c(\a / H) \to C_c(\a / L)$ the canonical embeddings. We want to prove that $\Phi_3 = \Phi_2 \circ \Phi_1$. For this it is enough to check this equality on elements of the form $[a]_{xH}$. We have
\begin{eqnarray*}
 \Phi_2 \circ \Phi_1 ([a]_{xH}) & = & \sum_{[h] \in \mathcal{S}_x \backslash H / K} \Phi_2([\alpha_{h^{-1}}(a)]_{x h K})\\
& = & \sum_{[h] \in \mathcal{S}_x \backslash H / K} \sum_{[k] \in \mathcal{S}_{xh} \backslash K / L} [\alpha_{k^{-1}h^{-1}}(a)]_{x h kL}\,.
\end{eqnarray*}
We claim that if $h_1, \dots, h_n \in H$ is a set of representatives for $\mathcal{S}_x \backslash H / K$, and if $k^i_1, \dots, k^i_{r_i}$ is a set of representatives of $\mathcal{S}_{xh_i} \backslash K / L$ for each $i=1, \dots, n$, then the set of all products of the form $h_ik^i_j$ is a set of representatives for $\mathcal{S}_x \backslash H / L$. Let us start by proving that every two such products correspond to distinct elements of $\mathcal{S}_x \backslash H / L$. In other words, we want to show that if $[h_ik^i_j] = [h_lk^l_p]$ in $S_x \backslash H / L$, then $h_i = h_l$ and $k^i_j = k^l_p$. To see this we notice that the equality $[h_ik^i_j] = [h_lk^l_p]$ means that $xh_ik^i_jL = x h_lk^l_p L $ (see \cite[Proposition 1.23]{palma crossed 11}), and therefore $xh_iK = x h_lK$, i.e. $[h_i] = [h_l]$ in $S_x \backslash H / K$, hence $h_i = h_l$ because these form a set a of representatives of $\mathcal{S}_x \backslash H / K$. Now, the equality $xh_ik^i_jL = x h_ik^i_p L$ means that $k^i_j = k^i_p$ for the same reasons. Now it remains to prove that any element of $[h] \in \mathcal{S}_x \backslash H / L$ has a representative of the form $h^i k^i_j$. To see this, first we take $h_i$ such that $xhK = xh_iK$, and we consider an element $k \in K$ such that $xh = xh_ik$, obtaining $xhL = xh_ikL$. Now we take $k^i_j$ such that $xh_ikL = xh_i k^i_j L$, and the result follows.

After proving the above claim we can now write
\begin{eqnarray*}
  \Phi_2 \circ \Phi_1 ([a]_{xH}) & = &  \sum_{[h] \in \mathcal{S}_x \backslash H / K} \sum_{[k] \in \mathcal{S}_{xh} \backslash K / L} [\alpha_{k^{-1}h^{-1}}(a)]_{x h kL}\\
& = &  \sum_{[\widetilde{h}] \in \mathcal{S}_x \backslash H / L} [\alpha_{\widetilde{h}^{-1}}(a)]_{x \widetilde{h} L}\\
& = & \Phi_3([a]_{xH})\,.
\end{eqnarray*}
This finishes the proof. \qed\\

Suppose now that $(G, \gm)$ is a Hecke pair for which the $G$-action on the Fell bundle $\a$ is $\gm$-good. We define the set $\mathcal{C}$\index{C@$\mathcal{C}$} as the set of all finite intersections of conjugates of $\gm$, i.e.
\begin{align}
 \mathcal{C} := \Big\{ \bigcap_{i=1}^n g_i \gm g_i^{-1} : n \in \mathbb{N}, g_1, \dots, g_n \in G \Big\}\,.
\end{align}
The set $\mathcal{C}$ becomes a directed set with respect to the partial order given by reverse inclusion of subgroups, i.e. $H_1 \leq H_2 \Leftrightarrow H_1 \supseteq H_2$, for any $H_1, H_2 \in \mathcal{C}$.

Since we are assuming that $(G,\gm)$ is a Hecke pair it is not difficult to see that for any $H_1, H_2 \in \mathcal{C}$ we have
\begin{align*}
 H_1 \leq H_2 \;\Longrightarrow\; [H_1 : H_2] < \infty\,.
\end{align*}
Also, since we are assuming that the $G$-action on $\a$ is $\gm$-good and this property passes to conjugates and subgroups, it follows automatically that the action is also $H$-good, for any $H \in \mathcal{C}$.

The observations in the previous paragraph together with Proposition \ref{compatibility of embeddings CcAH CcAK} imply that $\{C_c(\a / H) \}_{H \in \mathcal{C}}$ is a directed system of $^*$-algebras. Let us denote by $D(\a)$\index{D3 (A)@$\mathcal{D}(\a)$} the $^*$-algebraic direct limit of this directed system, i.e.
\begin{align}
 D(\a) := \lim_{H \in \mathcal{C}} C_c(\a / H)\,.
\end{align}

There is an equivalent way of defining the algebra $D(\a)$, by viewing it as the $^*$-subalgebra of $M(C_c(\a))$ generated by all the $C_c(\a / H)$ with $H \in \mathcal{C}$, as we prove in the next result. This characterization of $D(\a)$ is also a very useful one.\\

\begin{prop}
\label{D(A) embedds in MCcA}
Let $K \subseteq H$ be subgroups of $G$ such that $[H:K] < \infty$. Then the following diagram of canonical embeddings commutes:
 \begin{align}
 \label{diagram of embeddings CcaH CcAK MCcA}
\xymatrix{ C_c(\a / H) \ar[r] \ar[dr] &  C_c(\a / K) \ar[d]\\
  & M(C_c(\a))\,.} 
\end{align}
As a consequence, $D(\a)$ is $^*$-isomorphic to the $^*$-subalgebra of $M(C_c(\a))$ generated by all the $C_c(\a / H)$ with $H \in \mathcal{C}$.\\
\end{prop}

{\bf \emph{Proof:}} We have to show that, inside $M(C_c(\a))$, we have
\begin{align*}
[a]_{xH} = \sum_{[h] \in \mathcal{S}_x \backslash H / K} [\alpha_{h^{-1}}(a)]_{xh K}\,,
\end{align*}
for all $x, y \in X$, $a \in \a_x$ and $b \in \a_y$. This was proven in \cite[Proposition 2.22]{palma crossed 11}.

Commutativity of the diagram (\ref{diagram of embeddings CcaH CcAK MCcA}) then implies, by universal properties, that there exists a $^*$-homomorphism from $D(\a)$ to $M(C_c(\a))$ whose image is precisely the $^*$-subalgebra generated by all $C_c(\a / H)$, with $H \in \mathcal{C}$. This $^*$-homomorphism is injective since all the maps in the diagram (\ref{diagram of embeddings CcaH CcAK MCcA}) are injective. \qed\\

It is clear that the action $\overline{\alpha}$ gives rise to an action of $G$ on $\mathcal{D}(\a)$, which we will still denote by $\overline{\alpha}$. This can be seen either directly, or simply by noticing that the action $\overline{\alpha}$ on $M(C_c(\a))$ takes $\mathcal{D}(\a)$ to itself (since for a given $g \in G$ it takes $C_c(\a / H)$ to $C_c(\a / g H g^{-1})$).

The algebra $\mathcal{D}(\a)$ will play an essential role in the definition of the various $C^*$-crossed products by Hecke pairs, particularly the reduced ones. There are two reduced $C^*$-crossed products by Hecke pairs which are of particular interest to us, and these are $C^*_r(\a/ \gm) \times_{\alpha, r} G / \gm$ and $C^*(\a / \gm) \times_{\alpha, r} G / \gm$. These will be defined and studied, in a singular approach, in Section \ref{reduced C crossed products section}, but for that we need first to understand  how the canonical embeddings
\begin{align}
\label{canonical embedding CcAH into CcAK preliminary section of reduced C cross prod}
 C_c(\a / H) \to C_c(\a / K)\,,
\end{align}
defined in Proposition \ref{algebraic embeddings of CcA H into Cc A K} for $K \subseteq H$ such that $[H:K] < \infty$, behave with respect the full and reduced $C^*$-completions. The goal of next subsections is exactly to show that these embeddings always give rise to embeddings in the two canonical $C^*$-completions
\begin{align*}
 C^*_r(\a / H) \to C^*_r(\a / K) \qquad\;\text{and}\qquad\; C^*(\a / H) \to C^*(\a / K)\,,
\end{align*}
so that we are able to form the useful $C^*$-direct limits $\lim_{H \in C} C^*_r(\a / H)$ and $\lim_{H \in C} C^*(\a / H)$.\\

\subsection{Reduced completions $C^*_r(\a / H)$}
\label{reduced completions of Cc(a H) section}

The purpose of this subsection is to prove the following result:\\

\begin{thm}
\label{embedding of Cc(A H) in Cc(A K) completes to reduced}
 Let $K \subseteq H \subseteq G$ be subgroups such that $[H:K] < \infty$. The canonical embedding of $C_c(\a / H)$ into $C_c(\a / K)$ completes to an embedding of $C^*_r(\a / H)$ into $C^*_r(\a / K)$.\\
\end{thm}

In order to prove this result we need to establish some notation and some lemmas first. Even though Theorem \ref{embedding of Cc(A H) in Cc(A K) completes to reduced} is stated for subgroups $K \subseteq H$ for which we have a finite index $[H:K]$ we will state and prove the two following lemmas in greater generality, as it will be convenient later on.

Recall, from \cite[Proposition 2.21]{palma crossed 11}, that for any two subgroups $K \subseteq H$ of $G$ for which the $G$-action is $H$-good we have that, inside $M(C_c(\a))$, the algebra $C_c(\a / H)$ acts on $C_c(\a / K)$ in the following way:
\begin{align*}
 [a]_{xH} [b]_{yK} & = \begin{cases}
 [\alpha_{\widetilde{h}^{-1}}(a)b]_{x\widetilde{h}yK}\,, \qquad\; \text{if}\;\; H_{x,y} \neq \emptyset \\
 0, \qquad\qquad\quad \text{otherwise,}\;\;
\end{cases}
\end{align*}
where $\widetilde{h}$ is any element of $H_{x,y}$. As a consequence, this action of $C_c(\a / H)$ on $C_c(\a / K)$ defines a $^*$-homomorphism
\begin{align*}
 C_c(\a / H) \to M(C_c(\a / K))\,.
\end{align*}
It could be proven (in the same fashion as \cite[Theorem 2.16]{palma crossed 11}) that the $^*$-homomorphism above is in fact an embedding, but we will not need this fact here. We now make the following definition:\\

\begin{df}
\label{def of right A B bimodule and right-pre-hilbert}
 Suppose $A$ is $^*$-algebra and $B$ is a $C^*$-algebra. A \emph{right $A - B$ bimodule} $X$ is a (right) inner product $B$-module (in the sense of \cite[Definition 2.1]{morita equiv}) which is also a left $A$-module satisfying:
\begin{align*}
 a(xb) &=(ax)b\,,\\
\langle a x, y \rangle_B & = \langle x, a^*y\rangle_B\,,
\end{align*}
for all $x,y \in X$, $a \in A$ and $b \in B$.

Given a right $A - B$ bimodule $X$ we will say that $A$ \emph{acts by bounded operators} on $X$ if for any $a \in A$ there exists $C > 0$ such that
\begin{align*}
 \|ax\|_B \leq C \|x \|_B\,,
\end{align*}
for every $x \in X$, where $\| \cdot\|_B$ is the norm induced by $\langle \cdot\,,\, \cdot \rangle_B$.\\
\end{df}

If $A$ is a $^*$-algebra which has an enveloping $C^*$-algebra $C^*(A)$, then any right $A - B$ bimodule where $A$ acts by bounded operators can be completed to a right-Hilbert $C^*(A) - B$ bimodule.\\

\begin{lemma}
\label{embedding of Cc A H into multiplier of completion of Cc A}
 Let $K \subseteq H$ be subgroups of $G$ and let $D$ be a $C^*$-algebra. Suppose $C_c(\a / K)$ is an inner product $D$-module, denoted by $C_c(\a / K)_D$. Assume furthermore that $C_c(\a / K)_D$ is a right $C_c(\a / K) - D$ bimodule and also a right $C_c(\a / H) - D$ bimodule, where $C_c(\a / K)$ acts on itself by right multiplication and $C_c(\a  / H)$ acts on $C_c(\a / K)$ in the canonical way. 

If $C_c(\a / K)$ acts on $C_c(\a / K)_D$ by bounded operators, then $C_c(\a / H)$ also acts on $C_c(\a / K)_D$ by bounded operators.\\
\end{lemma}

{\bf \emph{Proof:}} Suppose that $C_c(\a / K)$ acts on $C_c(\a / K)_D$ by bounded operators. We need to show that $C_c(\a / H)$ also acts on $C_c(\a / K)_D$ by bounded operators, with respect to the norm $ \| \cdot \|_D$ induced by the $D$-valued inner product in $C_c(\a / K)_D$. For this it is enough to prove that the maps
\begin{align*}
 [a]_{xH}: C_c(\a / K) \to C_c(\a / K)\,,
\end{align*}
are bounded with respect to the norm $\| \cdot \|_D$. Moreover, from the fact that $([a]_{xH})^*[a]_{xH} = ([a^*a])_{s(x)H}$ it actually suffices to show that for any unit $u \in X^0$ the mapping $[a]_{uH} : C_c(\a / K) \to C_c(\a / K)$ is bounded with respect to the  norm $\| \cdot \|_D$, .

 As we have seen in \cite[end of Section 2.3]{palma crossed 11} we can write any element $f \in C_c(\a / K)$ as a sum of the form $f = \sum_{yK \in X / K} [f(y  )]_{yK}$. Furthermore, we can split the sum according the ranges of elements, i.e.
\begin{align*}
 f = \sum_{yK \in X / K} [f(y)]_{yK} = \sum_{vK \in X^0 / K} \sum_{\substack{yK \in X / K \\ \mathbf{r}(y)K = vK}} [f(y)]_{yK}\,.
\end{align*}
Applying the multiplier $[a]_{uH}$ to this element we get
\begin{eqnarray*}
[a]_{uH} f & = & [a]_{uH} \sum_{vK \in X^0 / K} \sum_{\substack{yK \in X / K \\ \mathbf{r}(y)K = vK}} [f(y)]_{yK} \\
 & = &  \sum_{vK \in X^0 / K} \sum_{\substack{yK \in X / K \\ \mathbf{r}(y)K = vK}} [a]_{uH}[f(y)]_{yK} \\
 & = &  \sum_{vK \in X^0 / K} \sum_{\substack{yK \in X / K \\ \mathbf{r}(y)K = vK}} [\alpha_{\widetilde{h_v}^{-1}}(a)]_{vH}[f(y)]_{yK} \,,
\end{eqnarray*}
where $h_v$ is any element of $H_{u,v}$. Hence, if $k_{v,y}$ is any element of $K_{v,y} \subseteq H_{v,y}$ we get
\begin{eqnarray*}
 &  = &  \sum_{vK \subseteq uH} \sum_{\substack{yK \in X / K \\ \mathbf{r}(y)K = vK}} [\alpha_{\widetilde{k_{v,y}}^{-1}}(\alpha_{\widetilde{h_v}^{-1}}(a))f(y)]_{yK}\\
  &  = &  \sum_{vK \subseteq uH} \sum_{\substack{yK \in X / K \\ \mathbf{r}(y)K = vK}} [\alpha_{\widetilde{h_v}^{-1}}(a)]_{vK}[f(y)]_{yK}\,.
\end{eqnarray*}
Since $f$ has compact support, there are a finite number of elements $v_1K, \dots, v_nK \subseteq uH$ such that
\begin{eqnarray*}
 [a]_{uH}f & = & \sum_{i=1}^n \sum_{\substack{y K \in X / K \\ \mathbf{r}(y)K = v_iK}} [\alpha_{\widetilde{h_{v_i}}^{-1}}(a)]_{v_iK}[f(y)]_{yK}\\
& = & \big( \sum_{i=1}^n [\alpha_{\widetilde{h_{v_i}}^{-1}}(a)]_{v_iK} \big) \big( \sum_{yK \in X / K } [f(y)]_{yK} \big)\\
& = & \big( \sum_{i=1}^n [\alpha_{\widetilde{h_{v_i}}^{-1}}(a)]_{v_iK} \big) \,f \,.
\end{eqnarray*}

Our assumptions say that left multiplication by elements of $C_c(\a / K)$ is continuous with respect to $\| \cdot \|_D$. Denoting by $\overline{C_c(\a / K)_D}$ the completion of $C_c(\a / K)_D$ as a Hilbert $D$-module, we have that  every element of $C_c(\a / K)$ uniquely defines an element of $\mathcal{L}\big(\overline{C_c(\a / K)_D} \big)$. Denoting by $\| \cdot \|_{\mathcal{L}\big(\overline{C_c(\a / K)_D}\big)}$ the operator norm in $\mathcal{L}\big(\overline{C_c(\a /K)_D}\big)$, we have
\begin{eqnarray*}
 \| [a]_{uH}f \;\|_D & = & \| \big( \sum_{i=1}^n [\alpha_{\widetilde{h_{v_i}}^{-1}}(a)]_{v_iK} \big) f\|_D\\
 & \leq & \|  \sum_{i=1}^n [\alpha_{\widetilde{h_{v_i}}^{-1}}(a)]_{v_iK} \|_{\mathcal{L}\big(\overline{C_c(\a / K)_D}\big)} \; \|  f\|_{D}\,,
\end{eqnarray*}

Now we notice that we can canonically see $ \sum_{i=1}^n [\alpha_{\widetilde{h_{v_i}}^{-1}}(a)]_{v_iK} $ as an element of the direct sum of $C^*$-algebras $(\a/K)_{v_1K} \oplus \dots \oplus (\a / K)_{v_nK}$, from which we must have, by uniqueness of $C^*$-norms on $C^*$-algebras,
\begin{align*}
 \| \sum_{i=1}^n [\alpha_{\widetilde{h_{v_i}}^{-1}}(a)]_{v_iK} \|_{\mathcal{L} \big(\overline{C_c(\a / K)_D}\big)} = \max_i \|[\alpha_{\widetilde{h_{v_i}}^{-1}}(a)]\| =  \max_i \|\alpha_{\widetilde{h_{v_i}}^{-1}}(a)\| = \|a\|\,.
\end{align*}
Hence we conclude that $\| [a]_{uH}f \|_D \leq \|a\| \,\|f\|_D$, i.e. $[a]_{uH}$ is bounded.  \qed\\

Let us now consider $C_c(\a / K)$ as the right $C_c(\a / K) - C_0(\a^0 / K)$ bimodule whose completion is the right-Hilbert bimodule $_{C^*(\a / K)} L^2(\a / K)_{C_0(\a^0 / K)}$. We claim that the canonical action of $C_c(\a / H)$ on $C_c(\a / K)$ makes $C_c(\a / K)$ into a right $C_c(\a / H) - C_0(\a / K)$ bimodule. The fact that $f_1 (\xi f_2) = (f_1 \xi) f_2$, for any $f_1 \in C_c(\a / H)$, $\xi \in C_c(\a / K)$ and $f_2 \in C_0(\a^0 / K)$, is obvious. Thus, we only need to check that $\langle f \xi\,,\,\eta \rangle_{C_0(\a^0 / K)} = \langle  \xi\,,\, f^*\eta \rangle_{C_0(\a^0 / K)}$, for any $f \in C_c(\a / H)$ and $\xi, \eta \in C_c(\a / K)$. This is also easy to see because, by definition,
\begin{eqnarray*}
 \langle f \xi\,,\,\eta \rangle_{C_0(\a^0 / K)} & = & ((f\xi)^*\eta)|_{C_0(\a / K)}\\
 & = & (\xi^*(f^*\eta))|_{C_0(\a / K)}\\
 & = & \langle  \xi\,,\, f^*\eta \rangle_{C_0(\a^0 / K)}\,.
\end{eqnarray*}
Hence, we are under the conditions of Lemma \ref{embedding of Cc A H into multiplier of completion of Cc A}, and therefore the action of $C_c(\a / H)$ on $C_c(\a / K)_{C_0(\a^0/K)}$ is by bounded operators. Hence,  the right $C_c(\a / H) - C_0(\a^0 / K)$ bimodule $C_c(\a / K)$ can be completed to a right-Hilbert bimodule $_{C^*(\a / H)} L^2(\a / K)_{C_0(\a^0 / K)}$.\\

 \begin{lemma}
\label{lemma C correspondeces H and K}
 The $^*$-homomorphism $\Phi: C^*(\a / H) \to \mathcal{L}(L^2(\a / K))$ associated with the right-Hilbert bimodule $_{C^*(\a / H)} L^2(\a / K) _{C_0(\a^0 / K)}$ has the same kernel as the canonical map $\Lambda: C^*(\a / H) \to C^*_r(\a / H)$.\\
\end{lemma}

{\bf \emph{Proof:}} The proof of this fact is essentially an adaptation of the proof of \cite[Proposition 2.10]{full}, and is achieved by exhibiting two isomorphic right-Hilbert  $C^*(\a / H) - C_0(\a^0 / K)$ bimodules $Y$ and $Z$ such that the $^*$-homomorphisms from $C^*(\a / H)$ into $\mathcal{L}(Y)$ and $\mathcal{L}(Z)$ have the same kernels as $\Lambda$ and $\Phi$ respectively.

We naturally have a  right-Hilbert bimodule $_{C_0(\a^0 / H)} C_0(\a^0 / K)_{C_0(\a^0 / K)}$, where the action of $C_0(\a^0 / H)$ on $C_0(\a^0 / K)$ extends the action of $C_c(\a^0 / H)$ on $C_c(\a^0 / K)$. We define $Y$ as the balanced tensor product of the right-Hilbert bimodules $_{C^*(\a / H)} L^2(\a / H) _{C_0(\a^0 / H)}$ and $_{C_0(\a^0 / H)} C_0(\a^0 / K)_{C_0(\a^0) / K}$, i.e.
\begin{align*}
 Y := L^2(\a / H) \otimes_{C_0(\a^0 / H)} C_0(\a^0 / K)\,.
\end{align*}
Since $C_0(\a^0 / H)$ acts faithfully on $C_0(\a^0 / K)$, the associated $^*$-homomorphism of $C^*(\a / H)$ to $\mathcal{L}(Y)$ has the same kernel as  $\Lambda$. We define $Z$ simply as
\begin{align*}
 _{C^*(\a / H)} Z_{C_0(\a^0 / K)}:=\; _{C^*(\a / H)} L^2(\a / K)_{C_0(\a^0 / K)}\,.
\end{align*}

We now want to define an isomorphism $\Psi: L^2(\a / H) \otimes_{C_0(\a^0 / H)} C_0(\a^0 / K) \to L^2(\a / K)$ of Hilbert $C^*(\a / H) - C_0(\a^0 / K)$ bimodules. We start by defining 
\begin{align*}
 \Psi_0 : C_c(\a / H) \otimes_{C_c(\a^0 / H)} C_c(\a^0 / K) &\longrightarrow L^2(\a / K)\,,\\
 \Psi_0 (f_1 \otimes f_2) &:= f_1 \cdot f_2\,.
\end{align*}
It is easy to see that $\Psi_0$ is well-defined. To see that $\Psi_0$ preserves the inner products it is enough to check on the generators. So let $[a]_{xH}, [b]_{yH} \in C_c(\a / H)$ and $[c]_{uK}, [d]_{vK} \in C_c(\a^0 / K)$, with $u, v \in X^0$. We have
\begin{eqnarray*}
& & \langle \Psi_0([a]_{xH} \otimes [c]_{uK})\,,\, \Psi_0([b]_{yH} \otimes [d]_{vK}) \rangle_{C_0(\a^0/ K)}\;\; =\;\;\\
 & = & \langle [a]_{xH} [c]_{uK}\,,\, [b]_{yH}[d]_{vK} \rangle_{C_0(\a^0 / K)}\\
& = & ([c^*]_{uK} [a^*]_{x^{-1}H} [b]_{yH} [d]_{vK})|_{C_0(\a^0 / K)}\,.
\end{eqnarray*}
Now the product $([c^*]_{uK} [a^*]_{x^{-1}H} [b]_{yH} [d]_{vK})|_{C_0(\a^0 / K)}$ is automatically zero unless $vK = uK$, $xH = yH$ and $vK \subseteq \mathbf{s}(y)H$, in which case we necessarily have that $([c^*]_{uK} [a^*]_{x^{-1}H} [b]_{yH} [d]_{vK})|_{C_0(\a^0 / K)} = [c^*]_{uK} [a^*]_{x^{-1}H} [b]_{yH} [d]_{vK}$. On the other hand,
\begin{eqnarray*}
 & & \langle [a]_{xH} \otimes [c]_{uK}\,,\, [b]_{yH} \otimes [d]_{vK} \rangle_{C_0(\a^0 / K)} \,\; =\;\;\\
 & = & \langle  [c]_{uK}\,,\, \langle [a]_{xH} \,,\, [b]_{yH} \rangle_{C_0(\a^0 / H)}\, [d]_{vK} \rangle_{C_0(\a^0 / K)}\\
& = & [c^*]_{uK} ([a^*]_{x^{-1}H} [b]_{yH})|_{C_0(\a^0 / H)}\, [d]_{vK}
\end{eqnarray*}
Now the product $[c^*]_{uK} ([a^*]_{x^{-1}H} [b]_{yH})|_{C_0(\a^0 / H)}\, [d]_{vK}$ is automatically zero unless $vK = uK$, $xH = yH$ and $vK \subseteq \mathbf{s}(y)H$, in which case we necessarily have that $[c^*]_{uK} ([a^*]_{x^{-1}H} [b]_{yH})|_{C_0(\a^0 / H)}\, [d]_{vK} = [c^*]_{uK} [a^*]_{x^{-1}H} [b]_{yH} [d]_{vK}$. Hence, we conclude that $\Psi_0$ preserves the inner products.

Now, if $f_1, f_2 \in C_c(\a / H)$ and $f_3 \in C_c(\a^0 / K)$ we have
\begin{align*}
 \Psi_0(f_1 (f_2 \otimes f_3))  =  \Psi_0(f_1 f_2 \otimes f_3) = f_1f_2 f_3 = f_1 \Psi_0(f_2 \otimes f_3)\,.
\end{align*}
Thus, $\Psi_0$ preserves the left module actions. Let us now check that $\Psi_0$ has a dense image in $L^2(\a / K)$. It is enough to prove that all generators $[a]_{xK} \in C_c(\a / K)$ are in closure of the image of $\Psi_0$, since their span is dense in $L^2(\a / K)$. To see this, let $\{e^{\lambda}\}_{\lambda}$ be an approximate identity of $\a_{\mathbf{s}(x)}$. We have
\begin{align*}
 \Psi_0([a]_{xH} \otimes [e^{\lambda}]_{\mathbf{s}(x)K})  =  [a]_{xH}[e^{\lambda}]_{\mathbf{s}(x)K} = [ae^{\lambda}]_{xK}\,.
\end{align*}
We then get
\begin{eqnarray*}
 \| [ae^{\lambda}]_{xK} - [a]_{xK} \|^2_{L^2(\a / K)} & = & \| [ae^{\lambda} - a]_{xK} \|^2_{L^2(\a / K)}\\
& = & \| \big([ae^{\lambda} - a]^*[ae^{\lambda} - a]\big)_{\mathbf{s}(x) K} \|_{C_0(\a^0 / K)}\\
& = & \| [ae^{\lambda} - a]^*[ae^{\lambda} - a] \|\\
& = & \| e^{\lambda} a^*a e^{\lambda} - e^{\lambda}a^*a - a^*ae^{\lambda} + a^*a\|\,.
\end{eqnarray*}
Noticing that $a^*a \in \a_{\mathbf{s}(x)}$, we then have that
\begin{eqnarray*}
 & \leq & \| e^{\lambda} a^*a e^{\lambda} - e^{\lambda}a^*a \| + \| - a^*ae^{\lambda} + a^*a\|\\
& \leq & \|  a^*a e^{\lambda} - a^*a \| + \| - a^*ae^{\lambda} + a^*a\|\\
 & \longrightarrow & 0\,.
\end{eqnarray*}

Thus, we conclude that $\Psi_0$ has dense range. Hence, from \cite[Lemma 2.9]{full}, it follows that $\Psi_0$ extends to an isomorphism of the right-Hilbert $C^*(\a / H) - C_0(\a^0 / K)$ bimodules $Y$ and $Z$. \qed\\

{\bf \emph{Proof of Theorem \ref{embedding of Cc(A H) in Cc(A K) completes to reduced}:}} The image of $C^*(\a  / H)$ in $\mathcal{L}(L^2(\a / K))$ is isomorphic to $C^*_r(\a / H)$ by Lemma \ref{lemma C correspondeces H and K}. On the other hand, the image of $C^*(\a / H)$ in $\mathcal{L}(L^2(\a / K))$ is simply the completion of $C_c(\a / H)$ as a subalgebra of $C^*_r(\a / K)$. Hence, we conclude that the canonical embedding of $C_c(\a / H)$ into $C_c(\a / K)$ completes to an embedding of $C^*_r(\a / H)$ into $C_r^*(\a / K)$. \qed\\

It follows from Theorem \ref{embedding of Cc(A H) in Cc(A K) completes to reduced} and Proposition \ref{compatibility of embeddings CcAH CcAK} that for any subgroups $L \subseteq K \subseteq H$ such that $[H:L] < \infty$ the following diagram of canonical embeddings commutes
\begin{displaymath}
\xymatrix{ C_r^*(\a / H) \ar[r] \ar@/_1.5pc/[rr] &  C_r^*(\a / K) \ar[r] & C_r^*(\a / L)} \,.
\end{displaymath}\\

 Hence, we have a direct system of $C^*$-algebras $\{C_r^*(\a / H)\}_{H \in \cl}$. Let us denote by $\mathcal{D}_r(\a)$\index{Dr(A)@$\mathcal{D}_r(\a)$} its corresponding $C^*$-algebraic direct limit
\begin{align}
 \d_r(\a) := \lim_{H \in \cl} C^*_r(\a / H)\,.
\end{align}
We notice that the algebra $\d(\a)$ is a dense $^*$-subalgebra of $\d_r(\a)$. We now want to show that the action $\overline{\alpha}$ of $G$ on $\d(\a)$ extends to $\d_r(\a)$.\\

\begin{thm}
\label{action extends to DrA}
 The action $\overline{\alpha}$ of $G$ on $\d(\a)$ extends uniquely to an action of $G$ on $\d_r(\a)$ and is such that $\overline{\alpha}_g$ takes $C^*_r(\a / H)$ to $C^*_r(\a / g H g^{-1})$, for any $g \in G$.\\
\end{thm}

{\bf \emph{Proof:}} We have a canonical isomorphism between the right-Hilbert bimodules $_{C^*(\a / H)} L^2(\a / H)_{C_0(\a^0/ H)}$ and $_{C^*(\a / gHg^{-1})} L^2(\a / gHg^{-1})_{C_0(\a^0/ gHg^{-1})}$, that is determined by the canonical isomorphisms
$C_c(\a / H) \to C_c(\a / gHg^{-1})$ and $C_c(\a^0/H) \to C_c(\a^0/ gHg^{-1})$ defined by $\overline{\alpha}_g$, i.e. defined respectively by
\begin{align*}
 [a]_{xH} \mapsto [\alpha_g(a)]_{xg^{-1}gHg^{-1}}\,, \qquad\text{and}\qquad [b]_{uH} \mapsto [\alpha_g(b)]_{ug^{-1}gHg^{-1}}\,,
\end{align*}
where $x \in X$, $u \in X^0$, $a \in \a_x$ and $b \in \a_u$. Since $C_r^*(\a/H)$ is the image of $C^*(\a / H)$ inside $\mathcal{L}(L^2(\a / H))$, and similarly for $C^*_r(\a / g H g^{-1})$, we conclude that the isomorphism $C_c(\a / H) \cong C_c(\a / g Hg^{-1})$ defined by $\overline{\alpha}_g$ extends  to an isomorphism $C^*_r(\a / H) \cong C^*_r(\a / g Hg^{-1})$. Since $C^*_r(\a / gHg^{-1})$ is embedded in $D_r(\a)$, we can see $\overline{\alpha}_g$ as an injective $^*$-homomorphism from $C^*_r(\a / H)$ into $D_r(\a)$.

A routine computation shows that the following diagram of canonical injections commutes:
 \begin{align*}
\xymatrix{ C_r^*(\a / H) \ar[r] \ar[dr]_{\overline{\alpha}_g} &  C_r^*(\a / K) \ar[d]^{\overline{\alpha}_g}\\
  & D_r(\a)\,.} 
\end{align*}
Hence, we obtain an injective $^*$-homomorphism from $\d_r(\a)$ to itself, which we still denote by $\overline{\alpha}_g$, and which extends the usual map $\overline{\alpha}_g$ from $D(\a)$ to itself. It is also clear that this map is surjective, and that for $g, h \in G$ we have $\overline{\alpha}_{gh} = \overline{\alpha}_g \circ \overline{\alpha}_h$, so that we get an action of $G$ on $D_r(\a)$ which extends the usual action of $G$ on $D(\a)$. \qed\\

\subsection{Maximal completions $C^*(\a / H)$}
\label{maximal completions of Cc(a H) section}

The purpose of this subsection is to prove the following result:\\

\begin{thm}
\label{embedding CcAH into CcAK passes to Calgebra}
 Let $K \subseteq H $ be subgroups of $G$ such that $[H:K] < \infty$. The canonical embedding of $C_c(\a / H)$ into $C_c(\a / K)$ completes to a nondegenerate embedding of $C^*(\a / H)$ into $C^*(\a / K)$.\\
\end{thm}

In order to prove this result we will need to know how to ``extend" a representation of $C_c(\a /H)$ to a representation of $C_c(\a / K)$ on a larger Hilbert space.\\

\begin{df}
 Let $K \subseteq H $ be subgroups of $G$ such that $[H:K] < \infty$. Let $\pi:C_c(\a / H) \to B(\mathscr{H})$ be a $^*$-representation. We define the map $\pi^K: C_c(\a / K) \to B(\mathscr{H} \otimes \ell^2(X^0 / K))$\index{piK@$\pi^K$} by
\begin{align}
\label{piK definition expression}
 \pi^K([a]_{xK})\,(\xi \otimes \delta_{uK}) := \begin{cases}
 \pi([a]_{xH})\xi \otimes \delta_{\mathbf{r}(x)K}\,, \qquad\; \text{if}\;\; uK = \mathbf{s}(x)K \\
 0, \qquad\qquad\qquad\qquad\quad\;\; \text{otherwise.}\;\; 
\end{cases}
\end{align}\\

\end{df}

\begin{prop}
 The map $\pi^K$ is a well-defined $^*$-representation.\\
\end{prop}

{\bf \emph{Proof:}} It is clear that the expression that defines $\pi^{K}([a]_{xH})$ defines a linear operator in the inner product space $\mathscr{H} \otimes C_c(X^0 / K)$, which is easily observed to be bounded. Thus, $\pi^{K}([a]_{xH}) \in B(\mathscr{H} \otimes \ell^2(X^0 / K))$.

It is clear that expression (\ref{piK definition expression}) defines a linear mapping $\pi^K$ on $C_c(\a / K)$, so that we only need to see that it preserves products and the involution. To see that it preserves products, consider two elements of the form $[a]_{xK}$ and $[b]_{yK}$. There are two cases to consider: either $\mathbf{r}(y) \in \mathbf{s}(x) K$ or $\mathbf{r}(y) \notin \mathbf{s}(x) K$.

In the second case, we have $[a]_{xK}[b]_{yK} = 0$ and thus $\pi^K([a]_{xK}[b]_{yK}) = 0$. But also $\pi^K([a]_{xK})\pi^K([b]_{yK}) = 0$, because for any vector $\xi \otimes \delta_{uK}$ we have that $\pi^K([b]_{yK})(\xi \otimes \delta_{uK})$ is either zero or equal to $\pi([b]_{yK})\xi \otimes \delta_{\mathbf{r}(y)K}$, and therefore we always have $\pi^K([a]_{xK})\pi^K([b]_{yK})(\xi \otimes \delta_{uK}) = 0$.

In the first case we have
\begin{eqnarray*}
& & \pi^K([a]_{xK}[b]_{yK})\,(\xi \otimes \delta_{uK}) \;\; = \;\;\\
 & = & \pi^K([\alpha_{\widetilde{k}^{-1}}(a)b]_{x\widetilde{k}yK})\,(\xi \otimes \delta_{uK})\\
& = & \begin{cases}
 \pi([\alpha_{\widetilde{k}^{-1}}(a)b]_{x\widetilde{k}yH})\xi \otimes \delta_{\mathbf{r}(x\widetilde{k}y)K}\,, \qquad\; \text{if}\;\; uK = \mathbf{s}(x\widetilde{k}y)K \\
 0, \qquad\qquad\qquad\qquad\qquad\qquad\qquad\quad \text{otherwise.}\;\; 
\end{cases}\\
& = & \begin{cases}
 \pi( [a]_{xH})\pi([b]_{yH})\xi \otimes \delta_{\mathbf{r}(x)K}\,, \qquad\; \text{if}\;\; uK = \mathbf{s}(y)K \\
 0, \qquad\qquad\qquad\qquad\qquad\qquad\quad\;\; \text{otherwise.}\;\; 
\end{cases}\\
& = & \begin{cases}
 \pi^K([a]_{xK})\,(\pi([b]_{yH})\xi \otimes \delta_{\mathbf{s}(x)K})\,, \quad\; \text{if}\;\; uK = \mathbf{s}(y)K \\
 0, \qquad\qquad\qquad\qquad\qquad\qquad\qquad \text{otherwise.}\;\; 
\end{cases}\\
& = & \begin{cases}
 \pi^K([a]_{xK})\,(\pi([b]_{yH})\xi \otimes \delta_{\mathbf{r}(y)K})\,, \quad\; \text{if}\;\; uK = \mathbf{s}(y)K \\
 0, \qquad\qquad\qquad\qquad\qquad\qquad\qquad \text{otherwise.}\;\; 
\end{cases}\\
& = &  \pi^K([a]_{xK})\pi^K([b]_{yK})\,(\xi \otimes \delta_{uK})\,.
\end{eqnarray*}

In both cases we have $\pi^K([a]_{xK}[b]_{yK}) = \pi^K([a]_{xK})\pi^K([b]_{yK})$, hence $\pi^K$ preserves products. Let us now check that it preserves the involution. We have

\begin{eqnarray*}
& & \langle \pi^K([a]_{xK})\,(\xi \otimes \delta_{uK})\,,\, \eta \otimes \delta_{vK} \rangle \;\; =\\
& = & \begin{cases}
\langle \pi([a]_{xH})\xi \otimes \delta_{\mathbf{r}(x)K}\,,\, \eta \otimes \delta_{vK} \rangle  \,, \qquad\; \text{if}\;\; uK = \mathbf{s}(x)K \\
 0, \qquad\qquad\qquad\qquad\qquad\qquad\qquad \text{otherwise.}\;\; 
\end{cases}\\
& = & \begin{cases}
\langle \pi([a]_{xH})\xi \,,\, \eta  \rangle  \,, \qquad\; \text{if}\;\; uK = \mathbf{s}(x)K \;\;\text{and}\;\; vK = \mathbf{r}(x)K\\
 0, \qquad\qquad\qquad\qquad \text{otherwise.}\;\; 
\end{cases}\\
& = & \begin{cases}
\langle \xi \,,\, \pi([a^*]_{x^{-1}H})\eta  \rangle  \,, \qquad\; \text{if}\;\; uK = \mathbf{s}(x)K \;\;\text{and}\;\; vK = \mathbf{r}(x)K\\
 0, \qquad\qquad\qquad\qquad\quad \text{otherwise.}\;\; 
\end{cases}\\
& = & \begin{cases}
\langle \xi \otimes \delta_{uK}\,,\, \pi([a^*]_{x^{-1}H})\eta \otimes \delta_{\mathbf{s}(x)K} \rangle  \,, \qquad\; \text{if}\;\;  vK = \mathbf{r}(x)K\\
 0, \qquad\qquad\qquad\qquad\qquad\qquad\qquad\qquad \text{otherwise.}\;\; 
\end{cases}\\
& = & \langle \xi \otimes \delta_{uK} \,,\, \pi^K([a^*]_{x^{-1}K})\,(\eta \otimes \delta_{vK}) \rangle\,.
\end{eqnarray*}
Hence, we conclude that $\pi^K([a]_{xK})^* = \pi^K(([a]_{xK})^*)$, and therefore $\pi^K$ preserves the involution. Hence, $\pi^K$ is a $^*$-representation. \qed\\

\begin{lemma}
\label{piK applied to Cc(A H)}
 Let us denote by $\delta_{uH} \in \ell^2(X^0 / K)$ the vector
\begin{align}
\label{delta uH inside L2 XOK}
 \delta_{uH} := \sum_{[h] \in \mathcal{S}_u \backslash H / K} \delta_{uh K}\,.
\end{align}
The map $\pi^K$ satisfies
\begin{align*}
 \pi^K([a]_{xH})\,(\xi \otimes \delta_{uH}):= \begin{cases}
 \pi([a]_{xH})\xi \otimes \delta_{\mathbf{r}(x)H}\,, \qquad\; \text{if}\;\; uH = \mathbf{s}(x)H, \\
 0, \qquad\qquad\qquad\qquad\quad\;\; \text{otherwise.}
\end{cases}\\
\end{align*}

\end{lemma}

{\bf \emph{Proof:}} We have
\begin{eqnarray*}
 \pi^K([a]_{xH})\,(\xi \otimes \delta_{uH}) & = & \sum_{[h] \in \mathcal{S}_x \backslash H / K} \sum_{[h']  \in \mathcal{S}_u \backslash H / K} \pi^K([\alpha_{h^{-1}}(a)]_{xhK}) \,(\xi \otimes \delta_{uh'K})\,,
\end{eqnarray*}
from which we see that, if $uH \neq \mathbf{s}(x)H$ then $\pi^K([a]_{xH})\,(\xi \otimes \delta_{uH}) = 0$. On the other hand, if $uH = \mathbf{s}(x)H$, then we have
\begin{eqnarray*}
 & & \pi^K([a]_{xH})\,(\xi \otimes \delta_{\mathbf{s}(x)H}) \;\;=\;\;\\
 & = & \sum_{[h] \in \mathcal{S}_x \backslash H / K} \sum_{[h']  \in \mathcal{S}_{\mathbf{s}(x)} \backslash H / K} \pi^K([\alpha_{h^{-1}}(a)]_{xhK}) \,(\xi \otimes \delta_{\mathbf{s}(x)h'K})\\
& = & \sum_{[h] \in \mathcal{S}_{\mathbf{s}(x)} \backslash H / K} \sum_{[h']  \in \mathcal{S}_{\mathbf{s}(x)} \backslash H / K} \pi^K([\alpha_{h^{-1}}(a)]_{xhK}) \,(\xi \otimes \delta_{\mathbf{s}(x)h'K})\\
& = & \sum_{[h] \in \mathcal{S}_{\mathbf{s}(x)} \backslash H / K}  \pi([\alpha_{h^{-1}}(a)]_{xhH}) \xi \otimes \delta_{\mathbf{r}(x)hK}\\
& = & \sum_{[h] \in \mathcal{S}_{\mathbf{s}(x)} \backslash H / K}  \pi([a]_{xH}) \xi \otimes \delta_{\mathbf{r}(x)hK}\\
& = & \sum_{[h] \in \mathcal{S}_{\mathbf{r}(x)} \backslash H / K}  \pi([a]_{xH}) \xi \otimes \delta_{\mathbf{r}(x)hK}\\
& = & \pi([a]_{xH}) \xi \otimes \delta_{\mathbf{r}(x)H}\,.
\end{eqnarray*}

This finishes the proof. \qed\\

{\bf \emph{Proof of Theorem \ref{embedding CcAH into CcAK passes to Calgebra}:}} In order to prove this statement we have to show that for any $f \in C_c(\a / H)$ we have $\|f\|_{C^*(\a / K)} = \|f\|_{C^*(\a / H)}$. Since we are viewing $C_c(\a / H)$ as a $^*$-subalgebra of $C_c(\a / K)$ we automatically have the inequality
\begin{align*}
  \|f\|_{C^*(\a / K)} \leq \|f\|_{C^*(\a / H)}\,.
\end{align*}
In order to prove the converse inequality, it suffices to prove that
\begin{align}
\label{ineq pi and piK}
 \|\pi(f)\| \leq \| \pi^K(f)\|\,,
\end{align}
for any nondegenerate $^*$-representation $\pi$ of $C_c(\a / H)$, because, since $\pi$ is arbitrary, this clearly implies that $\|f\|_{C^*(\a / H)} \leq \|f\|_{C^*(\a / K)}$. Let us then prove inequality (\ref{ineq pi and piK}).

 We can write any element $f \in C_c(\a / H)$ as $f = \sum_{xH \in X / H} [f(x)]_{xH}$. Furthermore we can split this sum according to the ranges of elements, i.e.
\begin{align*}
 f = \sum_{xH \in X / H} [f(x)]_{xH} \;=\; \sum_{vH \in X^0 /H}  \sum_{\substack{xH \in X / H \\ \mathbf{r}(x)H = vH}} [f(x)]_{xH}\,.
\end{align*}
 Suppose $\pi: C_c(\a / H) \to B(\mathscr{H})$ is a $^*$-representation and $\xi \in \mathscr{H}$ is a vector of norm one. We have
\begin{eqnarray*}
 \| \pi \Big( \sum_{xH \in X / H} [f(x)]_{xH} \;\Big)\xi \|^2 & = & \|  \sum_{vH \in X^0 /H} \pi \Big(  \sum_{\substack{xH \in X / H \\ \mathbf{r}(x)H = vH}} [f(x)]_{xH} \;\Big)\xi \|^2\,.
\end{eqnarray*}
For different units $vH \in X^0 / H$, the elements $\pi \big(  \sum_{\substack{xH \in X / H \\ \mathbf{r}(x)H = vH}} [f(x)]_{xH} \;\big)\xi$ are easily seen to be orthogonal, so that
\begin{eqnarray*}
& = &  \sum_{vH \in X^0 /H} \| \pi \Big(  \sum_{\substack{xH \in X / H \\ \mathbf{r}(x)H = vH}} [f(x)]_{xH} \;\Big)\xi \|^2\\
& = & \sum_{vH \in X^0 /H} \|   \sum_{\substack{xH \in X / H \\ \mathbf{r}(x)H = vH}} \pi([f(x)]_{xH}) \xi \|^2
\end{eqnarray*}
In the notation of (\ref{delta uH inside L2 XOK}), let $\delta_{uH}:= \sum_{[h] \in \mathcal{S}_u \backslash H / K} \delta_{uh K}$. Let us denote by $C_u$ the number of elements of $\mathcal{S}_u \backslash H / K$. It is not difficult to check that for any $r \in H$ the map $[h] \mapsto [r^{-1}h]$ is a well-defined bijection between $\mathcal{S}_u \backslash H / K$ and $\mathcal{S}_{ur} \backslash H / K$, so that $C_u = C_{ur}$. We have
\begin{eqnarray*}
 & = & \sum_{vH \in X^0 /H}  \frac{1}{C_v} \| \sum_{\substack{xH \in X / H \\ \mathbf{r}(x)H = vH}}\, \pi([f(x)]_{xH}) \xi \otimes \delta_{vH} \|^2\\
 & = & \| \sum_{vH \in X^0 /H}    \sum_{\substack{xH \in X / H \\ \mathbf{r}(x)H = vH}}\, \frac{1}{C_v} \pi([f(x)]_{xH}) \xi \otimes \delta_{vH} \|^2\\
 & = & \| \sum_{vH \in X^0 /H}    \sum_{\substack{xH \in X / H \\ \mathbf{r}(x)H = vH}}\, \frac{1}{C_{\mathbf{r}(x)}} \pi([f(x)]_{xH}) \xi \otimes \delta_{\mathbf{r}(x)H} \|^2\\
& = & \|  \sum_{xH \in X / H} \frac{1}{C_{\mathbf{r}(x)}}\, \pi([f(x)]_{xH}) \widetilde{\pi}(1_{\mathbf{s}(x)H})\xi \otimes \delta_{\mathbf{r}(x)H} \|^2\,.
\end{eqnarray*}
By Lemma \ref{piK applied to Cc(A H)} we have that 
\begin{eqnarray*}
& = & \| \sum_{xH \in X / H} \frac{1}{C_{\mathbf{r}(x)}}\, \pi^K([f(x)]_{xH})\,\big( \widetilde{\pi}(1_{\mathbf{s}(x)H})\xi \otimes \delta_{\mathbf{s}(x)H} \big) \|^2\,,
\end{eqnarray*}
and since  $\mathcal{S}_{\mathbf{s}(x)} \backslash H / K = \mathcal{S}_{\mathbf{r}(x)} \backslash H / K$, we get that $C_{\mathbf{s}(x)} = C_{\mathbf{r}(x)}$. Thus,
\begin{eqnarray*}
& = &  \| \sum_{xH \in X / H} \frac{1}{C_{\mathbf{s}(x)}}\, \pi^K([f(x)]_{xH})\,\big( \widetilde{\pi}(1_{\mathbf{s}(x)H})\xi \otimes \delta_{\mathbf{s}(x)H} \big) \|^2\\
& = &  \| \sum_{xH \in X / H}  \frac{1}{C_{\mathbf{s}(x)}}\, \pi^K([f(x)]_{xH})\,\big( \widetilde{\pi}(1_{\mathbf{s}(x)H})\xi \otimes \delta_{\mathbf{s}(x)H} \big) \|^2\\
& = & \| \sum_{xH \in X / H}  \pi^K([f(x)]_{xH})\,\Big( \frac{1}{C_{s(x)}}\, \widetilde{\pi}(1_{\mathbf{s}(x)H})\xi \otimes \delta_{\mathbf{s}(x)H} \Big) \|^2\,.
\end{eqnarray*}
Similarly as we did for ranges, we can split the sum $\sum_{xH \in X / H} [f(x)]_{xH}$ according to sources. In this way, since this sum is finite, there is a finite number of units $u_1H, \dots, u_nH \in X^0 /H$, which we assume to be pairwise different, such that we can write
\begin{align*}
 \sum_{xH \in X / H} [f(x)]_{xH} = \sum_{i = 1}^n \sum_{\substack{xH \in X / H \\ \mathbf{s}(x)H = u_iH}} [f(x)]_{xH}\,.
\end{align*}
By Lemma \ref{piK applied to Cc(A H)} we see that  $\pi^K([f(x)]_{xH}) \Big( \widetilde{\pi}(1_{u_iH})\xi \otimes \delta_{u_iH} \Big) = 0$ unless $\mathbf{s}(x)H = u_iH$. Hence we get
\begin{eqnarray*}
&  & \| \sum_{xH \in X / H}  \pi^K([f(x)]_{xH})\,\Big( \frac{1}{C_{s(x)}}\, \widetilde{\pi}(1_{\mathbf{s}(x)H})\xi \otimes \delta_{\mathbf{s}(x)H} \Big) \|^2\\
& = & \|  \sum_{xH \in X / H}  \pi^K([f(x)]_{xH})\,\big(  \sum_{i=1}^n \frac{1}{C_{u_i}}\, \widetilde{\pi}(1_{u_iH})\xi \otimes \delta_{u_iH} \big) \|^2\\
& = & \| \pi^K \Big( \sum_{xH \in X / H} [f(x)]_{xH} \Big)\big( \sum_{i=1}^n \frac{1}{C_{u_i}}\, \widetilde{\pi}(1_{u_iH})\xi \otimes \delta_{u_iH} \big) \|^2\,.
\end{eqnarray*}
We now notice that, since we are assuming $\xi$ to be of norm one, it follows that the vector
\begin{align*}
 \sum_{i=1}^n \frac{1}{C_{u_i}}\, \widetilde{\pi}(1_{u_iH})\xi \otimes \delta_{u_iH}\,,
\end{align*}
also has norm less or equal to one, because
\begin{eqnarray*}
 \| \sum_{i=1}^n \frac{1}{C_{u_i}}\, \widetilde{\pi}(1_{u_iH})\xi \otimes \delta_{u_iH}\|^2 & = &  \sum_{i=1}^n \| \frac{1}{C_{u_i}}\, \widetilde{\pi}(1_{u_iH})\xi \otimes \delta_{u_iH}\|^2\\
& = & \sum_{i=1}^n \| \widetilde{\pi}(1_{u_iH})\xi \|^2\\
& = & \| \widetilde{\pi} \Big(\sum_{i=1}^n  1_{u_iH} \Big)\xi \|^2\\
& \leq & \| \xi \|^2 \\
& = & 1\,.
\end{eqnarray*}

 Hence, taking the supremum over vectors $\xi$ of norm one, we immediately get the inequality
\begin{align*}
 \| \pi ( f ) \|\; \leq\; \| \pi^K ( f) \|\,.
\end{align*}
As we explained earlier, this proves that we get an embedding of $C^*(\a / H)$ into $C^*(\a / K)$. \qed\\

It follows from \ref{embedding CcAH into CcAK passes to Calgebra} that $\{C^*(\a / H)\}_{H \in \cl}$ is a direct system of $C^*$-algebras. Let us denote by $\mathcal{D}_{\max}(\a)$\index{Dmax(A)@$\mathcal{D}_{\max}(\a)$} its corresponding $C^*$-algebraic direct limit
\begin{align}
 \d_{\max}(\a) := \lim_{H \in \cl} C^*(\a / H)\,,
\end{align}
We notice that the algebra $\d(\a)$ is a dense $^*$-subalgebra of $\d_{\max}(\a)$. We now want to show that the action $\overline{\alpha}$ of $G$ on $\d(\a)$ extends to $\d_{\max}(\a)$.\\

\begin{thm}
\label{action extends to DmaxA}
 The action $\overline{\alpha}$ of $G$ on $\d(\a)$ extends uniquely to an action of $G$ on $\d_{\max}(\a)$ and is such that $\overline{\alpha}_g$ takes $C^*(\a / H)$ to $C^*(\a / g H g^{-1})$, for any $g \in G$.\\
\end{thm}

{\bf \emph{Proof:}} Since $\overline{\alpha}_g$ is a $^*$-isomorphism between $C_c(\a / H)$ and $C_c(\a/g Hg^{-1})$, it necessarily extends to a $^*$-isomorphism between the enveloping $C^*$-algebras $C^*(\a / H)$ and $C^*(\a / gHg^{-1})$. Since $C^*(\a / gHg^{-1})$ is embedded in $D_{\max}(\a)$, we can see $\overline{\alpha}_g$ as an injective $^*$-homomorphism from $C^*(\a / H)$ into $D_{\max}(\a)$.

A routine computation shows that the following diagram of canonical injections commutes:
 \begin{align*}
\xymatrix{ C^*(\a / H) \ar[r] \ar[dr]_{\overline{\alpha}_g} &  C^*(\a / K) \ar[d]^{\overline{\alpha}_g}\\
  & D_{\max}(\a)\,.} 
\end{align*}
Hence, we obtain an injective $^*$-homomorphism from $\d_{\max}(\a)$ to itself, which we still denote by $\overline{\alpha}_g$, and which extends the usual map $\overline{\alpha}_g$ from $D(\a)$ to itself. It is also clear that this map is surjective, and that for $g, h \in G$ we have $\overline{\alpha}_{gh} = \overline{\alpha}_g \circ \overline{\alpha}_h$, so that we get an action of $G$ on $D_{\max}(\a)$ which extends the usual action of $G$ on $D(\a)$. \qed\\

\section{Reduced $C^*$-crossed products}
\label{reduced crossed products chapter}

In this section we define reduced $C^*$-crossed products by Hecke pairs and study some of their properties. Since the algebra $C_c(\a /\gm)$ admits several possible  $C^*$-completions, we will be able to form several reduced $C^*$-crossed products, such as $C^*_r(\a / \gm) \times_{\alpha, r} G / \gm$ and $C^*(\a / \gm) \times_{\alpha, r} G / \gm$. As we shall see, many of the main properties of reduced $C^*$-crossed products by groups hold also in the Hecke pair case.

In subsection \ref{comparison with larsen laca neshveyev section} we also compare our construction of a reduced crossed product by a Hecke pair with that of Laca, Larsen and Neshveyev in \cite{phase}, and show that they agree whenever they are both definable.\\

\subsection{Regular representations}
\label{regular representations chapter}

In this subsection we introduce the notion of \emph{regular representations} in the context of crossed products by Hecke pairs. These are concrete $^*$-representations of $C_c(\a / \gm) \times^{alg}_{\alpha} G /\gm$ involving the regular representation of the Hecke algebra $\h(G, \gm)$ and are indispensable for defining reduced $C^*$-crossed products.

In the theory of crossed products by groups $A \times G$, regular representations are the integrated forms of certain covariant representations involving the regular representation of $G$. They are defined in the following way: one starts with a nondegenerate representation $\pi$ of $A$ on some Hilbert space $\mathscr{H}$ and constructs a new representation $\pi_{\alpha}$ of $A$ on the Hilbert $\mathscr{H} \otimes \ell^2(G)$, defined in an appropriate way, such that $\pi_{\alpha}$ together with the regular representation of $G$ form a covariant representation. Their integrated form is then called a \emph{regular representation}.

We are now going to make an analogous construction in the case of Hecke pairs. The main novelty here is that we have to start with a representation $\pi$ of $\d(\a)$, instead of $C_c(\a / \gm)$, so that we can construct the new representation $\pi_{\alpha}$ of $C_c(\a / \gm)$. This is because we need to take into account all algebras of the form $C_c(\a / H)$, where $H = g_1\gm g_1^{-1} \cap \dots \cap g_n\gm g_n^{-1}$ is a finite intersection of conjugates of $\gm$. Naturally, when $\gm$ is a normal subgroup, $\d(\a)$ is nothing but the algebra $C_c(\a / \gm)$ itself, so that we will recover the original definition of a regular representation for crossed products by groups.\\

\begin{df}
 Let $\pi: D(\a) \to B(\mathscr{H})$ be a nondegenerate $^*$-representation. We define the map $\pi_{\alpha}:C_c(\a / \gm) \to B(\mathscr{H} \otimes \ell^2(G / \gm))$\index{pialpha@$\pi_{\alpha}$} by
\begin{align*}
 \pi_{\alpha}(f)\;(\xi \otimes \delta_{h\gm}) := \pi(\overline{\alpha}_h(f))\xi \otimes \delta_{h\gm}\,.\\
\end{align*}

\end{df}

\begin{prop}
\label{pi alpha is a representation}
 Let $\pi: D(\a) \to B(\mathscr{H})$ be a nondegenerate $^*$-representation. Then, the map $\pi_{\alpha}$ is a nondegenerate $^*$-representation of $C_c(\a / \gm)$.\\
\end{prop}

\begin{lemma}
\label{restriction of pi to Cc(a H) is nondegenerate}
 Let $\pi:D(\a) \to B(\mathscr{H})$ be a nondegenerate $^*$-representation. Then the restriction of $\pi$ to $C_c(\a / H)$ is nondegenerate, for any $H \in \mathcal{C}$.\\
\end{lemma}

{\bf \emph{Proof:}} Let $\xi \in \mathscr{H}$ be such that $\pi(C_c(\a / H)) \xi = 0$. Take any $x \in X$, $a \in \a_x$ and $K \in \mathcal{C}$ such that $K \subseteq H$. We have that
\begin{eqnarray*}
 \|\pi([a]_{xK}) \xi\|^2 & = & \langle \pi([a^*a]_{\mathbf{s}(x)K})\xi\,,\, \xi \rangle\\
& = & \langle \pi([a^*]_{x^{-1}K} \cdot [a]_{xH}) \xi\,,\, \xi \rangle\\
& = & \langle \pi([a^*]_{x^{-1}K}) \pi([a]_{xH}) \xi\,,\, \xi \rangle\\
& = & 0\,.
\end{eqnarray*}
From this we conclude that $\pi(C_c(\a / K)) \xi = 0$, for any $K \in \mathcal{C}$ such that $K \subseteq H$. Since for any subgroup $L \in \mathcal{C}$ we have $C_c(\a / L) \subseteq C_c(\a / (L \cap H))$, and obviously $L \cap H \subseteq H$, we can in fact conclude that $\pi(C_c(\a / L) )\xi = 0$ for all $L \in \mathcal{C}$. In other words, we have proven that $\pi(D(\a)) \xi = 0$, which by nondegeneracy of $\pi$ implies that $\xi = 0$. \qed\\

{\bf \emph{Proof of Proposition \ref{pi alpha is a representation}:}} It is clear that the expression that defines $\pi_{\alpha}(f)$, for $f \in C_c(\a / \gm)$, defines a linear operator on the inner product space $\mathscr{H} \otimes C_c(G / \gm)$. Let us first check that this operator is indeed bounded. We have
\begin{eqnarray*}
 \|\pi_{\alpha}(f)\;\Big(\sum_{[h] \in G / \gm} \xi_{h\gm} \otimes \delta_{h\gm} \Big)\|^2 & = & \| \sum_{[h] \in G / \gm} \pi(\overline{\alpha}_{h}(f))\xi_{h\gm} \otimes \delta_{h\gm}\|^2\\
& = &  \sum_{[h] \in G / \gm} \| \pi(\overline{\alpha}_{h}(f))\xi_{h\gm}\|^2\\
& \leq & \sum_{[h] \in G / \gm} \| \pi(\overline{\alpha}_{h}(f))\|^2 \|\xi_{h\gm}\|^2\\
& \leq & \sum_{[h] \in G / \gm} \| \overline{\alpha}_{h}(f)\|_{C^*(\a / h\gm h^{-1})}^2 \|\xi_{h\gm}\|^2\,.
\end{eqnarray*}
Since $\overline{\alpha}_{h}$ gives an isomorphism between $C^*(\a / \gm)$ and $C^*(\a / h \gm h^{-1})$ we get
\begin{eqnarray*}
 & = & \sum_{[h] \in G / \gm} \| f\|_{C^*(\a / \gm)}^2 \|\xi_{h\gm}\|^2\\
& = & \| f\|_{C^*(\a / \gm)}^2\, \|\sum_{[h] \in G / \gm} \xi_{h \gm} \otimes \delta_{h\gm}\|^2\,.
\end{eqnarray*}
Hence, $\pi_{\alpha}(f)$ is bounded and thus defines uniquely an operator in $B(\mathscr{H} \otimes \ell^2(G / \gm))$. It is simple to check that $\pi$ is linear and preserves products. Let us then see that it preserves the involution. We have
\begin{eqnarray*}
 \langle \pi_{\alpha}(f)\,(\xi \otimes \delta_{h\gm})\,,\, \eta \otimes \delta_{g\gm} \rangle & = & \langle \pi(\overline{\alpha}_h(f))\xi \otimes \delta_{h\gm}\,,\, \eta \otimes \delta_{g\gm} \rangle\\
& = & \langle \pi(\overline{\alpha}_h(f))\xi \,,\, \eta \rangle \langle \delta_{h\gm}\,,\, \delta_{g\gm} \rangle\\
& = & \langle \xi \,,\, \pi(\overline{\alpha}_h(f^*))\eta \rangle \langle \delta_{h\gm}\,,\, \delta_{g\gm} \rangle\\
& = & \langle \xi \,,\, \pi(\overline{\alpha}_g(f^*))\eta \rangle \langle \delta_{h\gm}\,,\, \delta_{g\gm} \rangle\\
& = & \langle \xi \otimes \delta_{h\gm}\,,\, \pi_{\alpha}(f^*)\, (\eta \otimes \delta_{g\gm}) \rangle\,.
\end{eqnarray*}
Thus, $\pi_{\alpha}(f)^*=\pi_{\alpha}(f^*)$, and therefore $\pi_{\alpha}$ defines a $^*$-representation. It remains to check that this $^*$-representation is nondegenerate. To see this, we start by canonically identifying $\mathscr{H} \otimes \ell^2(G / \gm)$ with the Hilbert space $\ell^2(G / \gm, \mathscr{H})$. On this Hilbert space, it is easy to see that $\pi_{\alpha}(f)$ is given by
\begin{align*}
 [\pi_{\alpha}(f)\,(\zeta)]\,(h \gm) = \pi(\overline{\alpha}_{h}(f))\,\zeta(h \gm)
\end{align*}
for $\zeta \in \ell^2(G / \gm, \mathscr{H})$. Suppose now that $\zeta \in \ell^2(G / \gm, \mathscr{H})$ is such that $\pi_{\alpha}(f)\, \zeta = 0$ for all $f \in C_c(\a / \gm)$. Thus, for each $h \gm \in G / \gm$ we have $\pi(\overline{\alpha}_{h}(f))\,\zeta(h \gm) = 0$ for all $f \in C_c(\a / \gm)$. This can be expressed equivalently as $\pi(f)\,\zeta(h \gm) = 0$ for all $f \in C_c(\a / h\gm h^{-1})$. By Lemma \ref{restriction of pi to Cc(a H) is nondegenerate} the restriction of $\pi$ to $C_c(\a / h \gm h^{-1})$ is nondegenerate and therefore we have $\zeta(h\gm) = 0$. Thus, $\pi_{\alpha}$ is nondegenerate. \qed\\

\begin{df}
 Let $\pi: D(\a) \to B(\mathscr{H})$ be a nondegenerate $^*$-representation and $\rho:\h(G, \gm) \to B(\ell^2(G / \gm))$ the right regular representation of the Hecke algebra. The pair $(\pi_{\alpha}, 1 \otimes \rho)$ is called a \emph{regular covariant representation}.\\
\end{df}

\begin{rem}
We observe that when $\gm$ is a normal subgroup of $G$ we have $g \gm g^{-1} = \gm $ for all $g \in G$, so that the algebra $D(\a)$ coincides with $C_c(\a / \gm)$. For this reason our notion of a regular representation coincides with the usual notion of a regular covariant representation of the system $(C_c(\a / \gm), G / \gm, \overline{\alpha})$.\\
\end{rem}

\begin{thm}
\label{theorem about regular representations formula etc}
 Every regular covariant representation $(\pi_{\alpha}, 1 \otimes \rho)$ is a covariant $^*$-representation. Moreover, its integrated form is given by
\begin{align}
\label{expression for int form of reg rep}
[\pi_{\alpha} \times (1 \otimes \rho)](f)\;(\xi \otimes \delta_{h\gm}) = \sum_{[g] \in G / \gm} \Delta(g^{-1}h)^{\frac{1}{2}}\, \pi \big(\overline{\alpha}_{g}(f(g^{-1}h\gm)) \big) \xi \otimes \delta_{g \gm}\,,
\end{align}
for every $f \in C_c(\a  / \gm) \times_{\alpha}^{alg} G / \gm$.\\
\end{thm}

{\bf \emph{Proof:}} We shall first check that the expression (\ref{expression for int form of reg rep}) does indeed define a $^*$-representation of $C_c(\a  / \gm) \times_{\alpha}^{alg} G / \gm$. Afterwards we will show that the covariant pre-$^*$-representation associated to it is precisely $(\pi_{\alpha}, 1 \otimes \rho)$.

Let $\pi_{reg}:C_c(\a  / \gm) \times_{\alpha}^{alg} G / \gm \to B(\mathscr{H} \otimes \ell^2(G / \gm))$ be defined by
\begin{align*}
\pi_{reg}(f)\;(\xi \otimes \delta_{h\gm}) := \sum_{[g] \in G / \gm} \Delta(g^{-1}h)^{\frac{1}{2}}\, \pi \big(\overline{\alpha}_{g}(f(g^{-1}h\gm)) \big) \xi \otimes \delta_{g \gm}\,.
\end{align*}
It is not evident that $\pi_{reg}$ is a bounded operator for all $f  \in C_c(\a  / \gm) \times_{\alpha}^{alg} G / \gm$, but it is clear that $\pi_{reg}(f)$ is well-defined as a linear operator on the inner product space $\mathscr{H} \otimes C_c(G / \gm)$. Under the identification of $\mathscr{H} \otimes C_c(G / \gm)$ with $C_c(G / \gm, \mathscr{H})$, it is easy to see that $\pi_{reg}(f)$ is given by
\begin{align*}
 [\pi_{reg}(f)\,\eta] \,(g\gm) = \sum_{[h]\in G / \gm} \Delta(g^{-1}h)^{\frac{1}{2}}\, \pi \big(\overline{\alpha}_{g}(f(g^{-1}h \gm))\big) \eta(h\gm)\,,
\end{align*}
for any $\eta \in C_c(G / \gm, \mathscr{H})$. Let us now check that $\pi_{reg}(f)$ is indeed bounded. For any vector $\eta \in C_c(G / \gm, \mathscr{H})$ we have
\begin{eqnarray*}
 \| \pi_{reg} (f)\, \eta \|^2 & = & \sum_{[g] \in G / \gm} \|[\pi_{reg}(f)\, \eta](g\gm)\|^2\\
& = & \sum_{[g] \in G / \gm} \| \sum_{[h]\in G / \gm} \Delta(g^{-1}h)^{\frac{1}{2}}\, \pi \big(\overline{\alpha}_{g}(f(g^{-1}h \gm))\big) \eta(h\gm) \;\,\|^2\\
& \leq & \sum_{[g] \in G / \gm} \Big( \sum_{[h]\in G / \gm} \Delta(g^{-1}h)^{\frac{1}{2}}\, \|\pi \big(\overline{\alpha}_{g}(f(g^{-1}h \gm))\big)\| \|\eta(h\gm)\| \Big)^2\,.
\end{eqnarray*}
For each $h\gm \in G / \gm$ let us define $T^{h\gm} \in C_c(G / \gm)$ by
\begin{align*}
 T^{h\gm}(g\gm) := \Delta(g^{-1}h)^{\frac{1}{2}}\, \|\pi \big(\overline{\alpha}_{g}(f(g^{-1}h \gm))\big)\| \|\eta(h\gm)\|\,,
\end{align*}
and $T \in C_c(G / \gm)$ by $T := \sum_{[h] \in G/ \gm} T^{h\gm}$, which is clearly a finite sum since $\eta$ has finite support. Thus, we have
\begin{eqnarray*}
 \| \pi_{reg} (f)\, \eta \|^2 & \leq & \sum_{[g] \in G / \gm} \Big( \sum_{[h]\in G / \gm} T^{h \gm} (g \gm) \Big)^2\\
& = & \sum_{[g] \in G / \gm}  (T (g \gm))^2\\
& = & \| T \|^2_{\ell^2(G/\gm)}\\
& = & \| \sum_{[h] \in G / \gm} T^{h\gm}\; \|^2_{\ell^2(G/\gm)}\\
& \leq & \Big( \sum_{[h] \in G / \gm} \|T^{h\gm} \|_{\ell^2(G/\gm)} \Big)^2\\
& = &  \bigg( \sum_{[h] \in G / \gm} \sqrt{\sum_{[g] \in G / \gm} \Delta(g^{-1}h)\, \|\pi \big(\overline{\alpha}_{g}(f(g^{-1}h \gm))\big)\|^2 \|\eta(h\gm)\|^2 }\; \bigg)^2\\
& = &  \bigg( \sum_{[h] \in G / \gm} \|\eta(h\gm)\|\; \sqrt{\sum_{[g] \in G / \gm} \Delta(g^{-1}h)\, \|\pi \big(\overline{\alpha}_{g}(f(g^{-1}h \gm))\big)\|^2 }\; \bigg)^2\,.
\end{eqnarray*}
By the Cauchy-Schwarz inequality in $\ell^2(G /\gm)$ we get
\begin{eqnarray*}
& \leq &  \Big( \sum_{[h] \in G / \gm} \|\eta(h\gm)\|^2 \Big) \; \Big(\sum_{[h]\in G /\gm} \sum_{[g] \in G / \gm} \Delta(g^{-1}h)\, \|\pi \big(\overline{\alpha}_{g}(f(g^{-1}h \gm))\big)\|^2 \; \Big)\\
& = & \Big(\sum_{[h]\in G /\gm} \sum_{[g] \in G / \gm} \Delta(g^{-1}h)^{\frac{1}{2}}\, \|\pi \big(\overline{\alpha}_{g}(f(g^{-1}h \gm))\big)\|^2 \; \Big) \;\| \eta \|^2\,,
\end{eqnarray*}
which shows that $\pi_{reg}(f)$ is bounded.

Let us now check that $\pi_{reg}$ preserves products and the involution. Let $f_1, f_2 \in C_c(\a / \gm) \times_{\alpha} G / \gm$. We have
\begin{eqnarray*}
& & \pi_{reg}(f_1 * f_2 ) \, \big(\xi \otimes \delta_{h\gm} \big) =\\
 & = &\!\!\!\!\! \sum_{[g] \in G / \gm} \Delta(g^{-1}h)^{\frac{1}{2}}\, \pi \big(\overline{\alpha}_{g} \big((f_1* f_2)(g^{-1}h \gm) \big)\big) \xi \otimes \delta_{g \gm}\\
& = &\!\!\!\!\! \sum_{[g] \in G / \gm} \sum_{[s] \in G / \gm} \Delta(g^{-1}h)^{\frac{1}{2}}\, \pi \big(\overline{\alpha}_{g} \big((f_1(s \gm) \overline{\alpha}_s(f_2(s^{-1}g^{-1}h\gm)) \big)\big) \xi \otimes \delta_{g \gm}\\
& = &\!\!\!\!\! \sum_{[g] \in G / \gm} \sum_{[s] \in G / \gm} \Delta(g^{-1}h)^{\frac{1}{2}}\, \pi \big(\overline{\alpha}_{g} (f_1(s \gm)) \overline{\alpha}_{gs}(f_2(s^{-1}g^{-1}h\gm)) \big) \xi \otimes \delta_{g \gm}\\
& = & \!\!\!\!\!\sum_{[g] \in G / \gm} \sum_{[s] \in G / \gm} \Delta(g^{-1}h)^{\frac{1}{2}}\, \pi \big(\overline{\alpha}_{g} (f_1(g^{-1}s \gm)) \overline{\alpha}_s(f_2(s^{-1}h\gm)) \big) \xi \otimes \delta_{g \gm}\\
& = & \!\!\!\!\! \sum_{[s] \in G / \gm} \! \sum_{[g] \in G / \gm} \!\!\! \Delta(g^{-1}s)^{\frac{1}{2}}\Delta(s^{-1}h)^{\frac{1}{2}} \pi \big(\overline{\alpha}_{g} (f_1(g^{-1}s \gm)) \big) \pi \big( \overline{\alpha}_s(f_2(s^{-1}h\gm)) \big) \xi \otimes \delta_{g \gm}\\
& = &\!\!\!\!\! \sum_{[s] \in G / \gm}  \pi_{reg} (f_1)\, \Big(\Delta(s^{-1}h)^{\frac{1}{2}}\,\pi \big( \overline{\alpha}_s(f_2(s^{-1}h\gm)) \big) \xi \otimes \delta_{s \gm} \Big)\\
& = &\!\!\!\!\! \pi_{reg} (f_1) \pi_{reg} (f_2) \, \big(\xi \otimes \delta_{h \gm} \big)\,.
\end{eqnarray*}

Hence we conclude that $\pi_{reg}(f_1 * f_2 ) = \pi_{reg}(f_1)\pi_{reg}(f_2 )$. Let us now check that $\pi_{reg}$ preserves the involution. For $f \in C_c(\a / \gm) \times_{\alpha} G / \gm$ we have

\begin{eqnarray*}
& & \big\langle \pi_{reg} (f^*)\, \big(\xi \otimes \delta_{h\gm} \big)\, , \, \eta \otimes \delta_{s\gm} \big\rangle  =\\
& = & \sum_{[g] \in G / \gm} \Big\langle \Delta(g^{-1}h)^{\frac{1}{2}}\, \pi \big(\overline{\alpha}_{g}(f^*(g^{-1}h \gm))\big) \xi \otimes \delta_{g \gm} \, , \, \eta \otimes \delta_{s\gm} \Big\rangle\\
& = & \sum_{[g] \in G / \gm} \Big\langle \Delta(g^{-1}h)^{\frac{1}{2}}\Delta(h^{-1}g)\, \pi \big(\overline{\alpha}_{g}( \overline{\alpha}_{g^{-1}h}(f(h^{-1}g \gm))^*)\big) \xi \otimes \delta_{g \gm} \, , \, \eta \otimes \delta_{s\gm} \Big\rangle\\
& = & \sum_{[g] \in G / \gm} \big\langle \Delta(h^{-1}g)^{\frac{1}{2}}\, \pi \big( \overline{\alpha}_h(f(h^{-1}g \gm))\big)^* \xi \, , \, \eta \big\rangle \, \big\langle \delta_{g \gm} \, , \, \delta_{s\gm} \big\rangle\\
& = &  \big\langle \xi \, , \, \Delta(h^{-1}s)^{\frac{1}{2}}\, \pi \big( \overline{\alpha}_h(f(h^{-1}s \gm))\big) \eta \big\rangle\,.
\end{eqnarray*}
On the other side we also have
\begin{eqnarray*}
& & \big\langle \xi \otimes \delta_{h\gm} \, , \,  \pi_{reg} (f)\, \big(\eta \otimes \delta_{s\gm} \big) \big\rangle  =\\
& = & \sum_{[g] \in G / \gm} \Big\langle  \xi \otimes \delta_{h \gm} \, , \, \Delta(g^{-1}s)^{\frac{1}{2}}\, \pi \big(\overline{\alpha}_{g}(f(g^{-1}s \gm))\big) \eta \otimes \delta_{g\gm} \Big\rangle\\
& = & \sum_{[g] \in G / \gm} \big\langle  \xi \, , \, \Delta(g^{-1}s)^{\frac{1}{2}}\, \pi \big(\overline{\alpha}_{g}(f(g^{-1}s \gm))\big) \eta \big\rangle \, \big\langle \delta_{h \gm} \, , \, \delta_{g\gm} \big\rangle\\
& = &  \big\langle \xi \, , \, \Delta(h^{-1}s)^{\frac{1}{2}}\, \pi \big( \overline{\alpha}_h(f(h^{-1}s \gm))\big) \eta \big\rangle\,.
\end{eqnarray*}

Therefore we can conclude that $\pi_{reg}(f^*) = \pi_{reg}(f)^*$. Hence, $\pi_{reg}$ is a $^*$-representation.

The restriction of $\pi_{reg}$ to $C_c(\a / \gm)$ is precisely $\pi_{\alpha}$, and since $\pi_{\alpha}$ is nondegenerate, then so is $\pi_{reg}$. Hence, it follows from \cite[Theorem 3.40]{palma crossed 11} that $\pi_{reg}$ is the integrated form of a covariant pre-$^*$-representation $(\pi_{reg}|, \omega_{\pi_{reg}})$, as defined in  \cite[Proposition 3.36]{palma crossed 11}. As we pointed out above, $\pi_{reg}| = \pi_{\alpha}$. Thus, to finish the proof we only need to prove that $\omega_{\pi_{reg}} = 1 \otimes \rho$. For a vector of the form $\pi_{\alpha}([a]_{x\gm})(\xi \otimes \delta_{h \gm}) \in \pi_{\alpha}(C_c(\a / \gm)) ( \mathscr{H} \otimes \ell^2(G / \gm) )$ and a double coset $\gm g \gm$ we have
\begin{eqnarray*}
 \omega_{\pi_{reg}} (\gm g \gm) \, \pi_{\alpha}([a]_{x\gm}) (\xi \otimes \delta_{h\gm}) & = & \widetilde{\pi_{reg}}(\gm g \gm) \, \pi_{\alpha}([a]_{x\gm}) (\xi \otimes \delta_{h\gm})\\
& = & \widetilde{\pi_{reg}}(\gm g \gm)  \pi_{reg}([a]_{x\gm}) (\xi \otimes \delta_{h\gm})\\
& = & \pi_{reg}(\gm g \gm * [a]_{x\gm}) (\xi \otimes \delta_{h\gm})\,.
\end{eqnarray*}
Let us now compute $\pi_{reg}(f)(\xi \otimes \delta_{h\gm})$ for $f := \gm g \gm * [a]_{x\gm}$. By definition
\begin{align*}
 \pi_{reg}(f)\;(\xi \otimes \delta_{h\gm}) = \sum_{[s] \in G / \gm} \Delta(s^{-1}h)^{\frac{1}{2}}\, \pi \big(\overline{\alpha}_{s}(f(s^{-1}h\gm)) \big) \xi \otimes \delta_{s \gm}\,.
\end{align*}
It is clear that $f(s^{-1} h \gm)$ is nonzero if and only if $s^{-1}h \gm \subseteq \gm g \gm$, which is equivalent to $s\gm \subseteq h\gm g^{-1} \gm$. Hence,
\begin{eqnarray*}
  & = & \sum_{[s] \in\; h \gm g^{-1} \gm / \gm} \Delta(s^{-1}h)^{\frac{1}{2}}\, \pi \big(\overline{\alpha}_{s}(f(s^{-1}h \gm))\big) \xi \otimes \delta_{s \gm}\,.
\end{eqnarray*}
It is easy to see that $[\theta] \mapsto [h\theta g^{-1}]$ establishes a  well-defined bijection between $\gm / \gm^{g^{-1}}$  and $h \gm g^{-1} \gm / \gm$, so that
\begin{eqnarray*}
& = & \sum_{[\theta] \in \gm / \gm^{g^{-1}}} \Delta(g\theta^{-1}h^{-1}h)^{\frac{1}{2}}\, \pi \big(\overline{\alpha}_{h\theta g^{-1}}(f( g\theta^{-1}h^{-1}h \gm))\big) \xi \otimes \delta_{h \theta g^{-1} \gm}\\
& = & \sum_{[\theta] \in \gm / \gm^{g^{-1}}} \Delta(g)^{\frac{1}{2}}\, \pi \big(\overline{\alpha}_{h\theta g^{-1}}(f( g \gm))\big) \xi \otimes \delta_{h \theta g^{-1} \gm}\,.
\end{eqnarray*}
Now, it is easily seen that $f(g \gm) = [\alpha_g(a)]_{xg^{-1}g \gm g^{-1}}$. Hence, we get
\begin{eqnarray*}
 & = & \sum_{[\theta] \in \gm / \gm^{g^{-1}}} \Delta(g)^{\frac{1}{2}}\, \pi \big(\overline{\alpha}_{h\theta g^{-1}}([\alpha_g(a)]_{xg^{-1} g \gm g^{-1}})\big) \xi \otimes \delta_{h \theta g^{-1} \gm}\\
& = & \sum_{[\theta] \in \gm / \gm^{g^{-1}}} \Delta(g)^{\frac{1}{2}}\, \pi \big(\overline{\alpha}_h ([a]_{x \gm})\big) \xi \otimes \delta_{h \theta g^{-1} \gm}\\
& = & (1 \otimes \rho) (\gm g \gm)\, \Big( \pi \big(\overline{\alpha}_h ([a]_{x \gm})\big) \xi \otimes \delta_{h \gm} \Big)\\
& = & (1 \otimes \rho) (\gm g \gm)\, \pi_{\alpha} ([a]_{x \gm})\, \big(\xi \otimes \delta_{h \gm} \big)\,.
\end{eqnarray*}

This shows that $\omega_{\pi_{reg}} = 1 \otimes \rho$ in $\pi_{\alpha}(C_c(\a / \gm))( \mathscr{H} \otimes \ell^2(G / \gm) )$ and finishes the proof. \qed\\

\begin{rem}
 The proof of Theorem \ref{theorem about regular representations formula etc} may seem odd, since we did not first prove that the pair $(\pi_{\alpha}, 1 \otimes \rho)$ is a covariant $^*$-representation and then deduce that its integrated form $\pi_{\alpha} \times (1 \otimes \rho)$ is a $^*$-representation and is given by (\ref{expression for int form of reg rep}). Instead we followed the opposite approach. This is because we do not know how to prove directly that $(\pi_{\alpha}, 1 \otimes \rho)$ is a covariant $^*$-representation, due to the several difficult technicalities that arise in the computations.
\end{rem}

\subsection{Reduced $C^*$-crossed products}
\label{reduced C crossed products section}

We now want to define reduced $C^*$-norms in the $^*$-algebraic crossed product $C_c(\a / \gm) \times_{\alpha}^{alg} G / \gm $. Since $C_c(\a / \gm)$ admits several canonical $C^*$-completions one should expect that there are several reduced $C^*$-norms we can give to $C_c(\a / \gm) \times_{\alpha}^{alg} G / \gm $, which give rise to different reduced $C^*$-crossed products, as for example $C^*_r(\a / \gm) \times_{\alpha, r} G / \gm$ and $C^*(\a / \gm) \times_{\alpha, r} G / \gm$. We will treat in this section all these different reduced $C^*$-norms (and reduced $C^*$-crossed products) in a single approach, and for that the notion we need is that of a \emph{$\overline{\alpha}$-permissible} $C^*$-norm on $\d(\a)$:\\

\begin{df}
 A $C^*$-norm $\|\cdot \|_{\tau}$ in $\d(\a)$ is said to be \emph{$\overline{\alpha}$-permissible} if the action $\overline{\alpha}$ of $G$ on $\d(\a)$ extends to  $\d_{\tau}(\a)$\index{Dbeta(A)@$\d_{\tau}(\a)$}, the completion of $\d(\a)$ with respect to the norm $\| \cdot \|_{\tau}$. In other words, if for every $g \in G$ the automorphism $\overline{\alpha}_{g}$ of $\d(\a)$ is continuous with respect to $\| \cdot \|_{\tau}$.\\
\end{df}

\begin{df}
\label{norm beta r}
Let $\| \cdot \|_{\tau}$ be an $\overline{\alpha}$-permissible $C^*$-norm in $\d(\a)$ and let us denote by $\d_{\tau}(\a)$ and $C^*_{\tau}(\a / \gm)$\index{Cbeta(Agm)@$C^*_{\tau}(\a / \gm)$} the completions of $\d(\a)$ and $C_c(\a / \gm)$, respectively, with respect to the norm $\| \cdot \|_{\tau}$. We define the norm $\| \cdot \|_{\tau, r}$\index{1 norm beta r@$\Vert \cdot \Vert_{\tau,r}$} in $C_c(\a / \gm) \times_{\alpha}^{alg} G / \gm $ by
\begin{align*}
 \|f\|_{\tau,r} := \sup_{\pi \in R (\mathcal{D}_{\tau}(\a))} \| [\pi_{\alpha} \times (1 \otimes \rho)] (f)\|\,,
\end{align*}
where the supremum is taken over the class $R (\mathcal{D}_{\tau}(\a))$ of all nondegenerate $^*$-representations of $\mathcal{D}_{\tau}(\a )$. The completion of $C_c(\a / \gm) \times_{\alpha}^{alg} G / \gm $ with respect to this norm shall be denoted by $C_{\tau}^*(\a / \gm) \times_{\alpha, r} G / \gm $\index{crossed product (reduced)@$C_{\tau}^*(\a / \gm) \times_{\alpha, r} G / \gm $} and referred to as the \emph{reduced crossed product} of $C_{\tau}^*(\a / \gm)$ by the Hecke pair $(G, \gm)$.\\
\end{df}

Before we prove that $\| \cdot \|_{\tau}$ is indeed a $C^*$-norm, let us first look at the two main instances we have in mind, which arise when $C^*_{\tau}(\a / \gm)$ is $C^*_r(\a / \gm)$ or $C^*(\a / \gm)$. It is not obvious from the start that there exists a $C^*$-norm $\| \cdot \|_{\tau}$ in $\d(\a)$ whose restriction to $C_c(\a / \gm)$ will give the reduced or the maximal $C^*$-norm in $C_c(\a / \gm)$, but this is indeed the case from what we proved in the preliminary sections \ref{reduced completions of Cc(a H) section} and \ref{maximal completions of Cc(a H) section}:

\begin{itemize}
 \item \emph{For $C^*_r(\a / \gm)$:}

 As described in Section \ref{reduced completions of Cc(a H) section}, we can form the $C^*$-algebraic direct limit  $\d_r(\a) = \lim_{H \in \cl} C^*_r(\a / H)$, which contains $\d(\a)$ as a dense $^*$-subalgebra. Taking $\| \cdot \|_{\tau}$ to be the $C^*$-norm  $\| \cdot \|_r$ of $\d_r(\a)$, we see that $C^*_{\tau}(\a / \gm) = C^*_r(\a / \gm)$. The norm $\| \cdot \|_r$ is $\overline{\alpha}$-permissible because of Theorem \ref{action extends to DrA}.
 \item \emph{For $C^*(\a / \gm)$:} 

As described in Section \ref{maximal completions of Cc(a H) section}, we can form the $C^*$-algebraic direct limit  $\d_{\max}(\a) = \lim_{H \in \cl} C^*(\a / H)$, which contains $\d(\a)$ as a dense $^*$-subalgebra. Taking $\| \cdot \|_{\tau}$ to be the $C^*$-norm  $\| \cdot \|_{\max}$ of $\d_{\max}(\a)$, we see that $C^*_{\tau}(\a / \gm) = C^*(\a / \gm)$. The norm $\| \cdot \|_{\max}$ is $\overline{\alpha}$-permissible because of Theorem \ref{action extends to DmaxA}.\\
\end{itemize}

\begin{lemma}
\label{p alpha is continuous with repect to beta r norm}
 If $\pi:\d(\a) \to B(\mathscr{H})$ is a nondegenerate $^*$-representation which is continuous with respect to an $\overline{\alpha}$-permissible norm $\| \cdot \|_{\tau}$ in $\d(\a)$, then $\pi_{\alpha}$ is a representation of $C_c(\a / \gm)$ which is continuous with respect to the norm $\| \cdot \|_{\tau}$ as well.\\
\end{lemma}

{\bf \emph{Proof:}}  Let $f \in C_c(\a / \gm)$. We have
\begin{eqnarray*}
 \| \pi_{\alpha}(f) \big(\sum_{[h] \in G / \gm} \xi_{h \gm} \otimes \delta_{h\gm} \big)\|^2 & = & \|  \sum_{[h] \in G / \gm} \pi(\overline{\alpha}_h (f))\xi_{h \gm} \otimes \delta_{h\gm} \|^2\\
& = & \sum_{[h] \in G / \gm} \| \pi(\overline{\alpha}_h (f))\xi_{h \gm} \|^2\\
& \leq & \sum_{[h] \in G / \gm} \| \pi(\overline{\alpha}_h (f))\|^2\|\xi_{h \gm} \|^2\\
& \leq & \sum_{[h] \in G / \gm} \| \overline{\alpha}_h (f)\|_{\tau}^2\,\|\xi_{h \gm} \|^2\,.
\end{eqnarray*}

Since $\| \cdot \|_{\tau}$ is $\overline{\alpha}$-permissible we have that $\| \overline{\alpha}_h(f)\|_{\tau} = \|f \|_{\tau}$. Hence we have
\begin{eqnarray*}
\| \pi_{\alpha}(f) \big(\sum_{[h] \in G / \gm} \xi_{h \gm} \otimes \delta_{h\gm} \big)\|^2 & \leq & \sum_{[h] \in G / \gm} \| f\|_{\tau}^2\|\xi_{h \gm} \|^2\\
& = & \| f\|_{\tau}^2 \;\| \sum_{[h] \in G / \gm} \xi_{h \gm} \otimes \delta_{h \gm}\|^2\,.
\end{eqnarray*}
Hence, $\pi_{\alpha}$ is continuous with respect to the norm $\| \cdot \|_{\tau}$. \qed\\

\begin{prop} $\| \cdot \|_{\tau,r}$ is a well-defined $C^*$-norm on $C_c(\a / \gm) \times_{\alpha}^{alg} G / \gm $.\\
\end{prop}

{ \bf \emph{Proof:}} First we must show that the supremum in the definition of $\| \cdot \|_{\tau,r}$ is bounded. Given a $^*$-representation $\pi$ of $\d_{\tau}$ we have, by Lemma \ref{p alpha is continuous with repect to beta r norm}, that
\begin{eqnarray*}
&  &  \| [\pi_{\alpha} \times (1 \otimes \rho)] (f) \| \;\; \leq\;\;\\
 & \leq & \sum_{[g] \in \gm \backslash G / \gm}\; \sum_{x\gm^g \in X / \gm^g} \| \pi_{\alpha}\big( \big[f(g\gm) (x) \big]_{x \gm}\big) \| \| (1 \otimes \rho) \,(\gm g \gm) \| \| \widetilde{\pi_{\alpha}}(1_{xg\gm}) \| \\
& \leq & \sum_{[g] \in \gm \backslash G / \gm}\; \sum_{x\gm^g \in X / \gm^g} \| \big[f(g\gm) (x)\big]_{x \gm} \|_{\tau}\; \| \gm g \gm \|_{C^*_r(G, \gm)}\,.
\end{eqnarray*}
Thus, since $\| [\pi_{\alpha} \times (1 \otimes \rho)] (f) \|$ is finite and bounded by a number that does not depend on $\pi$, we conclude that $\| f\|_{\tau,r}$ is bounded by this same number.

It is clear from the definition and the above paragraph that $\| \cdot \|_{\tau,r}$ is $C^*$-seminorm. To prove that it is actually a $C^*$-norm it is enough to prove that if $\pi$ is a faithful nondegenerate $^*$-representation of $\d_{\tau}(\a)$, then $\pi_{\alpha} \times (1 \otimes \rho)$ is a faithful $^*$-representation of $C_c(\a / \gm) \times_{\alpha}^{alg} G / \gm $. Let us then prove this claim. Suppose $f \in C_c(\a / \gm) \times_{\alpha}^{alg} G / \gm$ is such that $[\pi_{\alpha} \times (1 \otimes \rho)](f) = 0$. Then, for every $\xi \otimes \delta_{h\gm} \in \mathscr{H} \otimes \ell^2(G / \gm)$ we have
\begin{align*}
 0\; = \;\; [\pi_{\alpha} \times (1 \otimes \rho)](f)(\xi \otimes \delta_{h \gm})\; = \;  \sum_{[g] \in G / \gm} \Delta(g^{-1}h)^{\frac{1}{2}}\, \pi \big(\overline{\alpha}_{g}(f(g^{-1}h \gm))\big) \xi \otimes \delta_{g \gm}\,.
\end{align*}
In particular, for $g\gm = \gm$, we have $\pi (f(h \gm)) \xi = 0$, and since this holds for every $\xi \in \mathscr{H}$ we have $\pi(f(h\gm))=0$. Now, since $\pi$ is a faithful $^*$-representation, it follows that $f(h\gm) =0$. Since this holds for every $h \gm \in G / \gm$, we have $f = 0$, i.e. $\pi_{\alpha} \times (1 \otimes \rho)$ is injective. \qed\\

The next result explains why we call the completion of $C_c(\a / \gm) \times_{\alpha}^{alg} G / \gm$ under the norm $\| \cdot \|_{\tau, r}$ the reduced crossed product of $C^*_{\tau}(\a / \gm)$ by the Hecke pair $(G, \gm)$ and justifies also the notation $ C^*_{\tau}(\a / \gm) \times_{\alpha, r} G / \gm $ chosen to denote this completion.\\

\begin{prop}
\label{restriction of rr norm to CcA gm}
The restriction of the norm $\| \cdot \|_{\tau,r}$ of $C_c(\a / \gm) \times_{\alpha}^{alg} G / \gm$ to $C_c(\a / \gm)$ is precisely the norm $\| \cdot \|_{\tau}$. Hence, the embedding $C_c(\a/ \gm) \to C_c(\a / \gm) \times_{\alpha}^{alg} G / \gm$ completes to an embedding $C_{\tau}^*(\a / \gm) \to C^*_{\tau}(\a / \gm) \times_{\alpha, r} G / \gm $.\\
\end{prop}

{\bf \emph{Proof:}} Let $\pi: \d_{\tau}(\a) \to B(\mathscr{H})$ be a nondegenerate $^*$-representation. From Lemma \ref{p alpha is continuous with repect to beta r norm} we have, for every $f \in C_c(\a / \gm)$,
\begin{align*}
 \| [\pi_{\alpha} \times (1 \otimes \rho)] (f) \| = \| \pi_{\alpha}(f) \| \leq \| f \|_{\tau}\,,
\end{align*}
 and therefore
\begin{align*}
 \| f \|_{\tau,r} \leq \| f \|_{\tau}\,.
\end{align*}
We now wish to prove the converse inequality. Let $\pi: \d_{\tau}(\a) \to B(\mathscr{H})$ be a faithful nondegenerate $^*$-representation. For any $f  \in C_c(\a / \gm)$ we have

\begin{eqnarray*}
 \| f \|_{\tau} & = &  \| \pi(f) \| \;\; = \;\;  \sup_{\|\xi\| = 1} \| \pi (f) \xi \|\\
 & = &  \sup_{\|\xi\| = 1} \| \pi_{\alpha} (f) (\xi \otimes \delta_{\gm}) \|\\
& \leq & \sup_{\|\zeta\| = 1} \| \pi_{\alpha} (f) \zeta \| \;\; = \;\; \| \pi_{\alpha}(f) \|\\
& = & \| [\pi_{\alpha} \times (1 \otimes \rho)](f) \| \;\; \leq \;\; \| f \|_{\tau, r}\,,
\end{eqnarray*}
thus proving the converse inequality. We conclude that
\begin{align*}
 \| f \|_{\tau,r} \;=\; \| f \|_{\tau}\,,
\end{align*}
for any $ f \in C_c(\a / \gm)$ and this finishes the proof. \qed\\

An important feature of reduced crossed products by groups $A \times_r G$ is the existence of a faithful conditional expectation onto $A$. We will now explain how this holds as well for reduced crossed products by Hecke pairs, with somewhat analogous proofs. The goal is to prove Theorem \ref{props of Eggm} bellow, and for that we follow closely the approach presented in \cite{phillips} in the case of groups.\\

\begin{prop}
\label{Eggm is continuous}
For every $g \gm \in G / \gm$ the map $E_{g \gm}$\index{Eggm@$E_{g\gm}$} defined by
\begin{align*}
 E_{g \gm}:C_c(\a / \gm) \times_{\alpha}^{alg} G / \gm & \longrightarrow C^*_{\tau}(\a / \gm^g)\\
E_{g \gm} (f) & := f(g\gm)\,.
\end{align*}
is linear and continuous with respect to the norm $\| \cdot \|_{\tau,r}$.\\
\end{prop}

Before we give a proof of the result above we need to set some notation. For each element $g \gm \in G / \gm$ we will denote by $\sigma_{g \gm}$\index{s1 igmaggm@$\sigma_{g\gm}$} the Hilbert space isometry $\sigma_{g \gm} : \mathscr{H} \to \mathscr{H} \otimes \ell^2(G / \gm)$ defined by
\begin{align}
 \sigma_{g \gm} (\xi) := \xi \otimes \delta_{g \gm}\,.
\end{align}\\

{\bf \emph{Proof of Proposition \ref{Eggm is continuous}:}} Let $\pi$ be a faithful nondegenerate $^*$-represen-tation of $\d_{\tau}(\a)$. It is easily seen that $\sigma_{\gm}^*\, [\pi \times (1 \otimes \rho)](f)\, \sigma_{g \gm} = \Delta(g)^{\frac{1}{2}}\, \pi(f(g\gm))$. Hence we have 
\begin{eqnarray*}
 \| E_{g\gm}(f) \|_{\tau}   & = & \| f(g \gm) \|_{\tau} \;\; = \;\;  \| \pi(f(g \gm)) \| \\
& = &\| \Delta(g^{-1})^{\frac{1}{2}}\, \sigma_{\gm}^*\, [\pi_{\alpha} \times (1 \otimes \rho)](f)\, \sigma_{g \gm}\|\\
& \leq & \Delta(g^{-1})^{\frac{1}{2}}\,\| [\pi_{\alpha} \times (1 \otimes \rho)](f) \| \\
& \leq & \Delta(g^{-1})^{\frac{1}{2}}\,\| f \|_{\tau,r}\,.
\end{eqnarray*}
This finishes the proof. \qed\\

We shall henceforward make no distinction of notation between the maps $E_{g \gm}$ defined on $C_c(\a / \gm) \times_{\alpha}^{alg} G / \gm$ and their extension to $C^*_{\tau}(\a / \gm) \times_{\alpha, r} G / \gm$.

The following result is of particular importance in theory of reduced $C^*$-crossed products. Analogously to the case of groups, it reveals two important features of reduced $C^*$-crossed products by Hecke pairs: the fact that every element of a reduced crossed product is uniquely described in terms of its coefficients (determined by the $E_{g \gm}$ maps); and the fact that $E_{\gm}$ is a faithful conditional expectation.\\

\begin{thm}
\label{props of Eggm}
 We have
\begin{itemize}
 \item[i)] If $f \in C^*_{\tau}(\a / \gm) \times_{\alpha, r} G / \gm$ and $E_{g \gm} (f) = 0$ for all $g \gm \in G / \gm$, then $f = 0$.
 \item[ii)] $E_{\gm}$ is a faithful conditional expectation of $C^*_{\tau}(\a / \gm) \times_{\alpha, r} G / \gm$ onto $C^*_{\tau}(\a / \gm)$.\\
\end{itemize}

\end{thm}

We start with the following auxiliary result:\\

\begin{lemma}
\label{expression with sigmas and Es lemma}
Let $\pi$ be a nondegenerate $^*$-representation of $\d_{\tau}(\a)$. For all $ f \in C_{\tau}^*(\a / \gm) \times_{\alpha, r} G / \gm$ we have
\begin{align}
\label{lemma expression about reg rep with sigmas}
 \sigma_{g \gm}^*\, [\pi \times (1 \otimes \rho)](f)\, \sigma_{h \gm} = \Delta(g^{-1} h)^{\frac{1}{2}}\, \pi(\overline{\alpha}_g(E_{g^{-1}h\gm}(f)))\,.
\end{align}\\

\end{lemma}

{\bf \emph{Proof:}} We notice that equality (\ref{lemma expression about reg rep with sigmas}) above holds for any $f \in C_c(\a / \gm) \times_{\alpha}^{alg} G / \gm$, following the definitions of the maps $E_{t \gm}$, $[\pi_{\alpha} \times (1 \otimes \rho)](f)$ and $\sigma_{t \gm}$, with $t \gm \in G / \gm$. By continuity, it follows readily that the equality must hold for every $f \in C_{\tau}^*(\a / \gm) \times_{\alpha, r} G / \gm$. \qed\\

{\bf \emph{Proof of Theorem \ref{props of Eggm}:}} $\;i)$ Let $f \in C_{\tau}^*(\a / \gm) \times_{\alpha, r} G / \gm$.  Suppose $E_{g \gm}(f) = 0$ for all $g \gm \in G / \gm$. Then, for any given nondegenerate $^*$-representation $\pi$ of $\d_{\tau}(\a)$ we have, by Lemma \ref{expression with sigmas and Es lemma}, that $\sigma_{g\gm}^*\,[\pi_{\alpha} \times (1 \otimes \rho)](f)\, \sigma_{h\gm}= 0$ for all $g\gm, h \gm \in G / \gm$. Hence, $[\pi_{\alpha} \times (1 \otimes \rho)](f)=0$. Since, this is true for any $\pi$, we must have $\| f \|_{\tau,r}=0$, i.e. $f = 0$.

$ii)$ Let us first prove that $E_{\gm}$ is a conditional expectation, i.e. $E_{\gm}$ is an idempotent, positive, $C^*_{\tau}(\a / \gm)$-linear map.

 If $f \in C_c(\a / \gm)$ then it is clear that $E_{\gm}(f) = f$. By continuity and Proposition \ref{restriction of rr norm to CcA gm} it follows that $E_{\gm}(f) = f$ for all $f \in C^*_{\tau}(\a / \gm)$. Thus, $E_{\gm}$ is idempotent. 

Suppose now that $f \in C_c(\a / \gm) \times_{\alpha}^{alg} G / \gm$. We have
\begin{eqnarray*}
 E_{\gm}(f^**f) & = & (f^**f)(\gm) \;\; = \;\; \sum_{[h] \in G / \gm} f^*(h\gm) \overline{\alpha}_h (f(h^{-1}\gm))\\
& = & \sum_{[h] \in G / \gm} \Delta(h^{-1})\, \overline{\alpha}_h( f(h^{-1}\gm))^* \overline{\alpha}_h (f(h^{-1}\gm)) \;\; \geq \;\; 0
\end{eqnarray*}
By continuity it follows that $E_{\gm}(f^**f) \geq 0$ for all $f \in C_{\tau}^*(\a / \gm) \times_{\alpha, r} G / \gm$, i.e. $E_{\gm}$ is positive. It remains to show that $E_{\gm}$ is $C^*_{\tau}(\a / \gm)$-linear. We recall that we see $C_c(\a / \gm)$ as a $^*$-subalgebra of $C_c(\a / \gm) \times_{\alpha}^{alg} G / \gm$ in the following way: an element $f \in C_c(\a / \gm)$ is identified with the element $F \in C_c(\a / \gm) \times_{\alpha}^{alg} G / \gm$ with support in $\gm$ and such that $F(\gm) = f$. For any $f \in C_c(\a / \gm)$ and $f_2 \in C_c(\a / \gm) \times_{\alpha}^{alg} G / \gm$ we have
\begin{eqnarray*}
 E_{\gm}(f * f_2 ) & = & (F* f_2) (\gm) \;\; = \;\; \sum_{[h] \in G / \gm} F(h\gm) \overline{\alpha}_h (f_2(h^{-1}\gm))\\
& = & F(\gm) f_2(\gm) \;\; = \;\; f \,E_{\gm}(f_2)\,,
\end{eqnarray*}
and similarly we get $E_{\gm}(f_2 * f) = E_{\gm}(f_2)\, f$. Once again by continuity we conclude that the same equalities hold for $f \in C^*_{\tau}(\a / \gm)$ and $f_2 \in C_{\tau}^*(\a / \gm) \times_{\alpha, r} G / \gm$. Thus, $E_{\gm}$ is a conditional expectation.

Let us now prove that $E_{\gm}$ is faithful. For any $f \in  C_c(\a / \gm) \times_{\alpha}^{alg} G / \gm$ we have (where the first equality was computed above):
\begin{eqnarray*}
 E_{\gm}(f^**f) & = & \sum_{[h] \in G / \gm} \Delta(h^{-1})\, \overline{\alpha}_h( f(h^{-1}\gm))^* \overline{\alpha}_h (f(h^{-1}\gm))\\
& = & \sum_{[h] \in G / \gm} \Delta(h^{-1})\, \overline{\alpha}_h( E_{h {-1}\gm}(f))^* \overline{\alpha}_h (E_{h {-1}\gm}(f))\,.
\end{eqnarray*}
Hence, we have $E_{\gm}(f^**f) \geq \Delta(h^{-1})\, \overline{\alpha}_h( E_{h {-1}\gm}(f))^* \overline{\alpha}_h (E_{h {-1}\gm}(f))$ for each $h \gm \in G / \gm$,. By continuity this inequality holds for every $f \in C_{\tau}^*(\a / \gm) \times_{\alpha, r} G / \gm$, and therefore if $f \in C_{\tau}^*(\a / \gm) \times_{\alpha, r} G / \gm$ is such that $E_{\gm}(f^**f) = 0$, then $E_{g \gm}(f) = 0$ for all $g \gm \in G / \gm$. Hence, by part $i)$, we conclude that $f = 0$. Thus, $E_{\gm}$ is faithful. \qed\\

The next result shows, like in crossed products by groups, that to define the norm $\| \cdot \|_{\tau,r}$ of the reduced crossed product $C_{\tau}^*(\a / \gm) \times_{\alpha, r} G / \gm$ we only need to start with a faithful nondegenerate $^*$-representation of $\d_{\tau}(\a)$, instead of taking the supremum over all nondegenerate $^*$-representations of $\d_{\tau}(\a)$.\\

\begin{thm}
\label{norm beta r only needs one faithful pi}
 Let $\pi: \d_{\tau}(\a) \to B(\mathscr{H})$ be a nondegenerate $^*$-representation. We have that
\begin{itemize}
 \item[i)]  If $\pi_{\alpha}:C^*_{\tau}(\a / \gm) \to B(\mathscr{H} \otimes \ell^2(G /\gm))$ is faithful, then $[\pi_{\alpha} \times (1 \otimes \rho)]$ is a faithful $^*$-representation of $C_{\tau}^*(\a / \gm) \times_{\alpha, r} G / \gm$. Consequently,
\begin{align*}
 \| f \|_{\tau,r} = \| [\pi_{\alpha} \times (1 \otimes \rho)](f) \|\,,
\end{align*}
for all $f \in C_{\tau}^*(\a / \gm) \times_{\alpha, r} G / \gm$.
 \item[ii)] If $\pi$ is faithful, then $\pi_{\alpha}$ is faithful.\\
\end{itemize}

\end{thm}

{\bf \emph{Proof:}} Let us prove $i)$ first. Suppose $\pi_{\alpha}$ is faithful as a $^*$-representation of $C^*_{\tau}(\a / \gm)$. Let $f \in C_{\tau}^*(\a / \gm) \times_{\alpha, r} G / \gm$ be such that $[\pi \times (1 \otimes \rho)](f) = 0$. Then, of course, $[\pi_{\alpha} \times (1 \otimes \rho)](f^**f) = 0$ and we have
\begin{eqnarray*}
0 & = & \sigma_{g \gm}^*\, [\pi_{\alpha} \times (1 \otimes \rho)](f^**f)\, \sigma_{g \gm} \;\; = \;\; \pi( \overline{\alpha}_g(E_{\gm}(f^**f)))\\
& = & \sigma_{g \gm}^*\, \pi_{\alpha}(E_{\gm}(f^**f))\, \sigma_{g \gm}\,.
\end{eqnarray*}
This implies that $\pi_{\alpha}(E_{\gm}(f^**f)) = 0$, i.e. $E_{\gm}(f^**f)=0$, and since $E_{\gm}$ is a faithful conditional expectation we have $f^**f = 0$, i.e. $f=0$. Thus, $\pi_{\alpha} \times (1 \otimes \rho)$ is faithful.

Let us now prove claim $ii)$. We know that $\pi_{\alpha}$, as a $^*$-representation of $C_c(\a / \gm)$, is given by
\begin{align*}
\pi_{\alpha}(f) \,(\xi \otimes \delta_{g \gm}) =  \pi(\overline{\alpha}_{g}(f)) \xi \otimes \delta_{g \gm}\,,
\end{align*}
By continuity the same expression holds for $f \in C^*_{\tau}(\a / \gm)$. Now suppose that $\pi_{\alpha}(f) = 0$ for some $f \in C^*_{\tau}(\a / \gm)$. Then, by the above expression, we have $\pi(f) = 0$. Since $\pi$ is faithful we must have $f = 0$. Thus, $\pi_{\alpha}$ is faithful. \qed\\

Another feature of reduced $C^*$-crossed products by groups $A \times_r G$ is the fact that the reduced $C^*$-algebra of the group is always canonically embedded in the multiplier algebra $M(A \times_r G)$. The same is true in the Hecke pair case as we now show:\\

\begin{prop}
 There is a unique embedding of the reduced Hecke $C^*$-algebra $C^*_r(G, \gm)$ into $M(C^*_{\tau}(\a / \gm) \times_{\alpha, r} G / \gm)$ extending the action of $\h(G, \gm)$ on $C_c(\a / \gm) \times_{\alpha}^{alg} G / \gm$.\\
\end{prop}

{ \bf \emph{Proof:}} Let us first see that the action of $\h(G, \gm)$ on $C_c(\a / \gm) \times_{\alpha}^{alg} G / \gm$ is continuous with respect to the norm $\| \cdot \|_{\tau,r}$, so that it extends uniquely to an action of $\h(G, \gm)$ on $C^*_{\tau}(\a / \gm) \times_{\alpha, r} G / \gm$.

Let $\pi$ be a faithful nondegenerate $^*$-representation of $\d_{\tau}(\a)$. From Theorem \ref{norm beta r only needs one faithful pi} we know that $\pi_{\alpha} \times (1 \otimes \rho)$ is also faithful. For $f_1 \in \h(G, \gm)$ and $f_2 \in C_c(\a / \gm) \times_{\alpha}^{alg} G / \gm$, we have
\begin{eqnarray*}
\| f_1 * f_2 \|_{\tau,r} & = & \| [\pi_{\alpha} \times (1 \otimes \rho)](f_1*f_2) \| \\
& \leq & \| (1 \otimes  \rho)(f_1) \| \|[\pi_{\alpha} \times (1 \otimes \rho)](f_2)\|\\
& = & \| \rho(f_1) \| \|f_2\|_{\tau,r}\,.
\end{eqnarray*}
Thus, the action of $\h(G, \gm)$ on $C_c(\a / \gm) \times_{\alpha}^{alg} G / \gm$ extends uniquely to an action on $C^*_{\tau}(\a / \gm) \times_{\alpha, r} G / \gm$, or in other words, we have an embedding of $\h(G, \gm)$ into $M(C^*_{\tau}(\a / \gm) \times_{\alpha, r} G / \gm)$. We now want to prove that this embedding extends to an embedding of $C^*_r(G,\gm)$ into the same multiplier algebra. For that it is enough to prove that, for any $f \in \h(G, \gm)$, we have
\begin{align*}
\|f \|_{M(C^*_{\tau}(\a / \gm) \times_{\alpha, r} G / \gm)} = \| f \|_{C^*_r(G, \gm)}\,.
\end{align*}
 Let $\widetilde{\pi_{\alpha} \times (1 \otimes \rho)}$ denote the extension of $\pi_{\alpha} \times (1 \otimes \rho)$ to $M(C^*_{\tau}(\a / \gm) \times_{\alpha, r} G / \gm)$, which is faithful since $\pi_{\alpha} \times (1 \otimes \rho)$ is faithful on $C^*_{\tau}(\a / \gm) \times_{\alpha, r} G / \gm$. We have that $\widetilde{\pi_{\alpha} \times (1 \otimes \rho)}$ and $(1 \otimes \rho)$ coincide in $\h(G, \gm)$ since they are given by the same expression on the dense subspace $[\pi_{\alpha} \times (1 \otimes \rho)](C_c(\a / \gm) \times_{\alpha}^{alg} G / \gm) \mathscr{H}$. Thus, we have
\begin{align*}
 \widetilde{[\pi_{\alpha} \times (1 \otimes \rho)]}(f) = (1 \otimes \rho) (f)\,,
\end{align*}
for any $ f \in \h(G, \gm)$. It then follows that
\begin{eqnarray*}
 \|f \|_{M(C^*_{\tau}(\a / \gm) \times_{\alpha, r} G / \gm)} & = & \| \widetilde{[\pi_{\alpha} \times (1 \otimes \rho)]}(f) \| \;\; = \;\; \| (1 \otimes \rho) (f)\|\\
& = & \| \rho(f) \| \;\; = \;\; \| f \|_{C^*_r(G, \gm)}\,.
\end{eqnarray*}
This finishes the proof. \qed\\

As it is known, reduced $C^*$-crossed products by discrete groups satisfy a universal property among all the $C^*$-completions of the $^*$-algebraic crossed product that have a certain conditional expectation. This universal property says that every such completion has a canonical surjective map onto the reduced $C^*$-crossed product. As a consequence, the reduced $C^*$-crossed product is the only $C^*$-completion of the $^*$-algebraic crossed product that has a certain faithful conditional expectation.

The next result explains how this holds in the Hecke pair case.\\

\begin{thm}
 Let $\| \cdot \|_{\tau}$ be an $\overline{\alpha}$-permissible $C^*$-norm on $D(\a)$ and $\| \cdot \|_{\omega}$ a $C^*$-norm on $C_c(\a/ \gm) \times_{\alpha}^{alg} G / \gm$ whose restriction to $C_c(\a / \gm)$ is just the norm $\| \cdot \|_{\tau}$. Let us denote by $C^*_{\tau}(\a / \gm) \times_{\alpha,  \omega} G / \gm$ the completion of $C_c(\a/ \gm) \times_{\alpha}^{alg} G / \gm$ under the norm $\| \cdot \|_{\omega}$.

If there exists a bounded linear map $F:C^*_{\tau}(\a / \gm) \times_{\alpha, \omega} G / \gm \to C^*_{\tau}(\a / \gm)$ such that
\begin{align*}
 F(f) = f(\gm)\,,
\end{align*}
for all $f \in C_c(\a/ \gm) \times_{\alpha}^{alg} G / \gm$, then:
\begin{itemize}
 \item[a)] There exists a surjective $^*$-homomorphism
\begin{align*}
 \Lambda: C^*_{\tau}(\a / \gm) \times_{\alpha,  \omega} G / \gm \to C^*_{\tau}(\a / \gm) \times_{\alpha, r} G / \gm\,,
\end{align*}
such that $\Lambda$ is the identity on $C_c(\a/ \gm) \times_{\alpha}^{alg} G / \gm$.
 \item[b)] $F$ is a conditional expectation.
 \item[c)] $F$ is faithful if and only if $\Lambda$ is an isomorphism.
\end{itemize}
\end{thm}

{\bf \emph{Proof:}} Let $X_0$ be the space $C_c(\a/ \gm) \times_{\alpha}^{alg} G / \gm$. It is easily seen that $X_0$ is a (right) inner product $C_c(\a / \gm)$-module, where $C_c(\a / \gm)$ acts on $X_0$ by right multiplication and the inner product is given by
\begin{align*}
\langle f_1\,,\, f_2 \rangle := (f_1^**f_2)(\gm)\,.
\end{align*}
Since for any $f \in X_0$ and $f_1 \in C_c(\a/ \gm)$ we have
\begin{eqnarray*}
 \|\langle f*f_1\,,\,f*f_1 \rangle\|_{\tau} & = & \|((f*f_1)^**(f*f_1))(\gm) \|_{\tau} \\
 & = & \|(f_1^**f^**f*f_1)(\gm)\|_{\tau}\\
& = & \|f_1^*((f^**f)(\gm))f_1\|_{\tau} \\
& = &  \| f_1\|^2_{\tau} \|\langle f\,,\,f \rangle\|_{\tau}\,,
\end{eqnarray*}
it follows that we can complete $X_0$ to a (right) Hilbert $C^*_{\tau}(\a / \gm)$-module, which we will denote by $X$. The inner product on $X$, which extends the inner product $\langle \cdot \,,\, \cdot \rangle$ above, will be denoted by $\langle \cdot \,,\, \cdot \rangle_{C^*_{\tau}(\a / \gm)}$.

The $^*$-algebra $C_c(\a/ \gm) \times_{\alpha}^{alg} G / \gm$ acts on $X_0$ by left multiplication and therefore it is easily seen that this action is compatible with the right module structure. Moreover, $C_c(\a/ \gm) \times_{\alpha}^{alg} G / \gm$ acts on $X_0$ by bounded operators, relatively to the norm induced by the inner product $\langle \cdot \,,\, \cdot \rangle_{C^*_{\tau}(\a / \gm)}$, as we now show. For this we recall the conditional expectation $E_{\gm}$ of $C^*_{\tau}(\a / \gm) \times_{\alpha ,r} G / \gm$ onto $C_{\tau}^*(\a /\gm)$ as defined in Proposition \ref{Eggm is continuous}. For any $f, f_1 \in C_c(\a/ \gm) \times_{\alpha}^{alg} G / \gm$ we have that inside $C^*_{\tau}(\a / \gm) \times_{\alpha ,r} G / \gm$ the following holds: 
\begin{eqnarray*}
 \langle f*f_1 \,,\, f*f_1 \rangle_{C^*_{\tau}(\a / \gm)}  & = & ((f*f_1)^**(f*f_1))(\gm) \\
& = & E_{\gm}((f*f_1)^**(f*f_1)) \\
& = & E_{\gm}(f_1^**f^**f*f_1 )\\
& \leq & \|f\|_{\tau, r}^2E_{\gm}(f_1^**f_1)\\
& = & \|f\|_{\tau, r}^2 \langle f_1 \,,\, f_1 \rangle_{C^*_{\tau}(\a / \gm)}\,,
\end{eqnarray*}
where we used the positivity of $E_{\gm}$ in $C^*_{\tau}(\a / \gm) \times_{\alpha ,r} G / \gm$. Since the norm $\| \cdot \|_{\tau}$ is just the restriction of the norm $\| \cdot \|_{\tau, r}$ we get
\begin{eqnarray}
\label{inequality showing the map of reduced crossed prod in special alg}
 \|\langle f*f_1 \,,\, f*f_1 \rangle_{C^*_{\tau}(\a / \gm)}\|_{\tau} & \leq & \|f\|_{\tau, r}^2 \|\langle f_1 \,,\, f_1 \rangle_{C^*_{\tau}(\a / \gm)}\|_{\tau}\,,
\end{eqnarray}
which shows that $C_c(\a/ \gm) \times_{\alpha}^{alg} G / \gm$ acts on $X_0$ by bounded operators. Moreover, inequality (\ref{inequality showing the map of reduced crossed prod in special alg}) shows that this action extends to an action of $C^*_{\tau}(\a/ \gm) \times_{\alpha, r} G / \gm$ on $X$ and thus gives rise to a $^*$-homomorphism $\Phi: C^*_{\tau}(\a/ \gm) \times_{\alpha, r} G / \gm \to \mathcal{L}(X)$. We will now show that $\Phi$ is injective. Firstly, we will prove that $\Phi$ is injective on $C_{\tau}^*(\a / \gm)$, which is the same as to show that
\begin{align}
\label{Phi is injective on CtauAgm norms}
 \| \Phi(f) \|_{\mathcal{L}(X)} = \| f \|_{\tau}\,,
\end{align}
for all $f \in C_c(\a / \gm)$. It is clear from inequality (\ref{inequality showing the map of reduced crossed prod in special alg}) that $\| \Phi(f) \|_{\mathcal{L}(X)} \leq \| f \|_{\tau}$. The converse inequality follows from the fact that, for any $f, f_1 \in C_c(\a / \gm)$, we have
\begin{eqnarray*}
 \| \langle f* f_1 \,,\, f*f_1 \rangle_{C^*_{\tau}(\a / \gm)} \|_{\tau} & = & \| f \cdot f_1 \|_{\tau}^{\frac{1}{2}}\,.
\end{eqnarray*}
Before we prove that $\Phi$ is injective in the whole of $C^*_{\tau}(\a/ \gm) \times_{\alpha, r} G / \gm$ we need to establish some notation and results.

As usual, $Y:=C^*_{\tau}(\a / \gm)$ is a Hilbert module over itself. We define the map $j_{\gm}:Y \to X$ simply by inclusion, i.e. $j_{\gm}(f) := f$. It is then easy to see that $j_{\gm}$ is adjointable with adjoint $j_{\gm}^*:X \to Y$ given by $j_{\gm}(f) = f(\gm)$, for any $f \in X_0$. It is also easy to see that, for any $f \in C_c(\a / \gm)$ we have
\begin{eqnarray*}
\langle j_{\gm}(f) \,,\, j_{\gm}(f) \rangle_{C^*_{\tau}(\a / \gm)} & = & \langle f \,,\, f \rangle_{C^*_{\tau}(\a / \gm)}\,,
\end{eqnarray*}
where the inner product on the left (respectively, right) hand side corresponds to the inner product in $X$ (respectively, in $Y$). Thus, $j_{\gm}$ is an isometry between $Y$ and $X$ and has therefore norm $1$.

Let $\widehat{E}: \Phi(C_c(\a/ \gm) \times_{\alpha}^{alg} G / \gm) \to C_{\tau}^*(\a / \gm)$ be the map defined by
\begin{align*}
\widehat{E}(\Phi(f)) := \Phi(f(\gm))\,.
\end{align*}
First let us say a few words about why $\widehat{E}$ is well-defined. This is the case because $\Phi$ is injective on $C_c(\a/ \gm) \times_{\alpha}^{alg} G / \gm$, which is easily seen to be true because $C_c(\a/ \gm) \times_{\alpha}^{alg} G / \gm$ is an essential $^*$-algebra (\cite[Theorem 3.9]{palma crossed 11}).

We claim that $\widehat{E}$ is continuous with respect to the norm of $\mathcal{L}(X)$. First we notice that for any $f \in C_c(\a/ \gm)$ we have that (as elements of $\mathcal{L}(Y)$)
\begin{align*}
f = j_{\gm}^* \Phi(f) j_{\gm}\,.
\end{align*}
Let $f \in C_c(\a/ \gm) \times_{\alpha}^{alg} G / \gm$. We have
\begin{align*}
\|\widehat{E}(\Phi(f))\|_{\mathcal{L}(X)} = \| \Phi(f(\gm)) \|_{\mathcal{L}(X)}\,.
\end{align*}
As we proved in (\ref{Phi is injective on CtauAgm norms}), the norm $\| \cdot \|_{\mathcal{L}(X)}$ when restricted to $\Phi(C_c(\a / \gm))$ is such that $\| \Phi(g) \|_{\mathcal{L}(X)} = \| g \|_{\tau}$, and moreover the norm $\| \cdot \|_{\tau}$ coincides with the norm $\| \cdot \|_{\mathcal{L}(Y)}$, since $\mathcal{L}(Y) = M(C^*_{\tau}(\a / \gm))$. Hence we have:
\begin{eqnarray*}
\|\widehat{E}(\Phi(f))\|_{\mathcal{L}(X)} & = &\| \Phi(f(\gm)) \|_{\mathcal{L}(X)} \;\; = \;\; \| f(\gm) \|_{\mathcal{L}(Y)}\\
& = & \|j_{\gm}^* \Phi(f) j_{\gm} \|_{\mathcal{L}(Y)} \;\; \leq \;\; \|\Phi(f) \|_{\mathcal{L}(X)}\,,
\end{eqnarray*}
which shows that $\widehat{E}$ is continuous with respect to the norm of $\mathcal{L}(X)$.

We can now prove that $\Phi$ is injective. First we notice that for any $f \in C_c(\a/ \gm) \times_{\alpha}^{alg} G / \gm$ we have $\widehat{E}(\Phi(f)) = \Phi(E_{\gm}(f))$. By continuity, this equality then holds for any $f \in C^*_{\tau}(\a/ \gm) \times_{\alpha, r} G / \gm$.
 Suppose now that $f \in C^*_{\tau}(\a/ \gm) \times_{\alpha, r} G / \gm$ is such that $\Phi(f) = 0$. Then we have
\begin{eqnarray*}
0 & = & \widehat{E}(\Phi(f^**f)) \;\; = \;\; \Phi(E_{\gm}(f^**f))\,.
\end{eqnarray*}
Since $\Phi$ is faithful on $C^*_{\tau}(\a / \gm)$, it then follows that $E_{\gm}(f^**f) = 0$,
and since $E_{\gm}$ is faithful this implies that $f^** f =  0$, i.e. $ f = 0$. Thus, $\Phi$ is injective.

 We will first prove part $b)$ of the theorem and only afterwards prove part $a)$. For that we need to show that $F$ is an idempotent, positive, $C^*_{\tau}(\a / \gm)$-linear map. The fact that $F$ is idempotent is obvious. Now, let $f \in C_c(\a/ \gm) \times_{\alpha}^{alg} G / \gm$. We have that
\begin{eqnarray*}
F(f^**f) & = & (f^**f) (\gm)\\
& = & \sum_{[h] \in G / \gm} f^*(h\gm) \overline{\alpha}_h (f(h^{-1}\gm))\\
& = & \sum_{[h] \in G / \gm} \Delta(h^{-1})\, \overline{\alpha}_h( f(h^{-1}\gm))^* \overline{\alpha}_h (f(h^{-1}\gm))\,,
\end{eqnarray*}
 which by continuity means that $F$ is positive. Moreover, for $f_1 \in C_c(\a / \gm)$ we have that
\begin{eqnarray*}
F(f_1*f) & = & (f_1*f)(\gm) \;\;=\;\; f_1 \cdot f(\gm)\\
& = & f_1 \cdot F(f)\,,
\end{eqnarray*}
and similarly $F(f*f_1) = F(f)\cdot f_1$. By continuity of $F$, it follows that $F(f_1*f)= f_1 \cdot F(f)$ and $F(f*f_1) = F(f)\cdot f_1$ for any $f \in C_{\tau}^*(\a/ \gm) \times_{\alpha,  \omega} G / \gm$ and $f_1 \in C^*_{\tau}(\a / \gm)$. Hence we have shown that $F$ is a conditional expectation, and therefore $b)$ is proven.

Now, let $f, g \in C_c(\a/ \gm) \times_{\alpha}^{alg} G / \gm$. We have that inside $C_{\tau}^*(\a/ \gm) \times_{\alpha,  \omega} G / \gm$ the following holds:
\begin{eqnarray*}
 \langle f*g \,,\, f*g \rangle_{C^*_{\tau}(\a / \gm)}  & = & ((f*g)^**(f*g))(\gm) \\
& = & F((f*g)^**(f*g)) \\
& = & F(g^**f^**f*g )\\
& \leq & \|f\|_{\omega}^2F(g^**g)\\
& = & \|f\|_{\omega}^2 \langle g \,,\, g \rangle_{C^*_{\tau}(\a / \gm)}\,,
\end{eqnarray*}
where we have used to the positivity of $F$. Since the norm $\| \cdot \|_{\tau}$ is just the restriction of the norm $\| \cdot \|_{\omega}$ we get
\begin{eqnarray}
\label{inequality showing the map of reduced omega crossed prod in special alg}
 \|\langle f*g \,,\, f*g \rangle_{C^*_{\tau}(\a / \gm)}\|_{\tau} & \leq & \|f\|_{\omega}^2 \|\langle g \,,\, g \rangle_{C^*_{\tau}(\a / \gm)}\|_{\tau}\,,
\end{eqnarray}
which shows that the action of $C_c(\a/ \gm) \times_{\alpha}^{alg} G / \gm$ on $X_0$ extends to an action of $C^*_{\tau}(\a/ \gm) \times_{\alpha, \omega} G / \gm$ on $X$ and thus gives rise to a $^*$-homomorphism from $ C^*_{\tau}(\a/ \gm) \times_{\alpha,  \omega} G / \gm$ to $\mathcal{L}(X)$. As the injectivity of $\Phi$ shows, the closure of the image of $C_c(\a/ \gm) \times_{\alpha}^{alg} G / \gm$ in $\mathcal{L}(X)$ is isomorphic to $C_{\tau}^*(\a/ \gm) \times_{\alpha, r} G / \gm$. Hence, we conclude that there is a map $\Lambda: C^*_{\tau}(\a/ \gm) \times_{\alpha, \omega} G / \gm \to C_{\tau}^*(\a/ \gm) \times_{\alpha, r} G / \gm$ such that $\Lambda(f) = f$, for $f \in C_c(\a/ \gm) \times_{\alpha}^{alg} G / \gm$, and so part $a)$ is proven.

Let us now prove $c)$. The direction $(\Longleftarrow)$ is clear, because $F$ is then nothing but the conditional expectation $E_{\gm}$, which is faithful. Let us now prove the direction $(\Longrightarrow)$. For any $f \in C_c(\a/ \gm) \times_{\alpha}^{alg} G / \gm$ we have that $E_{\gm} \circ \Lambda(f^**f) = F(f^**f)$. By continuity this formula holds for any $f \in C_{\tau}^*(\a/ \gm) \times_{\alpha, \omega} G / \gm$. Let  $f \in C_{\tau}^*(\a/ \gm) \times_{\alpha, \omega} G / \gm$ be such that $\Lambda(f) = 0$. Then we necessarily have that $0 = E_{\gm} \circ \Lambda(f^**f) = F(f^**f)$, and since $F$ is faithful we have that $f^**f = 0$, i.e. $f = 0$. \qed\\

\subsection{Alternative definition of $C_r^*(\a / \gm) \times_{\alpha, r} G / \gm$}

The $C^*$-direct limit $D_{r}(\a)$ played a key role in the definition of the reduced crossed product $C_{r}^*(\a / \gm) \rtimes_{\alpha, r} G / \gm$. In this section we will see that instead of $D_r(\a)$ one can use the more natural $C^*$-algebra $C_r^*(\a)$ to define the reduced crossed product $C_{r}^*(\a / \gm) \rtimes_{\alpha, r} G / \gm$. The algebra $C_r^*(\a)$ has several advantages over $D_{r}(\a)$. For instance $C_r^*(\a)$ appears more naturally in the setup for defining crossed products (recall that we start with the bundle $\a$ and then we form the various bundles $\a / \gm^g$ from it). Also, $C^*_r(\a)$, being a cross sectional algebra of a Fell bundle, seems to be structurally simpler  than $D_{r}(\a)$, which is a direct limit of cross sectional algebras of Fell bundles.

The question one might ask at this point is: can one similarly use $C^*(\a)$ instead of $D_{\max}(\a)$ in order to define $C^*(\a / \gm) \rtimes_{\alpha, r} G / \gm$ ? As we shall also see in this section, this is not possible in general. At the core of this problem lies the fact that one has always an embedding
\begin{align*}
 C^*_r(\a / H) \to M(C^*_r(\a))\,,
\end{align*}
extending the natural embedding of $C_c(\a / H)$ into $M(C_c(\a))$, whereas the analogous map
\begin{align*}
 C^*(\a / H) \to M(C^*(\a))\,,
\end{align*}
is not always injective. This implies that the while the algebra $D_r(\a)$ embeds naturally in $M(C_r(\a))$, the analogous map from $D_{\max}(\a)$ to $M(C^*(\a))$ is not an embedding in general.

We start with the following general result:\\

\begin{prop}
\label{cor mapping of C* A H to multiplier of completion of Cc A}
 Let $\| \cdot \|_{\tau}$ be any $C^*$-norm on $C_c(\a)$ and $C^*_{\tau}(\a)$ its completion. There is a unique mapping $C^*(\a / H) \to M(C^*_{\tau}(\a))$ which extends the action of $C_c(\a / H)$ on $C_c(\a)$.\\
\end{prop}

{\bf \emph{Proof:}} As is known $C^*_{\tau}(\a)$ is naturally a Hilbert $C^*_{\tau}(\a)$-module, whose algebra of adjointable operators $\mathcal{L}(C^*_{\tau}(\a))$ is precisely the multiplier algebra $M(C^*_{\tau}(\a))$. In particular $X:=C_c(\a)$ is an inner product $C^*_{\tau}(\a)$-module. Moreover, $X$ is also a right $C_c(\a) - C^*_{\tau}(\a)$ bimodule and a right $C_c(\a / H) - C^*_{\tau}(\a)$ bimodule (in the sense of Definition \ref{def of right A B bimodule and right-pre-hilbert}), under the canonical actions of $C_c(\a)$ and $C_c(\a / H)$ on $X$. Since $C_c(\a)$ acts on $X$ by bounded operators, it then follows from Lemma \ref{embedding of Cc A H into multiplier of completion of Cc A} (taking $K = \{e\}$) that $C_c(\a / H)$ acts on $X$ by bounded operators. Thus, by completion, we obtain a right-Hilbert bimodule $_{C^*(\a / H)} C^*_{\tau}(\a)_{C^*_{\tau}(\a)}$. Hence obtain a unique map $C^*(\a / H) \to M(C^*_{\tau}(\a))$ which extends the action of $C_c(\a / H)$ on $C_c(\a)$. \qed\\

As shall see later in this section the map $C^*(\a / H) \to M(C^*_{\tau}(\a))$ is not an embedding in general, not even when $C^*_{\tau}(\a) = C^*(\a)$. Nevertheless for the reduced norms we have the following result:\\

\begin{thm}
\label{embedding of Ccr in M Ccr}
 There is a unique embedding of $C^*_r(\a / H)$ into $M(C^*_r(\a))$ which extends the action of $C_c(\a / H)$ on $C_c(\a)$.\\
\end{thm}

{\bf \emph{Proof:}} From Proposition \ref{cor mapping of C* A H to multiplier of completion of Cc A} we know that there exists a unique $^*$-homo-morphism from $C^*(\a / H)$ to $M(C^*_r(\a))$, which extends the action of $C_c(\a / H)$ on $C_c(\a)$. Thus, we have a right-Hilbert bimodule $_{C^*(\a / H)} C^*_r(\a)_{C^*_r(\a)}$. Taking the balanced tensor product of this right-Hilbert bimodule with $_{C^*_r(\a )} L^2(\a)_{C_0(\a^0)}$ we get a $C^*(\a / H) - C_0(\a^0)$ right-Hilbert bimodule
\begin{align*}
 _{C^*(\a / H)}\Big(C^*_r(\a) \otimes_{C^*_r(\a)} L^2(\a) \Big) \,_{C_0(\a^0)}\,.
\end{align*}
Since the action of $C^*_r(\a)$ on $L^2(\a)$ is faithful, the kernels of the maps from $C^*(\a / H)$ to $M(C^*_r(\a))$ and $\mathcal{L}\Big(C^*_r(\a) \otimes_{C^*_r(\a)} L^2(\a) \Big)$ are the same.

Now, $C^*_r(\a) \otimes_{C^*_r(\a)} L^2(\a)$ is isomorphic to $L^2(\a)$ as a Hilbert $C^*(\a / H) - C_0(\a^0)$ bimodule. Hence, it follows that the kernel of the map from $C^*(\a / H)$ to $M(C^*_r(\a))$ is the same as the kernel of the map from $C^*(\a / H)$ to $\mathcal{L}(L^2(\a))$. Now, the latter map has the same kernel as the canonical map $\Lambda: C^*(\a / H) \to C^*_r(\a / H)$, by Lemma \ref{lemma C correspondeces H and K} applied when $K$ is the trivial subgroup. Thus, this gives an embedding of $C^*_r(\a/H)$ into $M(C^*_r(\a))$. $\qed$\\

The next result is a generalization of \cite[Proposition 2.10]{full} (see \cite[Example 2.15]{palma crossed 11}). Its proof relies ultimately on Lemma \ref{lemma C correspondeces H and K}, whose proof, we recall, was essentially an adaptation of the proof \cite[Proposition 2.10]{full} itself.\\

\begin{cor}
\label{true analogue of Echt Quigg}
 Suppose  $\a$ is amenable. Then, the kernel of the canonical map $C^*(\a / H) \to M(C^*(\a))$ is the same as the kernel of the canonical map $\Lambda: C^*(\a/ H) \to C^*_r(\a / H)$.\\
\end{cor}

{\bf \emph{Proof:}} In the proof of Proposition \ref{embedding of Ccr in M Ccr} we established that the kernel of the canonnical map $\Lambda: C^*(\a / H) \to C^*_r(\a / H)$ is the same as the kernel of the map $C^*(\a / H) \to M(C^*_r(\a))$, which is the same as the map $C^*(\a / H) \to M(C^*(\a))$ by amenability of $\a$. \qed\\

We now give an example where the map $C^*(\a / H) \to M(C^*(\a))$ is not injective:\\

\begin{ex}
 Let $\b$ be a non-amenable Fell bundle over the group $G$, and let $\a:= \b \times G$ be the associated Fell bundle over the transformation groupoid $G  \times G$. Following \cite[Example 2.15]{palma crossed 11}, we have a right $G$-action on $\a$ which entails the action on the groupoid $G \times G$, given by $(s,t)g := (s, tg)$. Moreover, since the $G$-action is free, it is $H$-good and satisfies the $H$-intersection property, for any subgroup $H \subseteq G$. In this example we will consider $H$ to be the whole group $G$. In this case the orbit groupoid $(G \times G) / G$ can be naturally identified with the group $G$, and moreover, the Fell bundle $\a / G$ is naturally identified with $\b$.

It is known that the bundle $\a$ is always amenable (see \cite[Remark 2.11]{full}), and therefore by Corollary \ref{true analogue of Echt Quigg} we have that the kernel of the map $C^*(\a / G) \to M(C^*(\a))$ is the same as the kernel of the canonical map $C^*(\a / G) \to C^*_r(\a / G)$. As we pointed out above, the bundle $\a / G$ is just $\b$, which is non-amenable by assumption. Hence, the canonical map $C^*(\a / G) \to C^*_r(\a / G)$ has a non-trivial kernel, and therefore the map $C^*(\a / G) \to M(C^*(\a))$ is not injective.\\
\end{ex}

We will now see that $D_r(\a)$ is canonically embedded in $M(C_r^*(\a))$, being the $C^*$-algebra generated by all the images of $C^*_r(\a / H)$ inside $M(C^*_r(\a))$, as in Proposition \ref{embedding of Ccr in M Ccr}, with $H \in \mathcal{C}$.\\

\begin{prop}
\label{embedding of DrA into MCrA}
 Let $K \subseteq H$ be subgroups of $G$ such that $[H:K] < \infty$. Then, the following diagram of canonical embeddings commutes:
 \begin{align}
 \label{diagram of embeddings CraH CrAK MCrA}
\xymatrix{ C_r^*(\a / H) \ar[r] \ar[dr] &  C_r^*(\a / K) \ar[d]\\
  & M(C_r^*(\a))\,.}
\end{align}
As a consequence $D_r(\a)$ embedds in $M(C_r^*(\a))$, being $^*$-isomorphic to the subalgebra of $M(C^*_r(\a))$ generated by all the $C^*_r(\a / H)$, with $H \in \mathcal{C}$.\\
\end{prop}

{\bf \emph{Proof:}} We have already proven in Proposition \ref{D(A) embedds in MCcA} that
\begin{align}
[a]_{xH} b_y = \sum_{[h] \in \mathcal{S}_x \backslash H / K} [\alpha_{h^{-1}}(a)]_{xh K} b_y\,,
\end{align}
for any $x, y \in X$, $a \in \a_x$ and $b \in \a_y$. Hence, by linearity, density and continuity, we conclude that diagram (\ref{diagram of embeddings CraH CrAK MCrA}) commutes. By the universal property of $D_r(\a)$ we then have a $^*$-homomorphism from $D_r(\a)$ to $M(C_r^*(\a))$, whose image is generated by all the images of $C^*_r(\a / H)$ inside $M(C^*_r(\a))$, for any $H \in \mathcal{C}$. This $^*$-homomorphism from $D_r(\a)$ to $M(C^*_r(\a))$ is injective because all the maps in diagram (\ref{diagram of embeddings CraH CrAK MCrA}) are injective. \qed\\

We can now give an equivalent definition for the reduced crossed product $C^*_r(\a / \gm) \times_{r, \alpha} G / \gm$, using the algebra $C^*_r(\a)$ instead of $D_r^*(\a)$. This can be advantageous as we observed in the opening paragraph of this subsection. Also, this equivalence of definitions will make the connection between our definition of a reduced crossed product by a Hecke pair and that of Laca, Larsen and Neshveyev in \cite{phase} more clear, as we shall see in the next subsection.\\

\begin{thm}
\label{alternative def faithful pi give faithful rep of red crossed prod}
 Let $\pi: C^*_r(\a) \to B(\mathscr{H})$ be a nondegenerate $^*$-representation, and $\widetilde{\pi}$ its extension to $M(C^*_r(\a))$. We have that
\begin{itemize}
 \item[i)] If $\widetilde{\pi}_{\alpha}: C^*_r(\a / \gm) \to B(\mathcal{H} \otimes \ell^2(G / \gm))$ is faithful, then $\widetilde{\pi}_{\alpha} \times (1 \otimes \rho)$ is a faithful representation of $C^*_r(\a / \gm) \times_{r, \alpha} G / \gm$. Consequently,
\begin{align*}
 \|f\|_{r, r}:= \| [\widetilde{\pi}_{\alpha} \times (1 \otimes \rho)](f)\|\,,
\end{align*}
for all $f \in C^*_r(\a / \gm) \times_{r, \alpha} G / \gm$.
 \item[ii)] If $\pi$ is faithful, then $\widetilde{\pi}_{\alpha}$ is faithful.\\
\end{itemize}
\end{thm}

{\bf \emph{Proof:}} By Proposition \ref{embedding of DrA into MCrA} $D_r(\a)$ is canonically embedded in $M(C^*_r(\a))$, so that $\widetilde{\pi}$ restricts to a $^*$-representation of $D_r(\a)$. This restriction is nondegenerate, because the restriction to $C^*_r(\a / \gm)$ is already nondegenerate, as follows from the following argument. Let $\xi \in \mathscr{H}$ be such that $\widetilde{\pi}(C_r^*(\a / \gm))\xi = 0$. For any $x \in X$ and $a \in \a_x$ we have
\begin{eqnarray*}
 \|\pi(a_x) \xi\|^2 & = & \langle \pi((a^*a)_{\mathbf{s}(x)})\xi\,,\, \xi \rangle\\
& = & \langle \pi(a^*_{x^{-1}} \cdot [a]_{x\gm}) \xi\,,\, \xi \rangle\\
& = & \langle \pi(a^*_{x^{-1}}) \widetilde{\pi}([a]_{x\gm}) \xi\,,\, \xi \rangle\\
& = & 0\,.
\end{eqnarray*}
Thus, by nondegeneracy of $\pi$ we get that $\xi = 0$, and therefore $\widetilde{\pi}$ restricted to $C^*_r(\a / \gm)$, and hence also $D_r(\a)$, is nondegenerate. We are now in the conditions of Theorem \ref{norm beta r only needs one faithful pi}.

Claim $ii)$ also follows from Theorem \ref{norm beta r only needs one faithful pi}, given the fact that a faithfull nondegenerate $^*$-representation of $C^*_r(\a)$ extends faithfully to $M(C^*_r(\a))$. \qed\\

\subsection{Comparison with Laca-Larsen-Neshveyev construction}
\label{comparison with larsen laca neshveyev section}

In \cite{phase}, Laca, Larsen and Neshveyev, based on the work of Connes-Marcolli \cite{connes marcolli} and Tzanev \cite{tzanev talk}, introduced an algebra which can be thought of as a reduced crossed product of an abelian algebra by an action of a Hecke pair.

The construction introduced by Laca, Larsen and Neshveyev was one of the motivations behind our definition of a crossed product by a Hecke pair. However, the setup for their construction in \cite{phase} is slightly different from ours, being on one side more particular, as it only allows one to take a crossed product by an abelian algebra, but also more general, as the underlying space is not assumed in \cite{phase} to be discrete. We will show in this section that when both setups agree, our crossed product is canonically isomorphic to the crossed product of \cite{phase}.

We will first briefly recall the setup and construction presented in \cite[Section 1]{phase}. In order to make a coherent and more meanigful comparison between our construction and that of \cite{phase} we will have to make a few simple modifications in the latter. Essentially, we will consider right actions of $G$ instead of left ones, and make the appropriate changes in the construction of \cite{phase} according to this.

Let $G$ be a group acting on the right on a locally compact space $X$. Let $\gm \subseteq G$ be a Hecke subgroup and consider the (right) action of $\gm \times \gm$ on $X \times G$, given by:
\begin{align}
 (x,g)(\gamma_1, \gamma_2):= (x\gamma_1, \gamma_1^{-1}g\gamma_2)\,.
\end{align}
Define $X \times_{\gm} G / \gm$ to be the quotient space of $X \times G$ by the action of $\gm \times \gm$. We assume that the space $X \times_{\gm} G / \gm$ is Hausdorff.\\

\begin{rem}
 In \cite{phase} the original assumption was that the action of $\gm$ on $X$ was proper (hence implying that $X \times_{\gm} G / \gm$ is Hausdorff), but as it was observed in \cite[Remark 1.4]{phase}, requiring that $X \times_{\gm} G / \gm$ is Hausdorff was actually enough for the construction to make sense, and this is an important detail for us as the actions we consider are not proper in general.\\
\end{rem}

 Let $C_c(X \times_{\gm} G / \gm)$ be the space of compactly supported continuous functions on $X \times_{\gm} G / \gm$. We will view the elements of $C_c(X \times_{\gm} G / \gm)$ as $(\gm \times \gm)$-invariant functions on $X \times G$. One can define a convolution product and involution in $C_c(X \times_{\gm} G / \gm)$ according to the following formulas:
\begin{align}
 (f_1 * f_2)(x, g) &:= \sum_{[h] \in G / \gm} f_1(x, h)f_2(xh, h^{-1}g)\,,\\
f^*(x,g) &:= \overline{f(xg,g^{-1})}\,.
\end{align}
For each given $x \in X$ we can define a $^*$-representation $\pi_x: C_c(X \times_{\gm} G / \gm) \to B(\ell^2(G / \gm))$ by
\begin{align}
 \pi_x(f) \delta_{h\gm} := \sum_{g \in G / \gm} f(xg,g^{-1}h) \delta_{g\gm}\,.
\end{align}
The $C^*$-algebra $C^*_r(X \times_{\gm} G / \gm)$ is defined as the completion of $C_c(X \times_{\gm} G / \gm)$ in the norm
\begin{align}
 \|f\|:= \sup_{x \in X} \|\pi_x(f)\|\,.
\end{align}

The setup behind this construction differs slightly from our own, so we will compare both constructions under the following assumptions:
\begin{itemize}
 \item $(G, \gm)$ is a Hecke pair;
 \item $X$ is a set (seen as both a discrete space and a discrete groupoid);
 \item There is a right action of $G$ on $X$;
 \item The $G$-action satisfies the $\gm$-intersection property.
\end{itemize}

We notice that since $X$ and $G$ are discrete the space $X \times_{\gm} G / \gm$ is also discrete and therefore Hausdorff, so that the necessary assumptions for the construction of \cite{phase} are satisfied. Also, since $X$ is just a set, the action $G$ on $X$ is necessarily $\gm$-good. Thus, the assumptions for our construction (\cite[Standing Assumption 3.1]{palma crossed 11}) are satisfied with respect to the trivial Fell bundle $\a$ over $X$ in which every fiber $\a_x$ is just $\mathbb{C}$. Recall that in this case $C_c(\a) = C_c(X)$ and $C_c(\a / \gm) = C_c(X / \gm)$.\\

\begin{thm}
 Let $(G, \gm)$ be a Hecke pair and $X$ a set. Assume that there is a right $G$-action on $X$ which satisfies the $\gm$-intersection property. Then, the map $\Phi:C_c(X / \gm) \times_{\alpha}^{alg} G/\gm \to C_c(X \times_{\gm} G / \gm)$ given by
\begin{align*}
 \Phi(f) \,(x, g) := \Delta(g)^{\frac{1}{2}}\, f(g \gm)(x)\,,
\end{align*}
is a $^*$-isomorphism. This map extends to a $^*$-isomorphism between the reduced completions $\Phi:C_0(X / \gm) \times_{\alpha, r} G/\gm \to C_r^*(X \times_{\gm} G / \gm)$. Moreover, under the $^*$-isomorphism $\Phi$, the $^*$-representation $\pi_x$ is just $(\widetilde{\varphi_x})_{\alpha} \times \rho$, where $\varphi_x$ is the $^*$-representation of $C_0(X)$ given by evaluation at $x$, i.e. $\varphi_x(f)= f(x)$.\\
\end{thm}

{\bf \emph{Proof:}} Let us first check that $\Phi$ is well-defined, i.e. $\Phi(f)$ is a $(\gm \times \gm)$-invariant function in $G \times X$, with compact support (as a function on $X \times_{\gm} G / \gm$). To see this, let $\gamma_1, \gamma_2 \in \gm$. We have that
\begin{eqnarray*}
 \Phi(f) \,(x\gamma_1, \gamma_1^{-1}g\gamma_2) & = & \Delta(\gamma_1^{-1}g\gamma_2)^{\frac{1}{2}}\,f(\gamma_1^{-1}g \gamma_2 \gm)(x \gamma_1) \\
& = & \Delta(g)^{\frac{1}{2}}\, \overline{\alpha}_{\gamma_1^{-1}}(f(g\gm))\,(x \gamma_1 )\\
& = & \Delta(g)^{\frac{1}{2}}\, f(g\gm)\,(x)\\
& = & \Phi(f) \,(x, g)\,,
\end{eqnarray*}
so that $\Phi(f)$ is $\gm \times \gm$-invariant. It is easy to see that $\Phi(f)$ has compact support (as a function on $X \times_{\gm} G / \gm$). Thus, $\Phi$ is  well-defined.

Let us now prove that $\Phi$ is a $^*$-homomorphism. It is clear that $\Phi$ is linear, so that we only need to check that $\Phi$ preserves products and the involution. For $f_1, f_2 \in C_c(X / \gm) \times_{\alpha}^{alg} G/\gm$ we have that
\begin{eqnarray*}
  & & \Phi(f_1*f_2) \,(x, g) \;\; =\;\;\\
 & = & \Delta(g)^{\frac{1}{2}}\,(f_1* f_2)(g \gm)\,(x) \\
& = & \sum_{[h] \in G / \gm}\Delta(g)^{\frac{1}{2}}\, f_1(h \gm)\overline{\alpha}_{h}(f_2(h^{-1}g\gm))\,(x )\\
& = & \sum_{[h] \in G / \gm} \Big(\Delta(h)^{\frac{1}{2}}\,f_1(h \gm) \, (x) \Big) \Big(\Delta(h^{-1}g)^{\frac{1}{2}}\,\overline{\alpha}_{h}(f_2(h^{-1}g\gm))\,(x )\Big)\\
& = & \sum_{[h] \in G / \gm} \Big(\Phi(f_1)\,(x, h)\Big) \Big(\Delta(h^{-1}g)^{\frac{1}{2}}\,f_2(h^{-1}g\gm)\,(xh)\Big)\\
& = & \sum_{[h] \in G / \gm} \Big(\Phi(f_1)\,(x, h)\Big) \Big(\Phi(f_2)(xh,h^{-1}g)\Big)\\
& = & \Phi(f_1)* \Phi(f_2)\,(x, g)\,.
\end{eqnarray*}
Also for $f \in C_c(X / \gm) \times_{\alpha}^{alg} G/\gm$ we have
\begin{eqnarray*}
 \Phi(f^*)\,(x,g) & = & \Delta(g)^{\frac{1}{2}}\,f^*(g\gm)\,(x) \;\; = \;\; \Delta(g)^{\frac{1}{2}}\Delta(g^{-1})\,\overline{\overline{\alpha}_g(f(g^{-1}\gm))\,(x )}\\
& = & \Delta(g^{-1})^{\frac{1}{2}}\, \overline{f(g^{-1}\gm)\,(xg)} \;\; = \;\; \overline{\Phi(f)\,(xg,g^{-1})}\\
& = & (\Phi(f))^*\,(x,g)\,.
\end{eqnarray*}
Hence, $\Phi$ is a $^*$-homomorphism. Let us now prove that $\Phi$ is injective. Suppose $\Phi(f) = 0$. Then for every $g \in G$ and $x \in X$ we have
\begin{eqnarray*}
 0 & = & \Phi(f) \,(x,g) \;\; = \;\; \Delta(g)^{\frac{1}{2}}\,f(g \gm) \,(x)\,.
\end{eqnarray*}
Hence, we conclude that $f(g\gm) = 0$ for all $g \in G$, and therefore $f = 0$, i.e. $\Phi$ is injective.

Let us now prove the surjectivity of $\Phi$. The elements of $C_c(X \times_{\gm} G / \gm)$ are simply linear combinations of characteristic functions of elements of $X \times_{\gm} G / \gm$, so in order to prove that $\Phi$ is surjective we only need to check that each of these characteristic functions belongs the image of $\Phi$. Let $[(x, g)] \in X \times_{\gm} G / \gm$. We claim that $\Phi(\Delta(g)^{-\frac{1}{2}}\,1_{x\gm}* \gm g \gm * 1_{xg\gm}) = 1_{[(x,g)]}$. To see this, we recall \cite[Lemma 3.15]{palma crossed 11} and notice that
\begin{eqnarray*}
 \Phi(1_{x\gm}* \gm g \gm * 1_{xg\gm})\,(x,g) & = & \Delta(g)^{\frac{1}{2}}\,.
\end{eqnarray*}
It is not difficult to see that $\Phi(1_{x\gm}* \gm g \gm * 1_{xg\gm})\,(y,h) = 0$ if $(y,h)$ does not belong to the $\gm \times \gm$-orbit of $(x, g)$, so that $\Phi(\Delta(g)^{-\frac{1}{2}}\,1_{x\gm}* \gm g \gm * 1_{xg\gm}) = 1_{[(x,g)]}$. Hence, we can conclude that $\Phi$ is surjective and therefore establishes a $^*$-isomorphism between $C_c(X / \gm) \times_{\alpha}^{alg} G/\gm$ and $C_c(X \times_{\gm} G / \gm)$.

We will now see that under the $^*$-isomorphism $\Phi$, the $^*$-representation $\pi_x$ is just $(\widetilde{\varphi_x})_{\alpha} \times \rho$, in other words $\pi_x \circ \Phi = (\widetilde{\varphi_x})_{\alpha} \times \rho$. This follows from the following computation:
\begin{eqnarray*}
 \pi_x \circ \Phi (f)\,\delta_{h\gm} & = & \sum_{g\gm \in G / \gm} \Phi(f)(xg, g^{-1}h)\,\delta_{g \gm}\\
& = & \sum_{g\gm \in G / \gm} \Delta(g^{-1}h)^{\frac{1}{2}}\, f(g^{-1} h \gm)(xg)\,\delta_{g \gm}\\
& = & \sum_{g\gm \in G / \gm} \Delta(g^{-1}h)^{\frac{1}{2}}\, \overline{\alpha}_g(f(g^{-1} h \gm))(x)\,\delta_{g \gm}\\
& = & \sum_{g\gm \in G / \gm} \Delta(g^{-1}h)^{\frac{1}{2}}\, \widetilde{\varphi_x}\big(\overline{\alpha}_g(f(g^{-1} h \gm)) \big) \,\delta_{g \gm}\\
& = & [(\widetilde{\varphi_x})_{\alpha} \times \rho](f)\,\delta_{g \gm}\,.
\end{eqnarray*}

Let us now prove that the $^*$-isomorphism $\Phi$ extends to a $^*$-isomorphism between $C_0(X / \gm) \times_{\alpha, r} G/\gm$ and  $C_r^*(X \times_{\gm} G / \gm)$.  Let $\pi:C_c(X \times_{\gm} G / \gm) \longrightarrow B(\ell^2(X))$ be the direct sum $^*$-representation $\pi:= \bigoplus_{x \in X} \pi_x$ on the Hilbert space $\bigoplus_{x \in X} \mathbb{C} \cong \ell^2(X)$. We then have that
\begin{eqnarray*}
 \pi\circ \Phi(f) & = & \bigoplus_{x \in X} \pi_x ( \Phi(f)) \\
& = & \bigoplus_{x \in X} [(\widetilde{\varphi_x})_{\alpha} \times \rho](f)\\
& = & [(\bigoplus_{x \in X} \widetilde{\varphi_x})_{\alpha} \times \rho](f)\\
& = & [(\widetilde{\bigoplus_{x \in X} \varphi_x})_{\alpha} \times \rho](f)\,.
\end{eqnarray*}
Now the $^*$-representation $\bigoplus_{x \in X} \varphi_x$ of $C_0(X)$ is obviously injective. Hence, by Theorem \ref{alternative def faithful pi give faithful rep of red crossed prod} $ii)$, it follows that $\pi \circ \Phi$ extends to a faithful $^*$-representation of $C_0(X / \gm) \times_{\alpha, r} G/\gm$. This implies that $\Phi$ extends to an isomorphism between  $C_0(X / \gm) \times_{\alpha, r} G/\gm$ and $C_r^*(X \times_{\gm} G / \gm)$, because
\begin{eqnarray*}
 \|\Phi(f)\| & = & \sup_{x \in X} \|\pi_x(\Phi(f))\| \;\; = \;\;  \|\pi \circ \Phi(f)\|\\
& = & \|f \|_{r,r}\,. 
\end{eqnarray*}\qed\\

\section{Other completions}
\label{other completions chapter}

Just like there are several canonical $C^*$-completions of a Hecke algebra, one can also consider different $C^*$-completions of crossed products by Hecke pairs. Especially interesting for this work are full $C^*$-crossed products, but we will also take a look at $C^*$-completions arising from a $L^1$-norm.\\

\subsection{Full $C^*$-crossed products}

In this section we define and study full $C^*$-crossed products by Hecke pairs. Just like in the reduced case,  several full $C^*$-crossed products can be considered, such as $C^*_r(\a / \gm) \times_{\alpha} G / \gm$ and $C^*(\a / \gm) \times_{\alpha} G / \gm$ where each of these is thought of as the full $C^*$-crossed product of $C^*_r(\a / \gm)$, respectively $C^*(\a / \gm)$, by the Hecke pair $(G, \gm)$. As is the case for Hecke algebras, full crossed products by Hecke  pairs do not have to exist in general.\\

\begin{df}
 Let $\|\cdot\|_{\tau}$ be an $\overline{\alpha}$-permissible $C^*$-norm in $D(\a)$. We will denote by
 $\|\cdot\|_{\tau, u}: C_c(\a / \gm) \times_{\alpha}^{alg} G / \gm  \longrightarrow \mathbb{R}^+_0 \cup \{ \infty\}$\index{1 norm tau u@$\Vert \cdot \Vert_{\tau, u}$} the function defined by
\begin{align}
\| f \|_{\tau, u} := \sup_{\Phi \in R_{\tau}} \| \Phi(f) \|\,,
\end{align}
where the supremum is taken over the class $R_{\tau}$ of $^*$-representations of $C_c(\a / \gm) \times_{\alpha}^{alg} G / \gm$ whose restrictions to $C_c(\a / \gm)$ are continuous with respect to $\| \cdot \|_{\tau}$.\\
\end{df}

\begin{prop}
 We have that $\| \cdot \|_{\tau, u}$ is a $C^*$-norm in $C_c(\a / \gm) \times_{\alpha}^{alg} G / \gm$ if and only if $\| f\|_{\tau, u} < \infty$ for all $f \in C_c(\a / \gm) \times_{\alpha}^{alg} G / \gm$.\\
\end{prop}

{\bf \emph{Proof:}} $(\Longrightarrow):$ This direction is trivial since a norm must take values in $\mathbb{R}^+_0$.

$(\Longleftarrow):$ It is clear in this case that $\| \cdot \|_{\tau, u}$ defines a $C^*$-seminorm. To check that it is a true $C^*$-norm it is enough to find a faithful $^*$-representation $\Phi \in R_{\tau}$. This is easy because since $\| \cdot \|_{\tau}$ is $\overline{\alpha}$-permissible we can take any nondegenerate faithful $^*$-representation $\pi$ of $D_{\tau}(\a)$ and take $\Phi: = \pi_{\alpha} \times (1 \otimes \rho)$, which is a faithful $^*$-representation by Theorem \ref{norm beta r only needs one faithful pi}. We have that $\Phi \in R_{\tau}$ because its restriction to $C_c(\a / \gm)$ is just $\pi_{\alpha}$, which is continuous with respect to $\| \cdot \|_{\tau}$ by Lemma \ref{p alpha is continuous with repect to beta r norm}. \qed\\

\begin{df}
 Let $\|\cdot\|_{\tau}$ be an $\overline{\alpha}$-permissible $C^*$-norm in $D(\a)$. When $\| \cdot \|_{\tau, u}$ is a $C^*$-norm we will call it \emph{the universal norm associated to} $\| \cdot \|_{\tau}$. The completion of $C_c(\a / \gm) \times_{\alpha}^{alg} G / \gm$ with respect to this norm will be denoted by $C_{\tau}^*(\a / \gm) \times_{\alpha} G / \gm$\index{crossed product (full)@$C_{\tau}^*(\a / \gm) \times_{\alpha} G / \gm$} and referred to as the \emph{full crossed product} of $C^*_{\tau}(\a / \gm)$ by the Hecke pair $(G, \gm)$.\\
\end{df}

It is clear that $\| \cdot \|_{\tau, r} \leq \| \cdot \|_{\tau, u}$ so that the identity map on $C_c(\a / \gm) \times_{\alpha}^{alg} G / \gm$ extends to a surjective $^*$-homomorphism
\begin{align}
 C_{\tau}^*(\a / \gm) \times_{\alpha} G / \gm \longrightarrow C_{\tau}^*(\a / \gm) \times_{\alpha, r} G / \gm\,,
\end{align}
in case $\| \cdot \|_{\tau, u}$ is a norm.

In general, full crossed products do not necessarily exist, as it is already clear from the fact  that a Hecke algebra (which is a particular case of crossed product by a Hecke pair) does not need to have an enveloping $C^*$-algebra. Nevertheless, for Hecke pairs whose Hecke algebras are $BG^*$-algebras one can always assure the existence of full $C^*$-crossed products, as we show below. We recall that a $^*$-algebra is called a \emph{$BG^*$-algebra} if all of its pre-$^*$-representations are normed. Most Hecke algebras for which it is known that a full Hecke $C^*$-algebra exists are known to be $BG^*$-algebras, as we discussed in \cite{palma}.\\

\begin{thm}
\label{Hecke BG algebra implies existence of full crossed prod}
 If $\h(G, \gm)$ is a $BG^*$-algebra, then the full crossed product $C_{\tau}^*(\a / \gm) \times_{\alpha} G / \gm$ always exists, for any $\overline{\alpha}$-permissible norm $\| \cdot \|_{\tau}$.\\
\end{thm}

{\bf \emph{Proof:}} We will prove that when $\h(G, \gm)$ is a $BG^*$-algebra we have
\begin{align}
\label{sup full norm proof}
\sup_{\Phi} \| \Phi(f)\| < \infty\,,
\end{align}
where the supremum runs over the class of all $^*$-representations of $C_c(\a / \gm) \times_{\alpha}^{alg} G / \gm$. To see this we first notice that it is enough to consider nondegenerate $^*$-representations. Secondly, from \cite[Theorem 3.40]{palma crossed 11}, any nondegenerate $^*$-representation $\Phi$ of $C_c(\a / \gm) \times_{\alpha}^{alg} G / \gm$ is the integrated form of a covariant pre-$^*$-representation $(\Phi|, \omega_{\Phi})$, so that we can write $\Phi = \Phi| \times \omega_{\Phi}$. Taking any element $[a]_{x\gm}* \gm g \gm * 1_{\mathbf{s}(x) g\gm}$ of the canonical spanning set of elements of the crossed product we then have
\begin{eqnarray*}
 \| \Phi([a]_{x\gm}* \gm g \gm * 1_{\mathbf{s}(x) g\gm})\| & = & \| \Phi|([a]_{x\gm}) \omega_{\Phi}( \gm g \gm ) \widetilde{\Phi|}( 1_{\mathbf{s}(x) g\gm})\|\\
& \leq & \| \Phi|([a]_{x\gm})\omega_{\Phi}( \gm g \gm )\|\,.
\end{eqnarray*}
Now, since $\h(G, \gm)$ is a $BG^*$-algebra we have that $\omega_{\Phi}$ is normed, i.e. $\omega_{\Phi}( \gm g \gm )$ is a bounded operator. Moreover, because it is a $BG^*$-algebra, $\h(G, \gm)$ has an enveloping $C^*$-algebra. Hence, we conclude that
\begin{eqnarray*}
 & \leq & \| \Phi|([a]_{x\gm})\| \|\omega_{\Phi}( \gm g \gm )\|\\
& \leq & \| [a]_{x\gm}\|_{C^*(\a / \gm)} \|\gm g \gm \|_{C^*(G, \gm)}\,.
\end{eqnarray*}
Thus, it is clear that
\begin{align*}
 \sup_{\Phi} \| \Phi([a]_{x\gm}* \gm g \gm * 1_{\mathbf{s}(x) g\gm}) \| < \infty\,.
\end{align*}
Since this is true for the elements of the canonical spanning set, it follows that (\ref{sup full norm proof}) holds for any $f \in C_c(\a / \gm) \times_{\alpha}^{alg} G / \gm$. \qed\\

Any $BG^*$-algebra necessarily has an enveloping $C^*$-algebra. Is it then possible to weaken the assumptions on Theorem \ref{Hecke BG algebra implies existence of full crossed prod} to cover all Hecke algebras with an enveloping $C^*$-algebra? In other words:\\

\begin{open}
If $\h(G, \gm)$ has an enveloping $C^*$-algebra, do the full crossed products $C^*_{\tau}(\a / \gm) \times_{\alpha} G / \gm$ always exist?\\
\end{open}

We do not know the answer to this question. In fact we do not even have an example of a Hecke algebra which has an enveloping $C^*$-algebra and is not a $BG^*$-algebra. More generally even, the author does not know any example of a $^*$-algebra that can be faithfully represented on a Hilbert space and has an enveloping $C^*$-algebra, but is not a $BG^*$-algebra.

Regarding the existence of full crossed products we will show, in the next section, that they can exist for Hecke pairs for which the Hecke algebra does not have an enveloping $C^*$-algebra. Namely, the full crossed product $C_0(G / \gm) \times_{\alpha} G / \gm$, arising from the action of $G$ on itself by translation, exists for all Hecke pairs $(G, \gm)$.\\

\subsection{$L^1$-norm and associated $C^*$-completion}

We now define a $L^1$-norm on $C_c(\a / \gm) \times_{\alpha}^{alg} G / \gm$, whose corresponding enveloping $C^*$-algebra can still be understood as a crossed product of $C^*_{\tau}(\a / \gm)$ by the Hecke pair $(G, \gm)$, for a $\overline{\alpha}$-permissible norm $\| \cdot \|_{\tau}$.\\

\begin{df}
 Let $\| \cdot \|_{\tau}$ be an $\overline{\alpha}$-permissible $C^*$-norm on $D(\a)$. We define the norm $\| \cdot \|_{\tau, L^1}$\index{1 norm tau L1@$\Vert \cdot \Vert_{\tau, L^1}$} in $C_c(\a / \gm) \times_{\alpha}^{alg} G / \gm$ by:
\begin{align}
 \| f \|_{\tau, L^1} := \sum_{[g] \in \gm \backslash G / \gm}  L(g)\, \| f(g\gm) \|_{\tau} \,.
\end{align}\\

\end{df}

Before we prove that  $\| \cdot \|_{\tau, L^1}$ is a norm we observe that $\| \cdot \|_{\tau, L^1}$ is well-defined, i.e. it does not depend on the chosen representative $g$ of $[g]$, because for any $\gamma \in \gm$ we have, using the fact that the $\| \cdot \|_{\tau}$ is $\overline{\alpha}$-permissible,
\begin{eqnarray*}
 \| f(\gamma g \gm) \|_{\tau} & = & \| \overline{\alpha}_{\gamma}(f( g \gm)) \|_{\tau} \;\; = \;\; \| f( g \gm) \|_{\tau}\,.
\end{eqnarray*}
With this observation at hand we can easily derive another formula for  $\| \cdot \|_{\tau, L^1}$, for which we have
\begin{align}
\label{another formula for norm tau L1}
 \| f \|_{\tau, L^1} = \sum_{[g] \in G / \gm}   \| f(g\gm) \|_{\tau}\,.
\end{align}\\

\begin{prop}
 The function $\| \cdot \|_{\tau, L^1}$ is a norm for which
\begin{align*}
 \| f_1 * f_2 \|_{\tau, L^1} \leq \| f_1\|_{\tau, L^1} \|f_2 \|_{\tau, L^1} \qquad \text{and} \qquad \|f^* \|_{\tau, L^1} = \| f \|_{\tau, L^1}\,.
\end{align*}
Thus, under this norm $C_c(\a / \gm) \times_{\alpha}^{alg} G / \gm$ becomes a normed $^*$-algebra.\\
\end{prop}

{\bf \emph{Proof:}} It is easy to check that $\| \cdot \|_{\tau, L^1}$ is a vector space norm in $C_c(\a / \gm) \times_{\alpha}^{alg} G / \gm$. Let us prove first that $\|f^* \|_{\tau, L^1} = \| f \|_{\tau, L^1}$. We have
\begin{eqnarray*}
 \|f^* \|_{\tau, L^1} &  = & \sum_{[g] \in \gm \backslash G / \gm}  L(g)\, \| f^*(g\gm) \|_{\tau}\\
& = & \sum_{[g] \in \gm \backslash G / \gm}  L(g)\Delta(g^{-1})\, \|\overline{\alpha}_g(f(g^{-1}\gm))^* \|_{\tau}\\
& = & \sum_{[g] \in \gm \backslash G / \gm}  L(g^{-1})\, \| f(g^{-1}\gm) \|_{\tau}\,.
\end{eqnarray*}
Since $[g] \mapsto [g^{-1}]$ is a bijection of the set $\gm \backslash G / \gm$, we get
\begin{eqnarray*}
 & = & \sum_{[g] \in \gm \backslash G / \gm}  L(g)\, \| f(g\gm) \|_{\tau}\\
 & = & \| f \|_{\tau, L^1}\,.
\end{eqnarray*}

Let us now prove that $\| f_1 * f_2 \|_{\tau, L^1} \leq \| f_1\|_{\tau, L^1} \|f_2 \|_{\tau, L^1}$. For this we will use the formula for $\| \cdot \|_{\tau, L^1}$ given by (\ref{another formula for norm tau L1}). We have that
\begin{eqnarray*}
 \| f_1 * f_2 \|_{\tau, L^1} & = & \sum_{[g] \in  G / \gm}  \| (f_1* f_2)(g\gm) \|_{\tau}\\
 & \leq & \sum_{[g] \in  G / \gm} \sum_{[h] \in G / \gm} \| f_1(h\gm) \|_{\tau}\| \overline{\alpha}_h(f_2(h^{-1}g\gm)) \|_{\tau}\,.
\end{eqnarray*}
Using the fact that $\| \cdot \|_{\tau}$ is $\overline{\alpha}$-permissible we have
\begin{eqnarray*}
 & = &   \sum_{[g] \in  G / \gm} \sum_{[h] \in G / \gm}  \| f_1(h\gm) \|_{\tau}\| f_2(h^{-1}g\gm) \|_{\tau}\\
 & = &    \sum_{[h] \in G / \gm} \sum_{[g] \in  G / \gm}  \| f_1(h\gm) \|_{\tau}\| f_2(h^{-1}g\gm) \|_{\tau}\\
 & = &    \sum_{[h] \in G / \gm} \sum_{[g] \in  G / \gm}  \| f_1(h\gm) \|_{\tau}\| f_2(g\gm) \|_{\tau}\\
& = &    \Big(\sum_{[h] \in G / \gm} \| f_1(h\gm) \|_{\tau} \Big) \Big( \sum_{[g] \in  G / \gm}  \| f_2(g\gm) \|_{\tau} \Big)\\
& = & \| f_1\|_{\tau, L^1} \|f_2 \|_{\tau, L^1}\,.
\end{eqnarray*}\qed\\

 Completing $C_c(\a / \gm) \times_{\alpha}^{alg} G / \gm$ in the norm $\| \cdot \|_{\tau, L^1}$ we obtain a Banach $^*$-algebra, and taking the enveloping $C^*$-algebra of this Banach $^*$-algebra we obtain a $C^*$-completion of $C_c(\a / \gm) \times_{\alpha}^{alg} G / \gm$, which we denote by $C^*_{\tau}(\a / \gm) \times_{\alpha, L^1} G / \gm$\index{crossed product (L1)@$C^*_{\tau}(\a / \gm) \times_{\alpha, L^1} G / \gm$}.
We notice that the restriction of the norm $\| \cdot \|_{\tau, L^1}$ to $C_c(\a / \gm)$ is precisely the norm $\| \cdot \|_{\tau}$, from which we can conclude that $\| \cdot \|_{\tau, u}$ is always greater or equal to the $C^*$-norm of $C^*_{\tau}(\a / \gm) \times_{\alpha, L^1} G / \gm$. This means that, if $\| \cdot \|_{\tau, u}$ is a norm, there is canonical map
\begin{align}
\label{map crossed product full to L1}
 C^*_{\tau}(\a / \gm) \times_{\alpha} G / \gm \to C^*_{\tau}(\a / \gm) \times_{\alpha, L^1} G / \gm\,.
\end{align}
In case the crossed product is just the Hecke algebra itself, the map (\ref{map crossed product full to L1}) is just the usual map
\begin{align*}
 C^*(G, \gm) \to C^*(L^1(G, \gm))\,.
\end{align*}

So far we have seen three canonical $C^*$-crossed products of $C^*_{\tau}(\a / \gm)$ by the Hecke pair $(G, \gm)$, and these are $C^*_{\tau}(\a / \gm) \times_{\alpha, r} G / \gm$, $C^*_{\tau}(\a / \gm) \times_{\alpha, L^1} G / \gm$ and $C^*_{\tau}(\a / \gm) \times_{\alpha} G / \gm$ if it exists. Each one of these corresponds respectively, in the Hecke algebra case, to the completions $C^*_r(G, \gm)$, $C^*(L^1(G, \gm))$ and $C^*(G, \gm)$. It is an interesting problem, which we will not explore here, to understand how the Schlichting completion construction and the remaining Hecke $C^*$-algebra $pC^*(\overline{G})p$ carry over to the setting of crossed products by Hecke pairs.\\

\section{Stone-von Neumann theorem for Hecke pairs}
\label{Stone von Neumann chapter}

A modern version of the Stone-von Neumann theorem in the language of crossed products by groups states that (see \cite[Theorem C.34]{morita equiv})
\begin{align*}
 C_0(G) \times_{\alpha} G\;\; \cong\;\; C_0(G) \times_{\alpha, r} G \;\;\cong\;\; \mathcal{K}(\ell^2(G))\,.
\end{align*}
More precisely, if $\alpha$ is the action of $G$ on $C_0(G)$ by right translation, $M:C_0(G) \to B(\ell^2(G))$ the $^*$-representation by pointwise multiplication and $\rho$ the right regular representation of $G$ on $\ell^2(G)$, then $(M, \rho)$ is a covariant representation of the system $(C_0(G), G)$ and $M \times \rho$ is a faithful $^*$-representation of $C_0(G) \times_{\alpha} G$ with range $\mathcal{K}(\ell^2(G))$.

It follows from this result that any covariant representation of $(C_0(G), G)$ is unitarily equivalent to an amplification $(1 \otimes M, 1 \otimes \rho)$ of $(M, \rho)$, since the algebra of compact operators has a trivial representation theory (\cite[Remark C.35]{morita equiv}).

 The goal of this section is to show how the Stone-von Neumann theorem generalizes to the setting of Hecke pairs and their crossed products. In the process we recover an Huef, Kaliszewski and Raeburn's  notion of a \emph{covariant pair} \cite{cov} and their version of the Stone-von Neumann theorem for Hecke pairs, which did not make use of crossed products and which we will now recall.

In \cite[Definition 1.1]{cov}, an Huef, Kaliszewski and Raeburn introduced the notion of a \emph{covariant pair} $(\pi, \mu)$ consisting of a nondegenerate $^*$-representation $\pi:C_0(G / \gm) \to B(\mathscr{H})$ and a unital $^*$-representation $\mu: \h(G, \gm) \to B(\mathscr{H})$ satisfying
\begin{align}
\label{covariant pairs equation}
 \mu(\gm g \gm) \pi(1_{x\gm}) \mu(\gm s\gm) = \sum_{\substack{[u] \in \gm g^{-1} \gm / \gm \\ [v] \in \gm  s \gm / \gm}} \pi(1_{xu\gm}) \mu( \gm u^{-1} v \gm ) \pi(1_{xv\gm}) \,.
\end{align}

The basic example of a covariant pair, computed in \cite[Example 1.5]{cov}, is that of $(M, \rho)$ where $M:C_0(G / \gm) \to B(\ell^2(G / \gm))$ is the $^*$-representation by pointwise multiplication and $\rho$ is the right regular representation of $\h(G, \gm)$. \\

\begin{rem}
One should note that the definition of the right regular representation $\rho$ used in \cite{cov} differs from ours, since in \cite{cov} the factor $\Delta^{\frac{1}{2}}$ is absent. Nevertheless, $(M, \rho)$ is still a covariant pair with our definition of $\rho$. Moreover, the results of \cite{cov} remain valid for our $\rho$ as well, up to multiplication by some factor in some of them.\\
\end{rem}

It was proven in \cite[Theorem 1.6]{cov} that all covariant pairs are unitarily equivalent to an amplification $(1 \otimes M, 1 \otimes \rho)$ of $(M, \rho)$, which can be seen as an analogue for Hecke pairs of the Stone-von Neumann theorem. It should be noted that this result was proven without any crossed product construction behind it.

In the following we will prove a Stone-von Neumann theorem for Hecke pairs in the language of crossed products, stating that
\begin{align*}
 C_0(G / \gm) \times G/ \gm \;\cong \; C_0(G / \gm) \times_r G/ \gm \;\cong \; \mathcal{K}(\ell^2(G / \gm))\,.
\end{align*}
We will also show that the covariant pairs of \cite{cov} coincide with our notion of a covariant $^*$-representation and we will recover an Huef, Kaliszewski and Raeburn's version of the Stone-von Neumann theorem (\cite[Theorem 1.6]{cov}) as a consequence of the above isomorphisms.

The case under consideration now is that when the groupoid $X$ is the set $G$ and $\a$ is the Fell bundle over (the set) $G$ whose fibers are $\mathbb{C}$. In this case we have $C_c(\a) = C_c(G)$ and, naturally, $C^*_r(\a) = C^*(\a) = C_0(G)$. We consider the action $\alpha$ of $G$ on $\a$ induced by the right multiplication of $G$ on itself. Since this action is free, it is $\gm$-good and satisfies the $\gm$-intersection property. Moreover the induced action $\overline{\alpha}$ of $G$ on $C_c(G)$ is simply the action by right translation. In this setting the groupoid $X / \gm$ is then nothing but the orbit set $G / \gm$, and $C_c(\a / \gm)= C_c(G / \gm)$. Moreover, $C^*_r(\a / \gm) = C^*(\a/ \gm) = C_0(G /\gm)$.\\

\begin{prop}
\label{matrix units svnt}
 Let $T_{g\gm, h\gm} \in C_c(G / \gm) \times_{\alpha}^{alg} G / \gm$ be the element
 \begin{align*}
 T_{g\gm, h\gm} := 1_{g \gm} * \gm g^{-1} h \gm * 1_{h \gm}\,.
 \end{align*}
 Then $\{ T_{g\gm, h \gm} \}_{g \gm , h \gm \in G / \gm}$ is a set of matrix units that span $C_c(G / \gm) \times_{\alpha}^{alg} G / \gm$.\\
\end{prop}

{\bf \emph{Proof:}} It is clear that $T_{g\gm, h \gm}^* = T_{h\gm, g\gm}$. Let us now compute the product $T_{g\gm, h \gm} * T_{s \gm, t \gm}$. If $h \gm \neq s \gm$, then $T_{g \gm, h \gm} * T_{s\gm, t \gm} = 0$. In case $h \gm = s \gm$, we get
\begin{eqnarray*}
T_{g \gm, h \gm} * T_{h \gm, t \gm} & = & 1_{g \gm} * \gm g^{-1} h \gm * 1_{h \gm}* \gm h^{-1} t \gm * 1_{t \gm}\\
& = & 1_{g \gm} * \big( \sum_{\substack{[u] \in \gm h^{-1} g \gm / \gm \\ [v] \in \gm h^{-1} t \gm / \gm}} 1_{hu\gm} * \gm u^{-1} v \gm * 1_{hv\gm} \;\big) * 1_{t \gm}\,.
\end{eqnarray*}
Now for the product $1_{g \gm} 1_{hu\gm}$ to be non-zero, we must have $h u \gm = g \gm$, i.e. $u \gm = h^{-1}g \gm$. Similarly, for the product $1_{hv \gm} 1_{t\gm}$ to be non-zero we must have $hv \gm = t \gm$, i.e. $v \gm = h^{-1} t \gm$. Thus,
\begin{eqnarray*}
T_{g \gm, h \gm} * T_{h \gm, t \gm} & = &  1_{g \gm} * 1_{hh^{-1} g\gm} * \gm (h^{-1}g)^{-1} h^{-1} t \gm * 1_{hh^{-1}t\gm} * 1_{t \gm}\\
& = & 1_{g \gm} * \gm g^{-1} t \gm * 1_{t \gm}\\
& = & T_{g \gm, t \gm}\,.
\end{eqnarray*}
Hence, $\{T_{g \gm, h\gm} \}_{g \gm, h\gm}$ is a set of matrix units. The fact that this set spans $C_c(G / \gm) \times_{\alpha}^{alg} G / \gm$ follows readily from \cite[Theorem 3.14]{palma crossed 11}, noting that for $x \in G$ and $g \gm \in G / \gm$ we have
\begin{align*}
1_{x \gm} * \gm g \gm * 1_{x g \gm} = T_{x \gm, x g \gm}\,.
\end{align*}
This finishes the proof. \qed\\

\begin{thm}
\label{stone von neumann theorem}
The full crossed product $C_0(G / \gm) \times_{\alpha} G / \gm$ exists and moreover
\begin{align*}
C_0(G / \gm) \times_{\alpha} G / \gm \;\;\cong\;\;  C_0(G / \gm) \times_{\alpha, r} G / \gm \;\; \cong \;\; \mathcal{K}(\ell^2(G / \gm))\,.\\
 \end{align*}
Denoting by $M:C_0(G / \gm) \to B(\ell^2(G))$ the $^*$-representation by pointwise multiplication, we have that $(M,\rho)$ is a covariant $^*$-representation and $M \times \rho$ is a faithful $^*$-representation of $C_0(G / \gm) \times_{\alpha} G / \gm$ with range $\mathcal{K}(\ell^2(G / \gm))$.\\
\end{thm}

{\bf \emph{Proof:}} By Proposition \ref{matrix units svnt} we have that $\{T_{g \gm, h\gm} \}_{g \gm , h \gm}$ is a set of matrix units that span $C_c(G / \gm) \times_{\alpha}^{alg} G / \gm$. Hence, the enveloping $C^*$-algebra of $C_c(G / \gm) \times_{\alpha}^{alg} G / \gm$ must exist. As it is known, there exists only one $C^*$-algebra, up to isomorphism, generated by a set of matrix units indexed by $G / \gm$, and that is $\mathcal{K}(\ell^2(G / \gm))$. Hence, we necessarily have
\begin{align*}
C_0(G / \gm) \times_{\alpha} G / \gm \;\;\cong\;\;  C_0(G / \gm) \times_{\alpha, r} G / \gm \;\; \cong \;\; \mathcal{K}(\ell^2(G / \gm))\,.
 \end{align*}
It has been shown in \cite[Example 1.5]{cov} that $(M, \rho)$ is a covariant pair, so that equality (\ref{covariant pairs equation}) holds. Since the action of $G$ on itself is free, it follows readily from \cite[Proposition 3.46]{palma crossed 11} that this means that $(M, \rho)$ is a covariant $^*$-representation.

Let us denote by $\phi: C_0(G) \to \mathbb{C}$ the $^*$-representation given by evaluation at the identity element, i.e.
\begin{align*}
 \phi(f) := f(e)\,,
\end{align*}
and let $\widetilde{\phi}$ be its extension to $M(C_0(G)) \cong C_b(G)$. We claim that $\widetilde{\phi}_{\alpha}$ restricted to $C_c(G / \gm)$ is nothing but the representation by multiplication, i.e. $\widetilde{\phi}_{\alpha} = M$, and this follows from the following computation, where $f \in C_c(G / \gm)$:
\begin{eqnarray*}
 \widetilde{\phi}_{\alpha}(f) \delta_{h\gm} & = & \widetilde{\phi}(\alpha_{h}(f)) \delta_{h\gm}\;\; = \;\; \alpha_{h}(f) \,(e)\, \delta_{h \gm}\\
& = & f(h\gm) \delta_{h\gm}\;\; = \;\; M(f) \delta_{h \gm}\,.
\end{eqnarray*}
Since $M = \widetilde{\phi}_{\alpha}$ is faithful, it now follows from Theorem \ref{alternative def faithful pi give faithful rep of red crossed prod} that $M \times \rho$ is a faithful $^*$-representation of $C_0(G / \gm) \times_{\alpha, r} G / \gm \cong C_0(G / \gm) \times_{\alpha} G / \gm$ in $B(\ell^2(G / \gm))$, whose image must necessarily be $\mathcal{K}(\ell^2(G / \gm))$. \qed\\

As a corollary of our Stone-von-Neumann theorem we recover \cite[Theorem 1.6]{cov} and we show that the covariant pre-$^*$-representations of $C_c(G / \gm) \times_{\alpha}^{alg} G / \gm$ coincide with the covariant pairs of \cite{cov}.\\

\begin{cor}
 Let $(G, \gm)$ be a Hecke pair, $\pi: C_0(G / \gm) \to B(\mathscr{H})$ a nondegenerate $^*$-representation and $\mu: \h(G, \gm) \to L(\pi(C_c(G / \gm)) \mathscr{H})$ a unital pre-$^*$-representation. Then $(\pi, \mu)$ is a covariant pre-$^*$-representation if and only if it is unitarily equivalent to an amplification $(1 \otimes M, 1 \otimes \rho)$ of $(M, \rho)$. In particular we have
\begin{itemize}
 \item[i)] All covariant pre-$^*$-representation are covariant $^*$-representations, and these are the same as the covariant pairs of \cite{cov}.
 \item[ii)] A $^*$-representation $\pi$ is equivalent to an amplification of $M$ if and only if there exists a $^*$-representation $\mu$ of $\h(G, \gm)$ such that $(\pi, \mu)$ is a covariant $^*$-representation.\\
\end{itemize}

\end{cor}

{\bf \emph{Proof:}} Let $(\pi, \mu)$ be a covariant pre-$^*$-representation of $C_c(G / \gm) \times_{\alpha}^{alg} G / \gm$. Then its integrated form $\pi \times \mu$ extends to a nondegenerate $^*$-representation of $C_0(G / \gm) \times_{\alpha} G / \gm$.  By Theorem \ref{stone von neumann theorem} $M \times \rho$ is a $^*$-isomorphism between $C_0(G / \gm) \times_{\alpha} G / \gm$ and $\mathcal{K}(\ell^2(G / \gm))$, so that $(\pi \times \mu)\circ (M \times \rho)^{-1}$ is a nondegenerate $^*$-representation of $\mathcal{K}(\ell^2(G / \gm))$. Since the algebra of compact operators has a trivial representation theory (see for example \cite[Lemma B.34]{morita equiv}) there exists a Hilbert space $\mathscr{H}$ such that $(\pi \times \mu)\circ (M \times \rho)^{-1}$ is unitarily equivalent to the representation $1 \otimes \mathrm{id}$ in $\mathscr{H} \otimes \ell^2(G / \gm)$. Hence, $(\pi \times \mu)$ is unitarily equivalent to $1 \otimes (M \times \rho)$. Now given the fact that $(M, \rho)$ is a covariant $^*$-representation, it is not difficult to see that $(1 \otimes M, 1 \otimes \rho)$ is also a covariant $^*$-representation and moreover
\begin{align*}
 1 \otimes (M \times \rho) = (1 \otimes M) \times (1 \otimes \rho)\,.
\end{align*}
By Proposition \ref{unitary equivalence and bijection} it follows that $(\pi, \mu)$ is unitarily equivalent to $(1 \otimes M, 1 \otimes \rho)$.

The converse is easier: suppose now that $(\pi, \mu)$ is equivalent ot an amplification $(1 \otimes M, 1 \otimes \rho)$ of $(M, \rho)$. Since $(1 \otimes M, 1 \otimes \rho)$ is a covariant $^*$-representation, it follows that $(\pi, \mu)$ must also be a covariant $^*$-representation.

Let us now check $i)$. As we have just proven, every covariant pre-$^*$-represen-tation is unitarily equivalent to an amplification $(1 \otimes M, 1 \otimes \rho)$ of $(M, \rho)$. Since, $(1 \otimes M, 1 \otimes \rho)$ is a covariant $^*$-representation it follows that every covariant pre-$^*$-representation is actually a covariant $^*$-representation.

Let us now prove $ii)$. Suppose $\pi:C_0(G / \gm) \to B(\mathscr{H})$ is equivalent to an amplification of $M$, i.e. there exists a Hilbert space $\mathscr{H}_0$ and a unitary $U: \mathscr{H} \to \mathscr{H}_0 \otimes \ell^2(G / \gm)$ such that $\pi = U (1 \otimes M) U^*$. As $(U(1 \otimes M)U^*, U(1 \otimes \rho)U^*)$ is a covariant $^*$-representation, we conclude that there exists a $^*$-representation $\mu$ such that $(\pi, \mu)$ is a covariant $^*$-representation. The converse follows easily from what we proved above: if there exists a $^*$-representation $\mu$ of $\h(G, \gm)$ such that $(\pi, \mu)$ is a covariant $^*$-representation, then $(\pi, \mu)$ is unitarily equivalent to an amplification $(1 \otimes M, 1 \otimes \rho)$ of $(M, \rho)$, and therefore $\pi$ is unitarily equivalent to an amplification of $M$. \qed\\

\section{Towards Katayama duality}
\label{towards katayama duality chapter}

The theory of crossed products by Hecke pairs we have developed is intended for applications in non-abelian crossed product duality. We have already taken the first step in this direction, having established a Stone-von Neumann theorem for Hecke pairs which reflects the results of an Huef, Kaliszewski and Raeburn \cite{cov}. We believe that this theory of crossed products by Hecke pairs can be further applied and bring insight into the emerging theory of crossed products by coactions of homogeneous spaces (\cite{full}, \cite{echt kal rae}). The basic idea here is to obtain duality results for ``actions'' and ``coactions'' of homogeneous spaces (those coming from Hecke pairs).

In this section we will explain how our construction of a crossed product of a Hecke pair seems very suitable for obtaining a form of Katayama duality for homogeneous spaces arising from Hecke pairs, with respect to what we would call the \emph{Echterhoff-Quigg crossed product}.

Let $\delta$ be a coaction of a discrete group $G$ on a $C^*$-algebra $B$ and $B \times_{\delta} G$ the corresponding crossed product. We follow the conventions and notation of \cite{full} for coactions and their crossed products. As it is known, there is an action $\widehat{\delta}$ of $G$ on $B \times_{\delta} G$, called the \emph{dual action}, determined by
\begin{align*}
 \widehat{\delta}_s \big(j_B(a)j_G(f) \big) := j_B(a)j_G(\sigma_{s}(f))\,, \qquad\qquad \forall a \in B, f \in C_0(G), s \in G\,,
\end{align*}
where $\sigma$ denotes the action of right translation on $C_0(G)$, i.e. $\sigma_s(f) (t) := f(ts)$.

Katayama's duality theorem (the original version comes from \cite[Theorem 8]{kat}) is an analogue for coactions of the duality theorem of Imai and Takai. A general version of it states that we have a canonical isomorphism
\begin{align}
 (B \times_{\delta} G) \times_{\widehat{\delta}, \omega} G \cong B  \otimes \mathcal{K}(\ell^2(G))\,,
\end{align}
for some $C^*$-completion of the $^*$-algebraic crossed product $(B \times_{\delta} G) \times_{\widehat{\delta}}^{alg} G$. This $C^*$-completion $(B \times_{\delta} G) \times_{\widehat{\delta}, \omega} G$ lies in between the full and the reduced crossed products, and the coaction $\delta$ is called \emph{maximal} (respectively, \emph{normal}) if this $C^*$-crossed product is the full (respectively, the reduced) crossed product.

 We would like to extend this duality result for coactions of homogeneous spaces $G / \gm$. In spirit we should obtain an isomorphism of the type
\begin{align}
\label{katayama}
 (B \times_{\delta } G/\gm) \times_{\widehat{\delta }, \omega} G/\gm \; \cong \; B  \otimes \mathcal{K}(\ell^2(G/ \gm))\,.
\end{align}
Of course, the expression on the left hand side makes no sense unless $\gm$ is normal in $G$ (in which case, the above is just Katayama's result), and there are a few reasons for that. First, it does not make sense in general for a homogeneous space to coact on a $C^*$-algebra, which consequently makes it difficult to give meaning to $B \times_{\delta} G / \gm$. Secondly, it also does not make sense in general for a homogeneous space $G / \gm$ to act (namely, by $\widehat{\delta}$) on a $C^*$-algebra.

The second objection can be overcomed by simply using our definition of a crossed product by (an ``action'' of) a Hecke pair $(G, \gm)$. The first objection can be overcomed because, even though there is no definition of a coaction of a homogeneous space, it is possible to define $C^*$-algebras $B \times_{\delta} G / \gm$  which can be thought of as crossed products of $B$ by a coaction of $G / \gm$ (\cite{full}, \cite{echt kal rae}). In this way the iterated crossed product in expression (\ref{katayama}) may have a true meaning. This is the approach we suggest towards a generalization of Katayama's result.

It is our point of view that such a Katayama duality result can hold when $B \times_{\delta} G / \gm$ is a certain $C^*$-completion of the algebra $C_c(\b  \times G / \gm)$ defined by Echterhoff and Quigg in \cite{full}. The full completion $C^*(\b \times G / \gm)$ has already been dubbed the \emph{Echterhoff and Quigg's crossed product} by the restricted coaction of $G / \gm$ in \cite{cov} (in case we start with a maximal coaction of $G$ on $B$). Hence, we have to ensure that Echterhoff and Quigg's algebra $C_c(\b \times G / \gm)$ falls in our set up for defining crossed products by Hecke pairs, and that is what we explain now.

We recall briefly the construction of Echterhoff and Quigg, and the reader is advised to read our \cite[Example 2.15]{palma crossed 11} again. We start with a coaction $\delta$ of a discrete group $G$ on a $C^*$-algebra $B$, and we denote by $\b$ its associated Fell bundle. Following \cite[Section 3]{ind} we denote by $ \b \times G$ the corresponding Fell bundle over the groupoid $ G \times G$. Elements of $\b \times G$ have the form $(b_s, t)$, with $b_s \in \b_s$ and $s, t \in G$. Any such element lies in the fiber $(\b \times G)_{(s,t)}$ over $(s,t) \in G \times G$.

 We recall that the multiplication and inversion in $G \times G$ are given by
\begin{align*}
  (s,tr)(t,r)=(st,r)\qquad \text{and}\qquad (s,t)^{-1}=(s^{-1}, st)\,,
\end{align*}
and the corresponding multiplication and involution on $\b \times G$ are given by
\begin{align*}
  (b_s,tr)(c_t,r)=((bc)_{st},r)\qquad \text{and}\qquad (b_s,t)^{-1}=(b^*_{s^{-1}}, st)\,.
\end{align*}
An important property of $C_c(\b \times G / \gm)$ is that it embeds densely in the coaction crossed product $B \times_{\delta} G$, by identifying $(a_s, t)$ with $j_B(a)j_G(1_t)$. In this setting we have that $B \times_{\delta} G \cong C^*(\b \times G) \cong C^*_r(\b \times G)$, as stated in \cite[Corollary 3.4]{ind}.

The dual action $\widehat{\delta}$ of $G$ on $B \times_{\delta} G$ is determined by $\widehat{\delta}_g(j_B(a)j_G(1_t)) = j_B(a)j_G(1_{tg^{-1}})$, which on the generators of $C_c(\b \times G)$ means
\begin{align}
\label{dual action expression}
 \widehat{\delta}_g(a_s, t) := (a_s, tg^{-1})\,.
\end{align}

Now let $H \subseteq G$ be a subgroup. Following \cite{full}, one can define a Fell bundle $\b \times G / H$  over the groupoid $G \times G / H$. We recall from \cite{full} that the operations on the groupoid $G \times G / H$ are defined by
\begin{align*}
  (s,trH)(t,rH)=(st,rH)\qquad \text{and}\qquad (s,tH)^{-1}=(s^{-1}, stH)\,,
\end{align*}
and the corresponding operations on the Fell bundle $\b \times G / H$ are defined by
\begin{align*}
  (b_s,trH)(c_t,rH)=((bc)_{st},rH)\qquad \text{and}\qquad (b_s,tH)^{-1}=(b_{s^{-1}}, stH)\,.
\end{align*}
The \emph{Echterhoff and Quigg algebra} is defined as the algebra $C_c(\b \times G / H)$ of finitely supported sections of this Fell bundle.

Let us now consider the case of a Hecke pair $(G, \gm)$ to see that the conditions of our definition of crossed products by Hecke pairs are met, and see that it makes sense to define $C_c(\b \times G /\gm) \times_{\widehat{\delta}}^{alg} G / \gm$.

For this we take the bundle $\a := \b \times G$ over the groupoid $X:= G \times G$, as above. We observe that there is a natural $G$-action $\widehat{\delta}$ on $\a$ given by (\ref{dual action expression}), which of course gives precisely the dual action of $G$ on $C_c(\a)$. This action also entails the canonical right action of $G$ on the groupoid $X$, given by
\begin{align}
\label{action of G on GG}
 (s,t)g := (s,tg)\,.
\end{align}
Since this action is free, it is $H$-good and satisfies the $H$-intersection property for any subgroup $H \subseteq G$.  Moreover, the orbit groupoid $X / H$ is canonically identified with the groupoid $G \times G / H$, simply by $(s, t)H \mapsto (s, tH)$. This canonical identification is easily seen to be a groupoid isomorphism, so that $X / H$ and $G \times G / H$ are ``the same'' groupoid. Under this identification, the Fell bundle $\a / H$ is just the Fell bundle $\b \times G / H$, and therefore we can canonically identify $C_c(\a / H)$ with $C_c(\b \times G / H)$.

We can now conclude that all our conditions are met and therefore we can define the $^*$-algebraic crossed product $C_c(\b \times G/\gm) \times_{\widehat{\delta}}^{alg} G / \gm$. We expect that there is a $C^*$-completion of the Echterhoff and Quigg algebra $C_c(\b \times G / \gm)$, which we would like to call \emph{the} Echterhoff and Quigg's crossed product, for which a form of Katayama duality as in (\ref{katayama}) holds.\\

\vskip1cm

\begin{flushleft}
 \emph{Department of Mathematics, University of Oslo,}\\   \emph{P.O. Box 1053 Blindern, NO-0316 Oslo, Norway}\\
\vskip0.5cm
 \emph{E-mail:} \texttt{ruip@math.uio.no}
\end{flushleft}

\end{document}